# A simple recursive representation of the Faulhaber series

Dietmar Pfeifer
Institut für Mathematik
Fakultät V
Carl-von-Ossietzky Universität Oldenburg

January 2025

## Abstract

We present a simple elementary recursive representation of the so called Faulhaber series $\sum_{k=1}^{n} k^N$ for integer $n$ and $N$, without reference to Bernoulli numbers or polynomials.

## 1. Introduction

A well-known historical problem is the explicit evaluation of the so called Faulhaber series $\sum_{k=1}^{n} k^N$ for integer $n$ and $N$, see e.g. Knuth [2] who also coined the wording Faulhaber series or formula. In the modern literature on this topic, an explicit representation of this expression is given on the basis of Bernoulli numbers and Bernoulli polynomials, see [1], [2] and [3]. In this note we present a simple recursive representation of the Faulhaber series without reference to Bernoulli numbers or polynomials.

## 2. Main result

Denote $s(n,N) := \sum_{k=1}^{n} k^N$ for integer $n$ and $N$. Then there holds

$$s(n,N) = \sum_{k=0}^{n} k^N = \frac{(n+1)^{N+1} - \sum_{j=0}^{N-1} \binom{N+1}{j} s(n,j)}{N+1} \qquad (1)$$

Proof:

---

D. Pfeifer (✉)
Institut für Mathematik, Schwerpunkt Versicherungs- und Finanzmathematik, Carl von Ossietzky Universität Oldenburg, Oldenburg, Deutschland
E-Mail: dietmar.pfeifer@uni-oldenburg.de

$$(n+1)^{N+1} = \sum_{k=0}^{n}\left((k+1)^{N+1} - k^{N+1}\right) = \sum_{k=0}^{n}\left[\sum_{j=0}^{N+1}\binom{N+1}{j}k^j - k^{N+1}\right] = \sum_{k=0}^{n}\sum_{j=0}^{N}\binom{N+1}{j}k^j$$

$$= \sum_{k=0}^{n}(N+1)k^N + \sum_{k=0}^{n}\sum_{j=0}^{N-1}\binom{N+1}{j}k^j = (N+1)\sum_{k=0}^{n}k^N + \sum_{j=0}^{N-1}\binom{N+1}{j}\sum_{k=0}^{n}k^j \qquad (2)$$

$$= (N+1)\sum_{k=0}^{n}k^N + \sum_{j=0}^{N-1}\binom{N+1}{j}s(n,j),$$

which by rearrangement leads to

$$s(n,N) = \sum_{k=0}^{n}k^N = \frac{(n+1)^{N+1} - \sum_{j=0}^{N-1}\binom{N+1}{j}s(n,j)}{N+1} \qquad (3)$$

Obviously,

$$s(n,N) = \sum_{k=0}^{n}k^N = \frac{(n+1)^{N+1} - \binom{N+1}{N-1}s(n,N-1) - \sum_{j=0}^{N-2}\binom{N+1}{j}s(n,j)}{N+1} \qquad (4)$$

which implies that $s(n,N)$ is a polynomial in $n$ of degree $N+1$ with leading term $\dfrac{(n+1)^{N+1}}{N+1}$ or more precisely, $s(n,N) = \dfrac{n^{N+1}}{N+1} + \dfrac{n^N}{2} + \mathcal{O}\left(n^{N-1}\right)$, cf. [3], p.3.

An evaluation of [1] for consecutive values of $N$ can easily be performed e.g. with the computer algebra system MAPLE, see the Appendix attached.

Seemingly, $s(n,N)$ can be represented as $s(n,N) = n^2(n+1)^2 P(n,N-3)$ for uneven $N \geq 3$ and $s(n,N) = n(n+1)(2n+1)P(n,N-2)$ for even $N \geq 2$ where $P(n,K)$ is a polynomial of degree $K$. This can be shown empirically at least for $N = 2, \cdots, 100$.

### References

[1] J. H. Conway and R. K. Guy (1996): The Book of Numbers. Springer, N.Y.

[2] D.E. Knuth (1993): Johann Faulhaber and Sums of Powers. Mathematics of Computation (61), no. 203, $277-294$.

[3] H. Richter and B. Schiekel (2004): Potenzsummen, Bernoulli-Zahlen und Eulersche Summenformel. Universität Ulm, doi:10.18725/OPARU-1819.

⌐ Appendix MAPLE Worksheet
⌐ > **restart:**
⌐ > **s[n,0]:=n+1;**

$$s_{n,0} := n+1$$

⌐ > **for N from 1 to 100 do**
  **s[n,N]:=sort(factor(simplify(((n+1)^(N+1)-sum(binomial(N+1,j)*s[n,j],j=0..N-1))/(N+1)))) od;**

$$s_{n,1} := \frac{(n+1)\,n}{2}$$

$$s_{n,2} := \frac{(n+1)\,(2\,n+1)\,n}{6}$$

$$s_{n,3} := \frac{(n+1)^2\,n^2}{4}$$

$$s_{n,4} := \frac{(n+1)\,(2\,n+1)\,(3\,n^2+3\,n-1)\,n}{30}$$

$$s_{n,5} := \frac{(2\,n^2+2\,n-1)\,(n+1)^2\,n^2}{12}$$

$$s_{n,6} := \frac{(n+1)\,(2\,n+1)\,(3\,n^4+6\,n^3-3\,n+1)\,n}{42}$$

$$s_{n,7} := \frac{(3\,n^4+6\,n^3-n^2-4\,n+2)\,(n+1)^2\,n^2}{24}$$

$$s_{n,8} := \frac{(2\,n+1)\,(n+1)\,(5\,n^6+15\,n^5+5\,n^4-15\,n^3-n^2+9\,n-3)\,n}{90}$$

$$s_{n,9} := \frac{(n^2+n-1)\,(2\,n^4+4\,n^3-n^2-3\,n+3)\,(n+1)^2\,n^2}{20}$$

$$s_{n,10} := \frac{(2\,n+1)\,(n+1)\,(n^2+n-1)\,(3\,n^6+9\,n^5+2\,n^4-11\,n^3+3\,n^2+10\,n-5)\,n}{66}$$

$$s_{n,11} := \frac{(2\,n^8+8\,n^7+4\,n^6-16\,n^5-5\,n^4+26\,n^3-3\,n^2-20\,n+10)\,(n+1)^2\,n^2}{24}$$

$$s_{n,12} := (2\,n+1)\,(n+1)\,(105\,n^{10}+525\,n^9+525\,n^8-1050\,n^7-1190\,n^6+2310\,n^5+1420\,n^4$$
$$-3285\,n^3-287\,n^2+2073\,n-691)\,n\,/\,2730$$

$$s_{n,13} := (30\,n^{10}+150\,n^9+125\,n^8-400\,n^7-326\,n^6+1052\,n^5+367\,n^4-1786\,n^3+202\,n^2$$
$$+1382\,n-691)\,(n+1)^2\,n^2\,/\,420$$

$$s_{n,14} := (2\,n+1)\,(n+1)\,(3\,n^{12}+18\,n^{11}+24\,n^{10}-45\,n^9-81\,n^8+144\,n^7+182\,n^6-345\,n^5$$
$$-217\,n^4+498\,n^3+44\,n^2-315\,n+105)\,n\,/\,90$$

$$s_{n,15} := (3\,n^{12}+18\,n^{11}+21\,n^{10}-60\,n^9-83\,n^8+226\,n^7+203\,n^6-632\,n^5-226\,n^4+1084\,n^3$$
$$-122\,n^2-840\,n+420)\,(n+1)^2\,n^2\,/\,48$$

$$s_{n,16} := (2n+1)(n+1)(15n^{14} + 105n^{13} + 175n^{12} - 315n^{11} - 805n^{10} + 1365n^9 + 2775n^8$$
$$- 4845n^7 - 6275n^6 + 11835n^5 + 7485n^4 - 17145n^3 - 1519n^2 + 10851n - 3617)n / 510$$

$$s_{n,17} := (10n^{14} + 70n^{13} + 105n^{12} - 280n^{11} - 565n^{10} + 1410n^9 + 2165n^8 - 5740n^7 - 5271n^6$$
$$+ 16282n^5 + 5857n^4 - 27996n^3 + 3147n^2 + 21702n - 10851)(n+1)^2 n^2 / 180$$

$$s_{n,18} := (2n+1)(n+1)(105n^{16} + 840n^{15} + 1680n^{14} - 2940n^{13} - 9996n^{12} + 16464n^{11}$$
$$+ 48132n^{10} - 80430n^9 - 167958n^8 + 292152n^7 + 380576n^6 - 716940n^5 - 454036n^4$$
$$+ 1039524n^3 + 92162n^2 - 658005n + 219335)n / 3990$$

$$s_{n,19} := (42n^{16} + 336n^{15} + 616n^{14} - 1568n^{13} - 4263n^{12} + 10094n^{11} + 22835n^{10} - 55764n^9$$
$$- 87665n^8 + 231094n^7 + 213337n^6 - 657768n^5 - 236959n^4 + 1131686n^3 - 127173n^2$$
$$- 877340n + 438670)(n+1)^2 n^2 / 840$$

$$s_{n,20} := (2n+1)(n+1)(165n^{18} + 1485n^{17} + 3465n^{16} - 5940n^{15} - 25740n^{14} + 41580n^{13}$$
$$+ 163680n^{12} - 266310n^{11} - 801570n^{10} + 1335510n^9 + 2806470n^8 - 4877460n^7$$
$$- 6362660n^6 + 11982720n^5 + 7591150n^4 - 17378085n^3 - 1540967n^2 + 11000493n$$
$$- 3666831)n / 6930$$

$$s_{n,21} := (30n^{18} + 270n^{17} + 585n^{16} - 1440n^{15} - 5020n^{14} + 11480n^{13} + 35355n^{12} - 82190n^{11}$$
$$- 190745n^{10} + 463680n^9 + 733035n^8 - 1929750n^7 - 1783781n^6 + 5497312n^5$$
$$+ 1981107n^4 - 9459526n^3 + 1062932n^2 + 7333662n - 3666831)(n+1)^2 n^2 / 660$$

$$s_{n,22} := (2n+1)(n+1)(15n^{20} + 150n^{19} + 400n^{18} - 675n^{17} - 3615n^{16} + 5760n^{15}$$
$$+ 29220n^{14} - 46710n^{13} - 189702n^{12} + 307908n^{11} + 933064n^{10} - 1553550n^9 - 3269646n^8$$
$$+ 5681244n^7 + 7413782n^6 - 13961295n^5 - 8845327n^4 + 20248638n^3 + 1795584n^2$$
$$- 12817695n + 4272565)n / 690$$

$$s_{n,23} := (30n^{20} + 300n^{19} + 750n^{18} - 1800n^{17} - 7776n^{16} + 17352n^{15} + 69212n^{14}$$
$$- 155776n^{13} - 493131n^{12} + 1142038n^{11} + 2666455n^{10} - 6474948n^9 - 10250315n^8$$
$$+ 26975578n^7 + 24943119n^6 - 76861816n^5 - 27701758n^4 + 132265332n^3 - 14861886n^2$$
$$- 102541560n + 51270780)(n+1)^2 n^2 / 720$$

$$s_{n,24} := (2n+1)(n+1)(273n^{22} + 3003n^{21} + 9009n^{20} - 15015n^{19} - 97097n^{18} + 153153n^{17}$$
$$+ 969969n^{16} - 1531530n^{15} - 8030022n^{14} + 12810798n^{13} + 52402714n^{12} - 85009470n^{11}$$
$$- 258027882n^{10} + 429546558n^9 + 904376004n^8 - 1571337285n^7 - 2050706147n^6$$
$$+ 3861727863n^5 + 2446689429n^4 - 5600898075n^3 - 496674885n^2 + 3545461365n$$
$$- 1181820455)n / 13650$$

$$s_{n,25} := (42n^{22} + 462n^{21} + 1309n^{20} - 3080n^{19} - 16079n^{18} + 35238n^{17} + 175833n^{16}$$
$$- 386904n^{15} - 1589210n^{14} + 3565324n^{13} + 11359537n^{12} - 26284398n^{11} - 61459521n^{10}$$
$$+ 149203440n^9 + 236279941n^8 - 621763322n^7 - 574962926n^6 + 1771689174n^5$$
$$+ 638548653n^4 - 3048786480n^3 + 342572785n^2 + 2363640910n - 1181820455)(n+1)^2$$

$$n^2 / 1092$$

$$s_{n,26} := (2n+1)(n+1)(7n^{24} + 84n^{23} + 280n^{22} - 462n^{21} - 3542n^{20} + 5544n^{19} + 42790n^{18}$$
$$- 66957n^{17} - 438977n^{16} + 691944n^{15} + 3655360n^{14} - 5829012n^{13} - 23884796n^{12}$$
$$+ 38741700n^{11} + 117639298n^{10} - 195829797n^9 - 412342529n^8 + 716428692n^7$$
$$+ 935010264n^6 - 1760729742n^5 - 1115557926n^4 + 2553701760n^3 + 226457058n^2$$
$$- 1616536467n + 538845489)n / 378$$

$$s_{n,27} := (2n^{24} + 24n^{23} + 76n^{22} - 176n^{21} - 1089n^{20} + 2354n^{19} + 14321n^{18} - 30996n^{17}$$
$$- 159536n^{16} + 350068n^{15} + 1447750n^{14} - 3245568n^{13} - 10356931n^{12} + 23959430n^{11}$$
$$+ 56043471n^{10} - 136046372n^9 - 215462444n^8 + 566971260n^7 + 524305554n^6$$
$$- 1615582368n^5 - 582288225n^4 + 2780158818n^3 - 312388431n^2 - 2155381956n$$
$$+ 1077690978)(n+1)^2 n^2 / 56$$

$$s_{n,28} := (2n+1)(n+1)(15n^{26} + 195n^{25} + 715n^{24} - 1170n^{23} - 10478n^{22} + 16302n^{21}$$
$$+ 150436n^{20} - 233805n^{19} - 1870583n^{18} + 2922777n^{17} + 19312501n^{16} - 30430140n^{15}$$
$$- 161044508n^{14} + 256781832n^{13} + 1052627806n^{12} - 1707332625n^{11} - 5184837923n^{10}$$
$$+ 8630923197n^9 + 18173840941n^8 - 31576223010n^7 - 41210314958n^6 + 77603583942n^5$$
$$+ 49167912016n^4 - 112553659995n^3 - 9981035393n^2 + 71248383087n - 23749461029)n$$
$$/ 870$$

$$s_{n,29} := (2n^{26} + 26n^{25} + 91n^{24} - 208n^{23} - 1502n^{22} + 3212n^{21} + 23353n^{20} - 49918n^{19}$$
$$- 313712n^{18} + 677342n^{17} + 3511303n^{16} - 7699948n^{15} - 31896622n^{14} + 71493192n^{13}$$
$$+ 228229813n^{12} - 527952818n^{11} - 1235045022n^{10} + 2998042862n^9 + 4748228223n^8$$
$$- 12494499308n^7 - 11554321208n^6 + 35603141724n^5 + 12832102985n^4$$
$$- 61267347694n^3 + 6884212818n^2 + 47498922058n - 23749461029)(n+1)^2 n^2 / 60$$

$$s_{n,30} := (2n+1)(n+1)(231n^{28} + 3234n^{27} + 12936n^{26} - 21021n^{25} - 217217n^{24}$$
$$+ 336336n^{23} + 3653650n^{22} - 5648643n^{21} - 54097043n^{20} + 83969886n^{19} + 677256580n^{18}$$
$$- 1057869813n^{17} - 7003032113n^{16} + 11033483076n^{15} + 58417981930n^{14}$$
$$- 93143714433n^{13} - 381864885017n^{12} + 619369184742n^{11} + 1880950772008n^{10}$$
$$- 3131110750383n^9 - 6593111576555n^8 + 11455222740024n^7 + 14950298960254n^6$$
$$- 28153059810393n^5 - 17837160922265n^4 + 40832271288594n^3 + 3620925455812n^2$$
$$- 25847523828015n + 8615841276005)n / 14322$$

$$s_{n,31} := (231n^{28} + 3234n^{27} + 12397n^{26} - 28028n^{25} - 233233n^{24} + 494494n^{23} + 4228301n^{22}$$
$$- 8951096n^{21} - 67317019n^{20} + 143585134n^{19} + 909110951n^{18} - 1961807036n^{17}$$
$$- 10187430547n^{16} + 22336668130n^{15} + 92565988487n^{14} - 207468645104n^{13}$$
$$- 662371629682n^{12} + 1532211904468n^{11} + 3584398350146n^{10} - 8701008604760n^9$$
$$- 13780520650294n^8 + 36262049905348n^7 + 33533468142038n^6 - 103328986189424n^5$$
$$- 37241900394100n^4 + 177812786977624n^3 - 19979663280772n^2 - 137853460416080n$$

$$+ \, 68926730208040) \, (n+1)^2 \, n^2 \, / \, 7392$$

$$s_{n,32} := (2\,n+1)\,(n+1)\,(1785\,n^{30} + 26775\,n^{29} + 116025\,n^{28} - 187425\,n^{27} - 2211615\,n^{26}$$
$$+ \, 3411135\,n^{25} + 43060745\,n^{24} - 66296685\,n^{23} - 748192627\,n^{22} + 1155437283\,n^{21}$$
$$+ \, 11157717389\,n^{20} - 17314294725\,n^{19} - 139916611147\,n^{18} + 218532064083\,n^{17}$$
$$+ \, 1447330368449\,n^{16} - 2280261584715\,n^{15} - 12074416778197\,n^{14} + 19251755959653\,n^{13}$$
$$+ \, 78929233449419\,n^{12} - 128019728153955\,n^{11} - 388783145644957\,n^{10}$$
$$+ \, 647184582544413\,n^9 + 1362764431926059\,n^8 - 2367738939161295\,n^7$$
$$- \, 3090155271763537\,n^6 + 5819102377225953\,n^5 + 3686855873113679\,n^4$$
$$- \, 8439834998283495\,n^3 - 748427988281089\,n^2 + 5342559481563381\,n$$
$$- \, 1780853160521127)\,n \, / \, 117810$$

$$s_{n,33} := (210\,n^{30} + 3150\,n^{29} + 13125\,n^{28} - 29400\,n^{27} - 278957\,n^{26} + 587314\,n^{25} + 5828849\,n^{24}$$
$$- \, 12245012\,n^{23} - 108432253\,n^{22} + 229109518\,n^{21} + 1736342717\,n^{20} - 3701794952\,n^{19}$$
$$- \, 23479879093\,n^{18} + 50661553138\,n^{17} + 263190804617\,n^{16} - 577043162372\,n^{15}$$
$$- \, 2391580917883\,n^{14} + 5360204998138\,n^{13} + 17113585671557\,n^{12} - 39587376341252\,n^{11}$$
$$- \, 92609724746533\,n^{10} + 224806825834318\,n^9 + 356045984051297\,n^8$$
$$- \, 936898793936912\,n^7 - 866400981994453\,n^6 + 2669700757925818\,n^5$$
$$+ \, 962215367678237\,n^4 - 4594131493282292\,n^3 + 516212586120019\,n^2$$
$$+ \, 3561706321042254\,n - 1780853160521127)\,(n+1)^2\,n^2 \, / \, 7140$$

$$s_{n,34} := (2\,n+1)\,(n+1)\,(3\,n^{32} + 48\,n^{31} + 224\,n^{30} - 360\,n^{29} - 4808\,n^{28} + 7392\,n^{27}$$
$$+ \, 107256\,n^{26} - 164580\,n^{25} - 2160340\,n^{24} + 3322800\,n^{23} + 37818560\,n^{22} - 58389240\,n^{21}$$
$$- \, 564958600\,n^{20} + 876632520\,n^{19} + 7087388420\,n^{18} - 11069398890\,n^{17} - 73320554770\,n^{16}$$
$$+ \, 115515531600\,n^{15} + 611693296160\,n^{14} - 975297710040\,n^{13} - 3998596474936\,n^{12}$$
$$+ \, 6485543567424\,n^{11} + 19695979462232\,n^{10} - 32786740977060\,n^9 - 69038450348116\,n^8$$
$$+ \, 119951046010704\,n^7 + 156549095132192\,n^6 - 294799165703640\,n^5$$
$$- \, 186778301351088\,n^4 + 427567034878452\,n^3 + 37915805126206\,n^2 - 270657225128535\,n$$
$$+ \, 90219075042845)\,n \, / \, 210$$

$$s_{n,35} := (2\,n^{32} + 32\,n^{31} + 144\,n^{30} - 320\,n^{29} - 3431\,n^{28} + 7182\,n^{27} + 81819\,n^{26} - 170820\,n^{25}$$
$$- \, 1757535\,n^{24} + 3685890\,n^{23} + 32898915\,n^{22} - 69483720\,n^{21} - 527564655\,n^{20}$$
$$+ \, 1124613030\,n^{19} + 7136365395\,n^{18} - 15397343820\,n^{17} - 79998828429\,n^{16}$$
$$+ \, 175395000678\,n^{15} + 726950581973\,n^{14} - 1629296164624\,n^{13} - 5201902840599\,n^{12}$$
$$+ \, 12033101845822\,n^{11} + 28149978993355\,n^{10} - 68333059832532\,n^9 - 108225003493471\,n^8$$
$$+ \, 284783066819474\,n^7 + 263354323392579\,n^6 - 811491713604632\,n^5$$
$$- \, 292478403204671\,n^4 + 1396448520013974\,n^3 - 156909809749917\,n^2$$
$$- \, 1082628900514140\,n + 541314450257070)\,(n+1)^2\,n^2 \, / \, 72$$

$$s_{n,36} := (2\,n+1)\,(n+1)\,(25935\,n^{34} + 440895\,n^{33} + 2204475\,n^{32} - 3527160\,n^{31}$$

$$- 52907400\, n^{30} + 81124680\, n^{29} + 1340320800\, n^{28} - 2051043540\, n^{27} - 30970228380\, n^{26}$$
$$+ 47480864340\, n^{25} + 628653368980\, n^{24} - 966720485640\, n^{23} - 11024870583272\, n^{22}$$
$$+ 17020666117728\, n^{21} + 164766614966044\, n^{20} - 255660255507930\, n^{19}$$
$$- 2067199674245342\, n^{18} + 3228629639121978\, n^{17} + 21386129044051954\, n^{16}$$
$$- 33693508385638920\, n^{15} - 178419587790654392\, n^{14} + 284476135878801048\, n^{13}$$
$$+ 1166317852320026224\, n^{12} - 1891714846419439860\, n^{11} - 5744960434807234172\, n^{10}$$
$$+ 9563298075420571188\, n^{9} + 20137266455993253364\, n^{8} - 34987548721700165640\, n^{7}$$
$$- 45662538101010717712\, n^{6} + 85987581512366159388\, n^{5} + 54479850561631897574\, n^{4}$$
$$- 124713566598630926055\, n^{3} - 110593542398966646061\, n^{2} + 78945814659160432119\, n$$
$$- 26315271553053477373\,)\, n\, /\, 1919190$$

$s_{n,37} := (2730\, n^{34} + 46410\, n^{33} + 224315\, n^{32} - 495040\, n^{31} - 5951400\, n^{30} + 12397840\, n^{29}$
$$+ 160599985\, n^{28} - 333597810\, n^{27} - 3943622137\, n^{26} + 8220842084\, n^{25} + 85271873869\, n^{24}$$
$$- 178764589822\, n^{23} - 1598608021493\, n^{22} + 3375980632808\, n^{21} + 25644176564377\, n^{20}$$
$$- 54664333761562\, n^{19} - 346916736318383\, n^{18} + 748497806398328\, n^{17}$$
$$+ 3889015856644177\, n^{16} - 8526529519686682\, n^{15} - 35339687005512323\, n^{14}$$
$$+ 79205903530711328\, n^{13} + 252883445504814817\, n^{12} - 584972794540340962\, n^{11}$$
$$- 1368473221645955993\, n^{10} + 3321919237832252948\, n^{9} + 5261212529888281117\, n^{8}$$
$$- 13844344297608815182\, n^{7} - 12802615116896386549\, n^{6} + 39449574531401588280\, n^{5}$$
$$+ 14218442943931098889\, n^{4} - 67886460419263786058\, n^{3} + 7627958656578415656\, n^{2}$$
$$+ 52630543106106954746\, n - 26315271553053477373\,)\, (n + 1)^{2}\, n^{2}\, /\, 103740$$

$s_{n,38} := (2\, n + 1)\, (n + 1)\, (105\, n^{36} + 1890\, n^{35} + 10080\, n^{34} - 16065\, n^{33} - 268821\, n^{32}$
$$+ 411264\, n^{31} + 7674072\, n^{30} - 11716740\, n^{29} - 201595044\, n^{28} + 308250936\, n^{27}$$
$$+ 4695486768\, n^{26} - 7197355620\, n^{25} - 95488672580\, n^{24} + 146831686680\, n^{23}$$
$$+ 1675348594780\, n^{22} - 2586438735510\, n^{21} - 25040704929590\, n^{20} + 38854276762140\, n^{19}$$
$$+ 314174275401760\, n^{18} - 490688551483710\, n^{17} - 3250296167441750\, n^{16}$$
$$+ 5120788526904480\, n^{15} + 27116515195781080\, n^{14} - 43235167057123860\, n^{13}$$
$$- 177259043534387060\, n^{12} + 287506148830142520\, n^{11} + 873129273588152560\, n^{10}$$
$$- 1453446984797300100\, n^{9} - 3060497502628908284\, n^{8} + 5317469746342012476\, n^{7}$$
$$+ 6939873613441083578\, n^{6} - 13068545293332631605\, n^{5} - 8279944417725463861\, n^{4}$$
$$+ 18954189273254511594\, n^{3} + 1680820293702618872\, n^{2} - 11998325077181184105\, n$$
$$+ 3999441692393728035\,)\, n\, /\, 8190$$

$s_{n,39} := (42\, n^{36} + 756\, n^{35} + 3906\, n^{34} - 8568\, n^{33} - 114716\, n^{32} + 238000\, n^{31} + 3477096\, n^{30}$
$$- 7192192\, n^{29} - 96759271\, n^{28} + 200710734\, n^{27} + 2392439483\, n^{26} - 4985589700\, n^{25}$$
$$- 51813810155\, n^{24} + 108613210010\, n^{23} + 971726952295\, n^{22} - 2052067114600\, n^{21}$$
$$- 15589385213405\, n^{20} + 33230837541410\, n^{19} + 210898865097985\, n^{18}$$

$$- 455028567737380\, n^{17} - 2364235215196373\, n^{16} + 5183498998130126\, n^{15}$$
$$+ 21483948243306721\, n^{14} - 48151395484743568\, n^{13} - 153734693796001895\, n^{12}$$
$$+ 355620783076747358\, n^{11} + 831931956289561819\, n^{10} - 2019484695655870996\, n^{9}$$
$$- 3198433689018862699\, n^{8} + 8416352073693596394\, n^{7} + 7783056712255106951\, n^{6}$$
$$- 23982465498203810296\, n^{5} - 8643776817855225318\, n^{4} + 41270019133914260932\, n^{3}$$
$$- 4637242797382218326\, n^{2} - 31995533539149824280\, n + 15997766769574912140\,)$$
$$(n+1)^{2}\, n^{2}\, /\, 1680$$

$$s_{n,40} := (2\,n+1)\,(n+1)\,(1155\, n^{38} + 21945\, n^{37} + 124355\, n^{36} - 197505\, n^{35} - 3664815\, n^{34}$$
$$+ 5595975\, n^{33} + 117089115\, n^{32} - 178431660\, n^{31} - 3469504500\, n^{30} + 5293472580\, n^{29}$$
$$+ 91890741580\, n^{28} - 140482848660\, n^{27} - 2144362484220\, n^{26} + 3286785150660\, n^{25}$$
$$+ 43628109306320\, n^{24} - 67085556534810\, n^{23} - 765538561163926\, n^{22}$$
$$+ 1181850620013294\, n^{21} + 11442470621246522\, n^{20} - 17754631241876430\, n^{19}$$
$$- 143564343447033586\, n^{18} + 224223830791488594\, n^{17} + 1485249935008802402\, n^{16}$$
$$- 2339986817908947900\, n^{15} - 12391121942332776676\, n^{14} + 19756676322453638964\, n^{13}$$
$$+ 81000031774122870332\, n^{12} - 131378385822411124980\, n^{11}$$
$$- 398983875027004427556\, n^{10} + 664165005451712203824\, n^{9}$$
$$+ 1398520464324836204642\, n^{8} - 2429863199213110408875\, n^{7}$$
$$- 3171234501955388675261\, n^{6} + 5971783352539638217329\, n^{5}$$
$$+ 3783591297996299967787\, n^{4} - 8661278623264269060345\, n^{3}$$
$$- 768065184528864688599\, n^{2} + 5482737088425431563071\, n - 1827579029475143854357\,)$$
$$n\, /\, 94710$$

$$s_{n,41} := (330\, n^{38} + 6270\, n^{37} + 34485\, n^{36} - 75240\, n^{35} - 1115235\, n^{34} + 2305710\, n^{33}$$
$$+ 37720705\, n^{32} - 77747120\, n^{31} - 1180558500\, n^{30} + 2438864120\, n^{29} + 33088904585\, n^{28}$$
$$- 68616673290\, n^{27} - 819512112573\, n^{26} + 1707640898436\, n^{25} + 17755418511501\, n^{24}$$
$$- 37218477921438\, n^{23} - 333019774043847\, n^{22} + 703258026009132\, n^{21}$$
$$+ 5342742357367833\, n^{20} - 11388742740744798\, n^{19} - 72278923239091857\, n^{18}$$
$$+ 155946589218928512\, n^{17} + 810267850484745933\, n^{16} - 1776482290188420378\, n^{15}$$
$$- 7362954737902436727\, n^{14} + 16502391765993293832\, n^{13} + 52687785170893554393\, n^{12}$$
$$- 121877962107780402618\, n^{11} - 285118807704098189037\, n^{10}$$
$$+ 692115577515976780692\, n^{9} + 1096163687661046924253\, n^{8}$$
$$- 2884442952838070629198\, n^{7} - 2667400664320069089360\, n^{6}$$
$$+ 8219244281478208807918\, n^{5} + 2962385715105745907749\, n^{4}$$
$$- 14144015711689700623416\, n^{3} + 1589270767419418748637\, n^{2}$$
$$+ 10965474176850863126142\, n - 5482737088425431563071\,)\,(n+1)^{2}\, n^{2}\, /\, 13860$$

$$s_{n,42} := (2\,n+1)\,(n+1)\,(1155\, n^{40} + 23100\, n^{39} + 138600\, n^{38} - 219450\, n^{37} - 4491410\, n^{36}$$

$$+ 6846840\, n^{35} + 159627930\, n^{34} - 242865315\, n^{33} - 5298343743\, n^{32} + 8068948272\, n^{31}$$
$$+ 158312542976\, n^{30} - 241503288600\, n^{29} - 4201104934152\, n^{28} + 6422409045528\, n^{27}$$
$$+ 98082587983884\, n^{26} - 150335086498590\, n^{25} - 1995762640613350\, n^{24}$$
$$+ 3068811504169320\, n^{23} + 35020420648813040\, n^{22} - 54065036725304220\, n^{21}$$
$$- 523452080787043660\, n^{20} + 812210639543217600\, n^{19} + 6567565171810524980\, n^{18}$$
$$- 10257453077487396270\, n^{17} - 67945000122467694310\, n^{16} + 107046226722445239600\, n^{15}$$
$$+ 566850624785162732000\, n^{14} - 903799050538966717800\, n^{13}$$
$$- 3705469102296744536096\, n^{12} + 6010103178714600163044\, n^{11}$$
$$+ 18252121572755687072582\, n^{10} - 30383233948490830690395\, n^{9}$$
$$- 63977436571954515249535\, n^{8} + 111157771832177188219500\, n^{7}$$
$$+ 145072924856398721112280\, n^{6} - 273188273200686675778170\, n^{5}$$
$$- 173086113855535433583474\, n^{4} + 396223307383646488264296\, n^{3}$$
$$+ 35136305036438958874578\, n^{2} - 250816111246481682444015\, n$$
$$+ 8360537041549389414005 )\, n\, /\, 99330$$

$$s_{n,\,43} := ( 210\, n^{40} + 4200\, n^{39} + 24500\, n^{38} - 53200\, n^{37} - 868357\, n^{36} + 1789914\, n^{35}$$
$$+ 32583789\, n^{34} - 66957492\, n^{33} - 1139297194\, n^{32} + 2345551880\, n^{31} + 35922733084\, n^{30}$$
$$- 74191018048\, n^{29} - 1008591595879\, n^{28} + 2091374209806\, n^{27} + 24990016265627\, n^{26}$$
$$- 52071406741060\, n^{25} - 541482624603305\, n^{24} + 1135036655947670\, n^{23}$$
$$+ 10156257492634225\, n^{22} - 2144755164121612\, n^{21} - 162940962727066475\, n^{20}$$
$$+ 347329477095349070\, n^{19} + 2204338405166124535\, n^{18} - 4756006287427598140\, n^{17}$$
$$- 24711284691428338745\, n^{16} + 54178575670284275630\, n^{15} + 224553006598680562705\, n^{14}$$
$$- 503284588867645401040\, n^{13} - 1606855011719696384023\, n^{12}$$
$$+ 3716994612307038169086\, n^{11} + 8695461095493742931051\, n^{10}$$
$$- 21107916803294524031188\, n^{9} - 33430445294446640034566\, n^{8}$$
$$+ 87968807392187804100320\, n^{7} + 81349521965020828487336\, n^{6}$$
$$- 250667851322229461074992\, n^{5} - 90345880548927993031941\, n^{4}$$
$$+ 431359612420085447138874\, n^{3} - 48469065379054935273427\, n^{2}$$
$$- 334421481661975576592020\, n + 167210740830987788296010 )\, (n + 1)^{2}\, n^{2}\, /\, 9240$$

$$s_{n,\,44} := ( 2\, n + 1 ) ( n + 1 ) ( 2415\, n^{42} + 50715\, n^{41} + 321195\, n^{40} - 507150\, n^{39} - 11393970\, n^{38}$$
$$+ 17344530\, n^{37} + 448021140\, n^{36} - 680703975\, n^{35} - 16554986805\, n^{34} + 25172832195\, n^{33}$$
$$+ 554177509935\, n^{32} - 843852681000\, n^{31} - 16591426102440\, n^{30} + 25309065494160\, n^{29}$$
$$+ 440492101059780\, n^{28} - 673392684336750\, n^{27} - 10285289120291850\, n^{26}$$
$$+ 15764630022606150\, n^{25} + 209288622974823350\, n^{24} - 321815249473538100\, n^{23}$$
$$- 3672493473576199468\, n^{22} + 5669647835101068252\, n^{21} + 54893043908306011856\, n^{20}$$
$$- 8517438978000955191\, n^{19} - 6887235663616748651\, n^{18}$$

$$+ 1075672544432517073662\ n^{17} + 7125217014202969778726\ n^{16}$$
$$- 1122566179352071320920\ n^{15} - 59444165059652102000128\ n^{14}$$
$$+ 9477907848662385096602652\ n^{13} + 388583002939858526135486\ n^{12}$$
$$- 6302640436529070440004555\ n^{11} - 19140529909419997578003 3\ n^{10}$$
$$+ 3186211508239378485672327\ n^{9} + 670914903565159087466695 1\ n^{8}$$
$$- 1165682930759707555483659 0\ n^{7} - 1521342407706123209001745 8\ n^{6}$$
$$+ 2864855076930938591244448 2\ n^{5} + 1815109507545321894644295 6\ n^{4}$$
$$- 4155091799787502137588667 5\ n^{3} - 368465383566797147610226 5\ n^{2}$$
$$+ 2630243975243946790209673 5\ n - 8767479917479822634032245\ )\ n\ /\ 217350$$

$$s_{n,\,45} := (210\ n^{42} + 4410\ n^{41} + 27195\ n^{40} - 58800\ n^{39} - 1051890\ n^{38} + 2162580\ n^{37}$$
$$+ 43560825\ n^{36} - 89284230\ n^{35} - 1691522070\ n^{34} + 3472328370\ n^{33} + 59597894755\ n^{32}$$
$$- 122668117880\ n^{31} - 1882510643580\ n^{30} + 3887689405040\ n^{29} + 52877130114725\ n^{28}$$
$$- 109641949634490\ n^{27} - 1310277510348023\ n^{26} + 2730196970330536\ n^{25}$$
$$+ 2839172977662790 1\ n^{24} - 59513656523586338\ n^{23} - 532529345390854447\ n^{22}$$
$$+ 1124572347305295232\ n^{21} + 8543596138141849733\ n^{20} - 18211764623588994698\ n^{19}$$
$$- 11558163581415225798 7\ n^{18} + 24937503625189351067 2\ n^{17}$$
$$+ 12957043728445468279 73\ n^{16} - 28407837819409871666 18\ n^{15}$$
$$- 11774147751814058476547\ n^{14} + 26389079285569104119712\ n^{13}$$
$$+ 84253373715523410121973\ n^{12} - 194895826716615924363658\ n^{11}$$
$$- 455935307411257641028692\ n^{10} + 1106766441539131206421042\ n^{9}$$
$$+ 1752882358512889766615083\ n^{8} - 4612531158564910739651208\ n^{7}$$
$$- 4265457449635488281288214\ n^{6} + 13143446057835887302227636\ n^{5}$$
$$+ 4737169929467804561883417\ n^{4} - 22617785916771496425994470\ n^{3}$$
$$+ 2541413040905925578964990\ n^{2} + 17534959834959645268064490\ n$$
$$- 8767479917479822634032245\ )\ (n+1)^{2}\ n^{2}\ /\ 9660$$

$$s_{n,\,46} := (2\,n+1)\ (n+1)\ (105\ n^{44} + 2310\ n^{43} + 15400\ n^{42} - 24255\ n^{41} - 595595\ n^{40}$$
$$+ 905520\ n^{39} + 25783450\ n^{38} - 39127935\ n^{37} - 1054801055\ n^{36} + 1601765550\ n^{35}$$
$$+ 39313919540\ n^{34} - 59771762085\ n^{33} - 1318680777721\ n^{32} + 2007907047624\ n^{31}$$
$$+ 39498842229812\ n^{30} - 60252216868530\ n^{29} - 1048793049628354\ n^{28}$$
$$+ 1603315682876796\ n^{27} + 24489539754399408\ n^{26} - 37535967473037510\ n^{25}$$
$$- 498325108428490910\ n^{24} + 766255646379255120\ n^{23} + 8744377882765253380\ n^{22}$$
$$- 13499694647337507630\ n^{21} - 130702945652764510670\ n^{20}$$
$$+ 20280426580281551982 0\ n^{19} + 16398836943151717087 60\ n^{18}$$
$$- 256122767437416532305 0\ n^{17} - 169654822912284389525 70\ n^{16}$$
$$+ 2672883727402974109038 0\ n^{15} + 1415393991503882935638 90\ n^{14}$$

$$-\ 2256735173625973108910025\,n^{13}-925234710033345038538729\,n^{12}$$

$$+\ 15006888237313162132536 06\,n^{11}+4557451692744892172701968\,n^{10}$$

$$-\ 7586521950982996365679755\,n^{9}-1597480465570993970677 9111\,n^{8}$$

$$+\ 2775546795905640774300854 4\,n^{7}+3622389016641992888317 5202\,n^{6}$$

$$-\ 6821356922915809719626707 5\,n^{5}-4321862528011248632030163 5\,n^{4}$$

$$+\ 9893472253474777807858599 0\,n^{3}+8773336966657925403261940\,n^{2}$$

$$-\ 6262736671736077144185905\,n+2087578890578692571472863 5\,)\ n\ /\ 9870$$

$$s_{n,47} := (\ 210\,n^{44}+4620\,n^{43}+30030\,n^{42}-64680\,n^{41}-1262730\,n^{40}+2590140\,n^{39}$$

$$+\ 57440010\,n^{38}-117470160\,n^{37}-2463942648\,n^{36}+5045355456\,n^{35}+96430075536\,n^{34}$$

$$-\ 197905506528\,n^{33}-3403769625356\,n^{32}+7005444757240\,n^{31}+107561432804676\,n^{30}$$

$$-\ 222128310366592\,n^{29}-3021570673520221\,n^{28}+6265269657407034\,n^{27}$$

$$+\ 74875404971187513\,n^{26}-156016079599782060\,n^{25}-1622446524951671865\,n^{24}$$

$$+\ 3400909129503125790\,n^{23}+30431454597465016365\,n^{22}-64263818324433158520\,n^{21}$$

$$-\ 488225012678151281025\,n^{20}+1040713843680735720570\,n^{19}$$

$$+\ 6604929486158917220925\,n^{18}-14250572815998570162420\,n^{17}$$

$$-\ 74043217241195862093357\,n^{16}+162337007298390294349134\,n^{15}$$

$$+\ 672835410195448263918689\,n^{14}-1508007827689286822186512\,n^{13}$$

$$-\ 4814671473829155768105594\,n^{12}+11137350775347598358397700\,n^{11}$$

$$+\ 26054490425935289154027994\,n^{10}-63246331627218176666453688\,n^{9}$$

$$-\ 100168721057994841172702694\,n^{8}+263583773743207859011859076\,n^{7}$$

$$+\ 243750195432215001334249782\,n^{6}-751084164607637861680358640\,n^{5}$$

$$-\ 270706274704615290050908260\,n^{4}+1292496714016868441782175160\,n^{3}$$

$$-\ 145229423269548003737600340\,n^{2}-1002037867477772434306974480\,n$$

$$+\ 501018933738886217153487240\,)\ (n+1)^{2}\,n^{2}\ /\ 10080$$

$$s_{n,48} := (2\,n+1)\,(n+1)\,(\ 3315\,n^{46}+76245\,n^{45}+533715\,n^{44}-838695\,n^{43}-22421113\,n^{42}$$

$$+\ 34051017\,n^{41}+1063856651\,n^{40}-1612810485\,n^{39}-47945788891\,n^{38}+72725088579\,n^{37}$$

$$+\ 1978679330277\,n^{36}-3004381539705\,n^{35}-73898631149759\,n^{34}+112350137494491\,n^{33}$$

$$+\ 2479949123866123\,n^{32}-3776098754546430\,n^{31}-74291542909828322\,n^{30}$$

$$+\ 113325363742015698\,n^{29}+1972683951301707774\,n^{28}-3015688608823569510\,n^{27}$$

$$-\ 4606291721571782697 0\,n^{26}+7060222012798852521 0\,n^{25}+9373123568392245541 10\,n^{24}$$

$$-\ 144126964532283109377 0\,n^{23}-164475294753677225141 98\,n^{22}$$

$$+\ 253919290357129993181 82\,n^{21}+245842620158526128116 106\,n^{20}$$

$$-\ 38145989475564569183325 0\,n^{19}-30845005977015499977410 30\,n^{18}$$

$$+\ 481748084393014784252817 0\,n^{17}+319108243772041857703416 60\,n^{16}$$

$$-\ 502749769877713525767765 75\,n^{15}-26622519992200385627948273 7\,n^{14}$$

$$+ 4244752883768914607076123393\, n^{13} + 1740298441323453326041803999\, n^{12}$$

$$- 2822685306173625719416512195\, n^{11} - 8572231447647794452315267229\, n^{10}$$

$$+ 1426968982455850454538181156941\, n^{9} + 3004743266054994512582250394 3\, n^{8}$$

$$- 5220599390310416995782433438 5\, n^{7} - 6813447324953176155304698915 1\, n^{6}$$

$$+ 1283047068258497273084826506 19\, n^{5} + 8129105555757417148403002825 7\, n^{4}$$

$$- 1860889367492861208802863678 45\, n^{3} - 1650200158287407398072566304 3\, n^{2}$$

$$+ 1177974707489541714112316784 87\, n - 3926582358298472380374389282 9\,) \, n \,/\, 324870$$

$$s_{n,49} := (1326\, n^{46} + 30498\, n^{45} + 208403\, n^{44} - 447304\, n^{43} - 9493055\, n^{42} + 19433414\, n^{41}$$

$$+ 472318327\, n^{40} - 964070068\, n^{39} - 22274214521\, n^{38} + 45512499110\, n^{37}$$

$$+ 963146250651\, n^{36} - 1971805000412\, n^{35} - 37764727720447\, n^{34} + 77501260441306\, n^{33}$$

$$+ 1333609437133935\, n^{32} - 2744720134709176\, n^{31} - 42147531802841998\, n^{30}$$

$$+ 87039783740393172\, n^{29} + 1184021855234441379\, n^{28} - 2455083494209275930\, n^{27}$$

$$- 29340593635684107751\, n^{26} + 61136270765577491432\, n^{25} + 635771072406920859487\, n^{24}$$

$$- 1332678415579419210406\, n^{23} - 1192485867581143388313 5\, n^{22}$$

$$+ 2518239576720228697667 6\, n^{21} + 1913156912867499981266 83\, n^{20}$$

$$- 4078137783407022832300 42\, n^{19} - 2588205525292395860241 919\, n^{18}$$

$$+ 5584224828925494003713 880\, n^{17} + 2901455121169004038807 7559\, n^{16}$$

$$- 6361332725230557477986 8998\, n^{15} - 2636570669460331044499 80418\, n^{14}$$

$$+ 5909274611443717836798 29834\, n^{13} + 1886675612067025273294 672275\, n^{12}$$

$$- 4364278685278422330269 174384\, n^{11} - 1020970422186233256764 0799967\, n^{10}$$

$$+ 2478368712900308746555 0774318\, n^{9} + 3925208275312916548143 3594231\, n^{8}$$

$$- 1032878526352614182841 7962780\, n^{7} - 9551587303043241647834 2146185\, n^{6}$$

$$+ 2943195986961262513851 02255150\, n^{5} + 1060788735671371178837 13911035\, n^{4}$$

$$- 5064773458304004871525 30077220\, n^{3} + 5690955500027662455754 5574465\, n^{2}$$

$$+ 3926582358298472380374 38928290\, n - 1963291179149236190187 19464145\,) \,(n+1)^{2}\, n^{2}$$

$$/\, 66300$$

$$s_{n,50} := (2\,n+1)\,(n+1)\,(429\, n^{48} + 10296\, n^{47} + 75504\, n^{46} - 118404\, n^{45} - 3433716\, n^{44}$$

$$+ 5209776\, n^{43} + 177855964\, n^{42} - 269388834\, n^{41} - 8790556346\, n^{40} + 13320528936\, n^{39}$$

$$+ 399690120544\, n^{38} - 606195445284\, n^{37} - 16528902129916\, n^{36} + 25096450917516\, n^{35}$$

$$+ 617619215226134\, n^{34} - 938977048297959\, n^{33} - 20729072134582867\, n^{32}$$

$$+ 31563096726023280\, n^{31} + 620996836434449376\, n^{30} - 947276803014685704\, n^{29}$$

$$- 1648962159067344132 0\, n^{28} + 2520807078751750483 2\, n^{27} + 3850359595475766105 512\, n^{26}$$

$$- 5901634286074079106 84\, n^{25} - 7834990698171041513 900\, n^{24}$$

$$+ 1204756776156026622 6192\, n^{23} + 1374848556426715926 27616\, n^{22}$$

$$- 2122510673447875220 54520\, n^{21} - 2054997867949212429 383888\, n^{20}$$

$$+ 3188622335596212405103092\, n^{19} + 2578333320106135602760 1406\, n^{18}$$

$$- 4026931096939014024395 3655\, n^{17} - 266742505880274126845610963\, n^{16}$$

$$+ 42024841430510626039039 3272\, n^{15} + 2225375820438958588738849 79936\, n^{14}$$

$$- 3548187937811031961277943540\, n^{13} - 14547150581890312152150257764\, n^{12}$$

$$+ 23594819841740984208864358416\, n^{11} + 71655262530103266642836795708\, n^{10}$$

$$- 119280303716025392068687372770\, n^{9} - 251166419012752676999004054970\, n^{8}$$

$$+ 436389780377141711532849768840\, n^{7} + 569535901810418645522117212000\, n^{6}$$

$$- 1072498742904198824049600702420\, n^{5} - 679511742411895276079901986300\, n^{4}$$

$$+ 1555516985069942326144653330660\, n^{3} + 1379401924596908381481940879 30\, n^{2}$$

$$- 984668781224507420294617797225\, n + 328222927074835806764872599075\,) \, n \,/\, 43758$$

$s_{n,51} := (\, 66\, n^{48} + 1584\, n^{47} + 11352\, n^{46} - 24288\, n^{45} - 558371\, n^{44} + 1141030\, n^{43}$

$$+ 30268271\, n^{42} - 61677572\, n^{41} - 1562497057\, n^{40} + 3186671686\, n^{39} + 74289274785\, n^{38}$$

$$- 151765221256\, n^{37} - 3218428879787\, n^{36} + 6588622980830\, n^{35} + 126251591452327\, n^{34}$$

$$- 259091805885484\, n^{33} - 4458896446218412\, n^{32} + 9176884698322308\, n^{31}$$

$$+ 140923231862795646\, n^{30} - 291023348423913600\, n^{29} - 3958886129259529743\, n^{28}$$

$$+ 8208795606942973086\, n^{27} + 9810313247849508 8131\, n^{26} - 204415060563933149348\, n^{25}$$

$$- 2125763309448991740055\, n^{24} + 4455941679461916629458\, n^{23}$$

$$+ 39871945161579457314971\, n^{22} - 84199832002620831259400\, n^{21}$$

$$- 639682961538023766846791\, n^{20} + 1363565755078668364952982\, n^{19}$$

$$+ 8653921549737439189873227\, n^{18} - 18671408854553546744699436\, n^{17}$$

$$- 97013026134161809379570652\, n^{16} + 212697461122877165503840740\, n^{15}$$

$$+ 881563521202411129283570222\, n^{14} - 1975824503527699424070981184\, n^{13}$$

$$- 6308286803340695856328069907\, n^{12} + 14592398110209091136727120998\, n^{11}$$

$$+ 34137157441076563228621373311\, n^{10} - 82866712992362217593969867620\, n^{9}$$

$$- 131243227004178215894426115585\, n^{8} + 345353167000718649382822098790\, n^{7}$$

$$+ 319366783298544816586293788225\, n^{6} - 984086733597808282555409675240\, n^{5}$$

$$- 354685221965912440868718871675\, n^{4} + 1693457177529633164292847418590\, n^{3}$$

$$- 190282734615144968616678511145\, n^{2} - 1312891708299343227059490396300\, n$$

$$+ 656445854149671613529745198150\,) \, (n+1)^{2}\, n^{2} \,/\, 3432$$

$s_{n,52} := (2\, n + 1)\, (n + 1)\, (\, 165\, n^{50} + 4125\, n^{49} + 31625\, n^{48} - 49500\, n^{47} - 1552100\, n^{46}$

$$+ 2352900\, n^{45} + 87436800\, n^{44} - 132331650\, n^{43} - 4719995830\, n^{42} + 7146159570\, n^{41}$$

$$+ 235379503810\, n^{40} - 356642335500\, n^{39} - 10724586768060\, n^{38} + 16265201319840\, n^{37}$$

$$+ 443730622962770\, n^{36} - 673728535104075\, n^{35} - 16582482557910665\, n^{34}$$

$$+ 25210588104418035\, n^{33} + 556572470200953855\, n^{32} - 847463999353639800\, n^{31}$$

$$- 16673795425082270920\, n^{30} + 2543442513730022 6280\, n^{29} + 442747952816706904240\, n^{28}$$

$$- 67683914179371046950\, n^{27} - 103383554188223770912220\, n^{26}$$

$$+ 158459526991304208715800\, n^{25} + 21037038349274128299686000\, n^{24}$$

$$- 32347855158867713493108000\, n^{23} - 3691483970386290590970096\, n^{22}$$

$$+ 56989652313737744539206840\, n^{21} + 5517692262264920340608450222\, n^{20}$$

$$- 8561487195542493831822787500\, n^{19} - 6922854267940571303854740010\, n^{18}$$

$$+ 10812355761687981647373249390\, n^{17} + 71620666014651725678655779670\, n^{16}$$

$$- 112837176902821579341670294200\, n^{15} - 59751593725201294895164963796\, n^{14}$$

$$+ 9526924943294302130983096040400\, n^{13} + 39059264662688022068693935899200\, n^{12}$$

$$- 6335235946567918416853241536900\, n^{11} - 192395194364263362786326315304600\, n^{10}$$

$$+ 32026897127923463626375569889140\, n^{9} + 67438468993777078770771143620802000\, n^{8}$$

$$- 117171152054627349969344493925660\, n^{7} - 1529210370003312978698181842143600\, n^{6}$$

$$+ 28796713152781062178939974594984\, n^{5} + 1824496752061106902233657770434200\, n^{4}$$

$$- 41765807857307134622974853854005\, n^{3} - 3703709846545067520662731869431\, n^{2}$$

$$+ 264384686984711685924815247311490\, n - 8812822899490389530827174910383\,) \, n \, /$$

$$17490$$

$$s_{n,53} := (\ 110\, n^{50} + 2750\, n^{49} + 20625\, n^{48} - 44000\, n^{47} - 1092212\, n^{46} + 2228424\, n^{45}$$

$$+ 64277939\, n^{44} - 130784302\, n^{43} - 3617750565\, n^{42} + 7366285432\, n^{41} + 188307789451\, n^{40}$$

$$- 383981864334\, n^{39} - 8970510220273\, n^{38} + 18325002304880\, n^{37} + 388810024013163\, n^{36}$$

$$- 795945050331206\, n^{35} - 15253855759326386\, n^{34} + 31303656568983978\, n^{33}$$

$$+ 538743653699984005\, n^{32} - 1108790963968951988\, n^{31} - 17027084867069172174\, n^{30}$$

$$+ 35162960698107296336\, n^{29} + 478334192882235201227\, n^{28}$$

$$- 991831346462577698790\, n^{27} - 11853359686082587015667\, n^{26}$$

$$+ 24698550718627751730124\, n^{25} + 256846431604288730959359\, n^{24}$$

$$- 538391413927205213648842\, n^{23} - 4817548135166250744710071\, n^{22}$$

$$+ 10173487684259706703068984\, n^{21} + 77290020899186163007933603\, n^{20}$$

$$- 164753529482632032718936190\, n^{19} - 1045614497302201495295795718\, n^{18}$$

$$+ 2255982524087035023310527626\, n^{17} + 11721648503306545698276631541\, n^{16}$$

$$- 25699279530700126419863790708\, n^{15} - 106515363361690183883514488010\, n^{14}$$

$$+ 238730006254080494186892766728\, n^{13} + 762201979659340932218961419379\, n^{12}$$

$$- 1763133965572762358624815605486\, n^{11} - 4124639508750796649377907998197\, n^{10}$$

$$+ 10012412983074355675380631601880\, n^{9} + 15857530032836817155922187377387\, n^{8}$$

$$- 41727473048747989969225006356654\, n^{7} - 38587654946078999255314129662369\, n^{6}$$

$$+ 118902782940905988479853265681392\, n^{5} + 42855023363214460583257935037035\, n^{4}$$

$$- 204612829667334909646369135755462\, n^{3} + 22991008738253949045739993684284\, n^{2}$$

$$+ 158630812190827011554889148386894\, n - 79315406095413505777444574193447\,)$$

$$(n+1)^2\, n^2\, /\, 5940$$

$$s_{n,54} := (2n+1)(n+1)(399\,n^{52} + 10374\,n^{51} + 82992\,n^{50} - 129675\,n^{49} - 4383015\,n^{48}$$
$$+ 6639360\,n^{47} + 267634380\,n^{46} - 404771250\,n^{45} - 15720292770\,n^{44} + 23782824780\,n^{43}$$
$$+ 856292812440\,n^{42} - 1296330631050\,n^{41} - 42791979501930\,n^{40} + 64836134568420\,n^{39}$$
$$+ 1950718834760010\,n^{38} - 2958496319424225\,n^{37} - 80721171329965185\,n^{36}$$
$$+ 122561005154659890\,n^{35} + 3016691339019095520\,n^{34} - 4586317511105973225\,n^{33}$$
$$- 101252631308301013149\,n^{32} + 154172105718004506336\,n^{31}$$
$$+ 3033331471245944778408\,n^{30} - 4627083259727919420780\,n^{29}$$
$$- 80545663205129115477324\,n^{28} + 123132036437557632926376\,n^{27}$$
$$+ 1880776141926041855145168\,n^{26} - 2882730231107841599180940\,n^{25}$$
$$- 38271039539766677120992100\,n^{24} + 58847924425203936481078620\,n^{23}$$
$$+ 671562829229741252680382890\,n^{22} - 1036768206057213847261113645\,n^{21}$$
$$- 10037906988397386370504348445\,n^{20} + 15575244585624686479387079490\,n^{19}$$
$$+ 125942077630456598837389520200\,n^{18} - 196700738738497241495777820045\,n^{17}$$
$$- 1302938806865213292937165127825\,n^{16} + 2052758579667068560153636601760\,n^{15}$$
$$+ 10870140501388441539972943939660\,n^{14} - 17331590041916196590036234210370\,n^{13}$$
$$- 71057467808971765662566205258482\,n^{12} + 115251996734415746788867424992908\,n^{11}$$
$$+ 350009541864541416367247791701304\,n^{10} - 582640311164019997945305400048410\,n^{9}$$
$$- 1226855364790082831861072161428842\,n^{8}$$
$$+ 2131603202767134246764260942167468\,n^{7}$$
$$+ 2781972921871008261160393798476694\,n^{6}$$
$$- 5238760984190079515122721168798775\,n^{5}$$
$$- 3319164360796625033312714915543927\,n^{4}$$
$$+ 7598127033289977307530432957715278\,n^{3} + 6737869887406757638458493001 57584\,n^{2}$$
$$- 4809743999756002299533990429094015\,n + 1603247999918667433177996809698005\,)\,n$$
$$/\,43890$$

$$s_{n,55} := (114\,n^{52} + 2964\,n^{51} + 23218\,n^{50} - 49400\,n^{49} - 1320120\,n^{48} + 2689640\,n^{47}$$
$$+ 84069452\,n^{46} - 170828544\,n^{45} - 5140289849\,n^{44} + 10451408242\,n^{43} + 291752917965\,n^{42}$$
$$- 593957244172\,n^{41} - 1215958081361\,n^{40} + 31025873406894\,n^{39} + 725193562452973\,n^{38}$$
$$- 1481412998312840\,n^{37} - 31435790122515900\,n^{36} + 64352993243344640\,n^{35}$$
$$+ 1233328650408922520\,n^{34} - 2531010294061189680\,n^{33} - 4355964668693117 6594\,n^{32}$$
$$+ 89650303667923542868\,n^{31} + 1376712363665813987958\,n^{30}$$
$$- 2843075030999551518784\,n^{29} - 38675373481215722060035\,n^{28}$$
$$+ 80193821993430995638854\,n^{27} + 958395124150254395740231\,n^{26}$$
$$- 1996984070293939787119316\,n^{25} - 20767139449542614017039271\,n^{24}$$
$$+ 435312629693791678211 97858\,n^{23} + 389519503417038936123261875\,n^{22}$$

$$- 8225702698034570400677721608 \, n^{21} - 6249230885722349199358370304 \, n^{20}$$
$$+ 13321032041248155438784462216 \, n^{19} + 84542432995731445864493420572 \, n^{18}$$
$$- 182405898032711047167771303360 \, n^{17} - 947745737856494036138660492086 \, n^{16}$$
$$+ 2077897373745699119445092287532 \, n^{15} + 8612225628104151560568312608922 \, n^{14}$$
$$- 19302348629954002240581717505376 \, n^{13} - 6167310989347986857038962202377 \, n^{12}$$
$$+ 142556970608649975954659641910130 \, n^{11} + 333494858985395639198309194452317 \, n^{10}$$
$$- 809546688579441254351278030814764 \, n^{9} - 1282149562631635957956294953690529 \, n^{8}$$
$$+ 3373845813842713170263867938195822 \, n^{7}$$
$$+ 3119978004748878762030557686920445 \, n^{6}$$
$$- 9613801823340470694324983312036712 \, n^{5}$$
$$- 3465013110360417724213790601854506 \, n^{4}$$
$$+ 16543828044061306142752564515745724 \, n^{3}$$
$$- 1858922022355983338664295019080842 \, n^{2}$$
$$- 12825983999349339465423974477584040 \, n + 6412991999674669732711987238792020 \, )$$
$$(n+1)^{2} \, n^{2} \, / \, 6384$$

$$s_{n, 56} := (2 \, n + 1) \, (n + 1) \, (435 \, n^{54} + 11745 \, n^{53} + 97875 \, n^{52} - 152685 \, n^{51} - 5547555 \, n^{50}$$
$$+ 8397675 \, n^{49} + 366019875 \, n^{48} - 553228650 \, n^{47} - 23313133350 \, n^{46} + 35246314350 \, n^{45}$$
$$+ 1381863886650 \, n^{44} - 2090418987150 \, n^{43} - 75430548446010 \, n^{42} + 114191032162590 \, n^{41}$$
$$+ 3771460909300620 \, n^{40} - 5714286880032225 \, n^{39} - 171947700282889095 \, n^{38}$$
$$+ 260778693864349755 \, n^{37} + 7115453224041031665 \, n^{36} - 10803569182993722375 \, n^{35}$$
$$- 265918957995913926105 \, n^{34} + 404280221585367750345 \, n^{33}$$
$$+ 8925356233432869698085 \, n^{32} - 13590174460941988422300 \, n^{31}$$
$$- 267386399199603752710980 \, n^{30} + 407874686029876623277620 \, n^{29}$$
$$+ 7100054020156870639127580 \, n^{28} - 10854018373250244270330180 \, n^{27}$$
$$- 165789342945054862880601060 \, n^{26} + 25411023604207416456066680 \, n^{25}$$
$$+ 3373570304755489938620070310 \, n^{24} - 5187410968935338616158138805 \, n^{23}$$
$$- 59197880419611088832396626259 \, n^{22} + 91390526113884302556674008791 \, n^{21}$$
$$+ 884835776989845444517627060693 \, n^{20} - 1372948928541710318054777595435 \, n^{19}$$
$$- 11101722326703555764348640233189 \, n^{18} + 1733905795432618880550349147501 \, n^{17}$$
$$+ 114853312848085896175010084842093 \, n^{16} - 180949498249291938665290301836890 \, n^{15}$$
$$- 958196686703077107832829201513654 \, n^{14}$$
$$+ 1527769779179261631081888953188926 \, n^{13}$$
$$+ 6263675268176186368927129153991338 \, n^{12}$$
$$- 10159397791853910368931638207581470 \, n^{11}$$
$$- 30853141529008321090502095422449114 \, n^{10}$$

$$+ 513594111894394368202189622374644406 \, n^9$$
$$+ 108146660081507959976500757686354424 8 \, n^8$$
$$- 187899606817339118057620846413913575 \, n^7$$
$$- 245229326695266064607874841003537649 \, n^6$$
$$+ 461793793451568655940622684712263261 \, n^5$$
$$+ 292582445713258754238915099404224183 \, n^4$$
$$- 669770565295672459328683991462467905 \, n^3$$
$$- 593939388431503634772427305007854 71 \, n^2$$
$$+ 423976190912561774880206091482412159 \, n$$
$$- 141325396970853924960068697160804053 \, ) \, n \, / \, 49590$$

$$s_{n,\,57} := ( 30 \, n^{54} + 810 \, n^{53} + 6615 \, n^{52} - 14040 \, n^{51} - 402805 \, n^{50} + 819650 \, n^{49} + 27674475 \, n^{48}$$
$$- 56168600 \, n^{47} - 1832134586 \, n^{46} + 3720437772 \, n^{45} + 112980991667 \, n^{44}$$
$$- 229682421106 \, n^{43} - 6425376582795 \, n^{42} + 13080435586696 \, n^{41} + 335267431165903 \, n^{40}$$
$$- 683615297918502 \, n^{39} - 15980741611356244 \, n^{38} + 32645098520630990 \, n^{37}$$
$$+ 692755286039705439 \, n^{36} - 1418155670600041868 \, n^{35} - 27179242382564252633 \, n^{34}$$
$$+ 55776640435728547134 \, n^{33} + 959938971156992087315 \, n^{32}$$
$$- 1975654582749712721764 \, n^{31} - 3033909836733995493 7362 \, n^{30}$$
$$+ 6265385131742962259648 8 \, n^{29} + 85230300193271617195545 1 \, n^{28}$$
$$- 176725985518286196650739 0 \, n^{27} - 21120495757547231434478343 \, n^{26}$$
$$+ 44008251370277324835464076 \, n^{25} + 457652874258317386315531591 \, n^{24}$$
$$- 959313999886912097466527258 \, n^{23} - 8583980523813162697124600046 \, n^{22}$$
$$+ 18127275047513237491715727350 \, n^{21} + 137716534753931080036433105471 \, n^{20}$$
$$- 293560344555375397564581938292 \, n^{19} - 1863091814919761833847631609977 \, n^{18}$$
$$+ 4019743974394899065259845158246 \, n^{17} + 20885811589374644003560198982835 \, n^{16}$$
$$- 45791367151344187072380243123916 \, n^{15} - 189790694539302284912857404256768 \, n^{14}$$
$$+ 425372756231748756898095051637452 \, n^{13}$$
$$+ 1358103080474189836525054296736339 \, n^{12}$$
$$- 3141578917180128429948203645110130 \, n^{11}$$
$$- 7349345412599959679324643622751979 \, n^{10}$$
$$+ 17840269742380047788597490890614088 \, n^9$$
$$+ 28255188206084728918824652677665903 \, n^8$$
$$- 74350646154549505626246796245945894 \, n^7$$
$$- 68756070502479694391051989062846470 \, n^6$$
$$+ 211862787159508894408350774371638834 \, n^5$$
$$+ 76359732454951258497306293304993927 \, n^4$$
$$- 364582252069411411402963360981626688 \, n^3$$

$$+ 4096572906385178074141298330009291\, n^2$$

$$+ 2826507939417078499201373943216081 06\, n$$

$$- 14132539697085392496006869716080404 53 \big) \, (n+1)^2\, n^2 \,/\, 1740$$

$$s_{n,58} := (2n+1)(n+1)(15\, n^{56} + 420\, n^{55} + 3640\, n^{54} - 5670\, n^{53} - 220878\, n^{52} + 334152\, n^{51}$$

$$+ 15701166\, n^{50} - 23718825\, n^{49} - 1081008045\, n^{48} + 1633371480\, n^{47} + 69486306240\, n^{46}$$

$$- 105046145100\, n^{45} - 4127534207460\, n^{44} + 6243824383740\, n^{43} + 225422170746270\, n^{42}$$

$$- 341255168311275\, n^{41} - 1272335097855215\, n^{40} + 1709130230938460\, n^{39}$$

$$+ 513942103259923080\, n^{38} - 779452720005353850\, n^{37} - 2126785745656748 4930\, n^{36}$$

$$+ 3229151254485390432 0\, n^{35} + 79482460003407330231 0\, n^{34}$$

$$- 1208382656323536905625\, n^{33} - 2667766181752117285444 5\, n^{32}$$

$$+ 4062068405444352773448 0\, n^{31} + 79921121036973002619600 0\, n^{30}$$

$$- 12191271575818168031612 40\, n^{29} - 2122188318321269001056134 4\, n^{28}$$

$$+ 3244238835360994341742263 6\, n^{27} + 49554018653081493795788641 8\, n^{26}$$

$$- 759531473973027378645540945\, n^{25} - 10083517017683486786914062445\, n^{24}$$

$$+ 15505041263511743869693864140\, n^{23} + 17694097971724451052974938796 0\, n^{22}$$

$$- 27316399020762263772947101401 0\, n^{21} - 26447519434019341983985242704 50\, n^{20}$$

$$+ 4103709910206712616462521912680\, n^{19} + 33182769574758425262926984115530\, n^{18}$$

$$- 51826009317240994202621737129635\, n^{17} - 343293671288274165753567828462895\, n^{16}$$

$$+ 5408535115910317457316626112591 60\, n^{15}$$

$$+ 2864025862537140018958549784093120\, n^{14}$$

$$- 4566465549601225901303655981769260\, n^{13}$$

$$- 1872196826762258367282591759798994 0\, n^{12}$$

$$+ 3036618517623448845989070438786954 0\, n^{11}$$

$$+ 9221926615470752161523091988382081 0\, n^{10}$$

$$- 15351199182017852665279173201966598 5\, n^{9}$$

$$- 32324747724363799931745217117365092 5\, n^{8}$$

$$+ 5616272117755462623025741227703093 80\, n^{7}$$

$$+ 7329843065150480407856646417724770 80\, n^{6}$$

$$- 13802900656603451923297840240438703 10\, n^{5}$$

$$- 874521591522803891217101947275201598\, n^{4}$$

$$+ 2001927420114378432990544932934737552\, n^{3}$$

$$+ 1775269934506742742902861156287810 26\, n^{2}$$

$$- 1267254200233200627930701639910540315\, n$$

$$+ 4224180667444002093102338799708010 5 \big) \, n \,/\, 1770$$

$$s_{n,59} := (6\, n^{56} + 168\, n^{55} + 1428\, n^{54} - 3024\, n^{53} - 92907\, n^{52} + 188838\, n^{51} + 6867211\, n^{50}$$

$$- 13923260\, n^{49} - 490744860\, n^{48} + 995412980\, n^{47} + 32769931430\, n^{46} - 66535275840\, n^{45}$$
$$- 2024879515745\, n^{44} + 4116294307330\, n^{43} + 115213581449685\, n^{42} - 234543457206700\, n^{41}$$
$$- 601240308428190\, n^{40} + 12259349625769080\, n^{39} + 286593493008401980\, n^{38}$$
$$- 585446335642573040\, n^{37} - 12423740112690588399\, n^{36} + 25432926561023749838\, n^{35}$$
$$+ 487428164963108614823\, n^{34} - 1000289256487240979484\, n^{33}$$
$$- 17215398653314921030340\, n^{32} + 35431086563117083040164\, n^{31}$$
$$+ 544096805910117035917422\, n^{30} - 1123624698383351154875008\, n^{29}$$
$$- 15285073644943827377821309\, n^{28} + 31693771988271005910517626\, n^{27}$$
$$+ 378771791644100622019777817\, n^{26} - 789237355276472249950073260\, n^{25}$$
$$- 8207477768557687881465429000\, n^{24} + 17204192892391848012880931260\, n^{23}$$
$$+ 153943814848700499185703399970\, n^{22} - 325091822589792846384287731200\, n^{21}$$
$$- 2469787608523860706144984035055\, n^{20} + 5264667039637514258674255801310\, n^{19}$$
$$+ 33412408222013407194347204995335\, n^{18} - 72089483483664328647368665791980\, n^{17}$$
$$- 374563001826792483210854203930874\, n^{16} + 821215487137249295069077073653728\, n^{15}$$
$$+ 3403677753258361299311713105623968\, n^{14}$$
$$- 7628570993653971893692503284901664\, n^{13}$$
$$- 24356016257084196284879896553716545\, n^{12}$$
$$+ 56340603507822364463452296392334754\, n^{11}$$
$$+ 131802054587578831589953619768984437\, n^{10}$$
$$- 319944712682980027643359535930303628\, n^{9}$$
$$- 506724293014722833281283166430856176\, n^{8}$$
$$+ 1333393298712425694205925868792015980\, n^{7}$$
$$+ 1233061020925587687863350805665501866\, n^{6}$$
$$- 3799515340563601069932627480123019712\, n^{5}$$
$$- 1369423950065778525954932832783768011\, n^{4}$$
$$+ 6538363240695158121842493145690555734\, n^{3}$$
$$- 734673219881177805059843293024197237\, n^{2}$$
$$- 506901680093280251172280655964216120\, n$$
$$+ 2534508400466401255861403279821080630\,) \, (n+1)^2 \, n^2 \, / \, 360$$

$s_{n,\,60} := (2\,n+1)\,(n+1)\,(465465\, n^{58} + 13498485\, n^{57} + 121486365\, n^{56} - 188978790\, n^{55}$
$$- 7874116250\, n^{54} + 11905663770\, n^{53} + 601438497660\, n^{52} - 908110578375\, n^{51}$$
$$- 44629028569125\, n^{50} + 67397598142875\, n^{49} + 3101090582164575\, n^{48}$$
$$- 4685334672318300\, n^{47} - 199764091895039100\, n^{46} + 301988805178717800\, n^{45}$$
$$+ 11872285438595932350\, n^{44} - 17959422560483257425\, n^{43} - 648478527418302272835\, n^{42}$$
$$+ 981697502407695037965\, n^{41} + 32428474698033785644845\, n^{40}$$

$$- 49133560798254525986250\ n^{39} - 147853022551251589538967 0\ n^{38}$$

$$+ 224236211866790110607763 0\ n^{37} + 611843832387017295042806 40\ n^{36}$$

$$- 92897755917386544809459775\ n^{35} - 228659027718040149921452404 5\ n^{34}$$

$$+ 3476334293729295521226515955\ n^{33} + 7674761109905525239038760825 5\ n^{32}$$

$$- 116859583795447526346194670360\ n^{31} - 2299210228925332011266947867040\ n^{30}$$

$$+ 35072451352857217800735191357 40\ n^{29} + 6105216068777944538911504487331 0\ n^{28}$$

$$- 93331863599312028973709326877835\ n^{27} - 142559446241726743513742683900174 5\ n^{26}$$

$$+ 21850576254255571671929949219415 35\ n^{25}$$

$$+ 29008759367932668465519102037951495\ n^{24}$$

$$- 446056678646117812818751505178980 10\ n^{23}$$

$$- 5090325423631685443981229508645766 46\ n^{22}$$

$$+ 78585164747705870723812200155581397 4\ n^{21}$$

$$+ 76085529076118472055627406156768656 12\ n^{20}$$

$$- 118057551851563001619631719242932054 05\ n^{19}$$

$$- 95461829061892613766414540045373027271\ n^{18}$$

$$+ 149095621185417070730603396030206143609\ n^{17}$$

$$+ 987604174894594863769298445626811520597\ n^{16}$$

$$- 1555954072934600831019249366455320352700\ n^{15}$$

$$- 8239370939272265249372559997566454277036\ n^{14}$$

$$+ 13137033445375698289568464679577341591904\ n^{13}$$

$$+ 53860282229980909740405860737990132254562\ n^{12}$$

$$- 87358940067671498755393023446773869177795\ n^{11}$$

$$- 265300935838293557798120208191840759923001\ n^{10}$$

$$+ 441630873791276086074876824011148074473399\ n^{9}$$

$$+ 929934294599971116813336468785035439975527\ n^{8}$$

$$- 1615716878795594718257443115183127197199990\ n^{7}$$

$$- 2108685425310108970308752206575884627746106\ n^{6}$$

$$+ 3970886577362960814591849867455390540219154\ n^{5}$$

$$+ 2515866871598946473428825155862684475516472\ n^{4}$$

$$- 5759243596079900117439162667521721983384285\ n^{3}$$

$$- 5107184156075444382675557656523998309954887\ n^{2}$$

$$+ 3645699421451266716120914982239460738124473\ n$$

$$- 1215233140483755572040304994079820246041491\ )\ n\ /\ 56786730$$

$$s_{n,\,61} := (\ 30030\ n^{58} + 870870\ n^{57} + 7692685\ n^{56} - 16256240\ n^{55} - 533583050\ n^{54}$$

$$+ 1083422340\ n^{53} + 42321019455\ n^{52} - 85725461250\ n^{51} - 3255349740500\ n^{50}$$

$$+ 6596424942250\ n^{49} + 234668074094025\ n^{48} - 475932573130300\ n^{47}$$

$$-15702126816828166\,n^{46} + 31880186206786632\,n^{45} + 970723971104629027\,n^{44}$$

$$-1973328128416044686\,n^{43} - 55239851935118556420\,n^{42} + 112453031998653157526\,n^{41}$$

$$+2882768486841101065493\,n^{40} - 5877990005680855288512\,n^{39}$$

$$-137414057938420418254814\,n^{38} + 280706105882521691798140\,n^{37}$$

$$+5956867976568109356228759\,n^{36} - 12194442059018740404255658\,n^{35}$$

$$-233709533272965412934725168\,n^{34} + 47961350860949566273705994\,n^{33}$$

$$+8254351116180732382900659865\,n^{32} - 16988315740966414332075025724\,n^{31}$$

$$-260880754818833277200050145882\,n^{30} + 538749825378632968732175317488\,n^{29}$$

$$+7328809013592237596434680716331\,n^{28} - 15196367852563108161601536750150\,n^{27}$$

$$-181611563638112122060356048133610\,n^{26} + 378419495128787352282313633017370\,n^{25}$$

$$+3935279512646590345887140353722145\,n^{24}$$

$$-8248978520421968044056594340461660\,n^{23}$$

$$-73812193929667251502108740230098504\,n^{22}$$

$$+155873366379756471048274074800658668\,n^{21}$$

$$+1184201145777796141775704281058497543\,n^{20}$$

$$-2524275657935348754599682636917653754\,n^{19}$$

$$-16020410809147305979226871302491263230\,n^{18}$$

$$+34565097276229960713053425241900180214\,n^{17}$$

$$+179593554685042217149519274748577747277\,n^{16}$$

$$-393752206646314395012091974739055674768\,n^{15}$$

$$-1631978021665050996906843064034028407916\,n^{14}$$

$$+3657708249976416388825778102807112490600\,n^{13}$$

$$+11678098253816346856773356298063032802691\,n^{12}$$

$$-27013904757609110102372490698933178095982\,n^{11}$$

$$-63195775831399914098030321469627104485702\,n^{10}$$

$$+15340545642040893829843133638187387067386\,n^{9}$$

$$+242961575446494859843361254428531387592305\,n^{8}$$

$$-639328607313398657985155642495250162251996\,n^{7}$$

$$-59122179930110373958834942741254144511644\,n^{6}$$

$$+1821772205915606137161854497337758451275284\,n^{5}$$

$$+656604399964058070345752359624651227947151\,n^{4}$$

$$-3134981005843722277853359216587060907169586\,n^{3}$$

$$+352257362438105566886374614213710207543302\,n^{2}$$

$$+24304662809675114480609988159640492082982\,n$$

$$-1215233140483755572040304994079820246041491\,)\,(\,n+1\,)^{2}\,n^{2}\,/\,1861860$$

$$
\begin{aligned}
s_{n,62} := {} & (2n+1)(n+1)\,(2145\,n^{60} + 64350\,n^{59} + 600600\,n^{58} - 933075\,n^{57} - 41490735\,n^{56} \\
& + 62702640\,n^{55} + 3396766230\,n^{54} - 5126500665\,n^{53} - 270920602953\,n^{52} \\
& + 408944154762\,n^{51} + 20290441881996\,n^{50} - 30640134900375\,n^{49} - 1412948891148795\,n^{48} \\
& + 2134743404173380\,n^{47} + 9106591485358890\,n^{46} - 137666358930125025\,n^{45} \\
& - 541288513806791385\,n^{44} + 8188160886565249590\,n^{43} + 29566767112951440720\,n^{42} \\
& - 447595587137554230375\,n^{41} - 1478557330563132411091\,n^{40} \\
& + 22402157752015763281560\,n^{39} + 67412851555592117437863 0\,n^{38} \\
& - 102239385220988964320872 5\,n^{37} - 278967291811286781670309 65\,n^{36} \\
& + 423562906977979620721508 10\,n^{35} + 104256010963925921925393954 0\,n^{34} \\
& - 158501830980778780991698471 5\,n^{33} - 349927142566974591377571894 3\,n^{32} \\
& + 532815805399501377756220707 72\,n^{31} + 104831415579987857629786128132 6\,n^{30} \\
& - 159911202396979293333460295737 5\,n^{29} - 278364473156035548411167031119 91\,n^{28} \\
& + 425542269853902287283423561466 74\,n^{27} + 64999313237617292321636730258787 2\,n^{26} \\
& - 9962668120569544991887221319551 45\,n^{25} \\
& - 13226408257746754253547565212911 165\,n^{24} \\
& + 20337745792648608629915708885344 320\,n^{23} \\
& + 23209100866927599400852608355578 8050\,n^{22} \\
& - 35830538590023829532774697977635 4235\,n^{21} \\
& - 34690841388138863392268771516223 33515\,n^{20} \\
& + 53827789011709486565041892173216 77390\,n^{19} \\
& + 43525374809440841455514137734804 403060\,n^{18} \\
& - 67979451664746736511523301210867 443285\,n^{17} \\
& - 45029350786690921406661868675826 0460065\,n^{16} \\
& + 70942998763273718935568968074282 4411740\,n^{15} \\
& + 37567026721589711094016167327586 89623710\,n^{14} \\
& - 59897690020548252587802699395094 46641435\,n^{13} \\
& - 24557343960837783900330987427944 303212195\,n^{12} \\
& + 39830900442284088479886616111671 178139010\,n^{11} \\
& + 12096272172850104388724998501252 3519596520\,n^{10} \\
& - 20135953284983610070818285574620 868464285\,n^{9} \\
& - 42399919547312173235729279695564 8990411377\,n^{8} \\
& + 73667855965712940357134833822078 3919849208\,n^{7} \\
& + 96144526449795006336217985458610 5608424914\,n^{6} \\
& - 18105071765754897968289439509895 50372561975\,n^{5} \\
& - 11470977419262780697427685103817 16060343671\,n^{4}
\end{aligned}
$$

$$+\ 262590020117716200302862474106734927679 6494\ n^3$$
$$+\ 23285967469783108503500778039768797696 6078092\ n^2$$
$$-\ 166223961263532762906682404113020658751 5385\ n$$
$$+\ 554079870878442543022274680376735529171 795\ )\ n\ /\ 270270$$

$s_{n,\,63} := (\ 2145\ n^{60} + 64350\ n^{59} + 589875\ n^{58} - 1244100\ n^{57} - 43531059\ n^{56} + 88306218\ n^{55}$

$$+\ 3695966703\ n^{54} - 7480239624\ n^{53} - 305206311267\ n^{52} + 617892862158\ n^{51}$$

$$+\ 23683817994551\ n^{50} - 47985528851260\ n^{49} - 1710806757286935\ n^{48}$$

$$+\ 3469599043425130\ n^{47} + 114530497817373715\ n^{46} - 232530594678172560\ n^{45}$$

$$-\ 7081271549420094565\ n^{44} + 14395073693518361690\ n^{43} + 402977690966319207585\ n^{42}$$

$$-\ 820350455626156776860\ n^{41} - 21030106380838221296945\ n^{40}$$

$$+\ 42880563217302599370750\ n^{39} + 1002451265517477452643845\ n^{38}$$

$$-\ 2047784894252257504658440\ n^{37} - 43456091499677845714774545\ n^{36}$$

$$+\ 88959967893607948934207530\ n^{35} + 1704940250188911334420928365\ n^{34}$$

$$-\ 3498840468271430617776064260\ n^{33} - 60216524498694135466022979837\ n^{32}$$

$$+\ 123931889465659701549822023934\ n^{31} + 1903157770473643033551144161169\ n^{30}$$

$$-\ 3930247430412945768652110346272\ n^{29} - 53464579456964361402186277283406\ n^{28}$$

$$+\ 110859406344341668573024664913084\ n^{27}$$

$$+\ 1324879098387388200551224915994598\ n^{26}$$

$$-\ 2760617603119118069675474496902280\ n^{25}$$

$$-\ 28708356828478522439994481086875630\ n^{24}$$

$$+\ 60177331260076162949664436670653540\ n^{23}$$

$$+\ 538469197642173654575436965966735190\ n^{22}$$

$$-\ 1137115726544423472100538368604123920\ n^{21}$$

$$-\ 8638895646765608804278985774115742930\ n^{20}$$

$$+\ 18414907020075641080658509916835609780\ n^{19}$$

$$+\ 116870902964549806703055007351283075770\ n^{18}$$

$$-\ 252156712949175254486768524619401761320\ n^{17}$$

$$-\ 1310157470535987771669186981175389708210\ n^{16}$$

$$+\ 2872471654021150797825142486970181177740\ n^{15}$$

$$+\ 11905484028003454305447673244546302059530\ n^{14}$$

$$-\ 26683439710028059408720488976062785296800\ n^{13}$$

$$-\ 85193189119309480423209031175737244915900\ n^{12}$$

$$+\ 197069817948647020255138551327537275128600\ n^{11}$$

$$+\ 461021098207196169773112113508357026996700\ n^{10}$$

$$-\ 1119112014363039359801362778344251329122000\ n^9$$

$$-\ 177243511707691278627649446514688654640 3604\ n^8$$

$$+ 466398224851686493235435170863802442192920 \, n^7$$

$$+ 431303706002451605782734921388336245770062 \, n^6$$

$$- 132900563685658970480090501364047493373304 \, n^5$$

$$- 479001131921702382825000501765777430283311 \, n^4$$

$$+ 228700790069999447045090601717202979429966 \, n^3$$

$$- 256976156944489166389813519983238050474962 \, n^2$$

$$- 177305558681101613767127897720555369334974 \, n$$

$$+ 88652779340550806883563948860277684667487) \, (n+1)^2 \, n^2 \, / \, 137280$$

$$s_{n,\,64} := (2\,n+1)\,(n+1)\,(2805\,n^{62} + 86955\,n^{61} + 840565\,n^{60} - 1304325\,n^{59} - 61767035\,n^{58}$$

$$+ 93302715\,n^{57} + 5407354645\,n^{56} - 8157683325\,n^{55} - 462397328195\,n^{54}$$

$$+ 697674833955\,n^{53} + 37225118770365\,n^{52} - 56186515572525\,n^{51} - 2794011837810675\,n^{50}$$

$$+ 4219111014502275\,n^{49} + 194666487975380325\,n^{48} - 294109287470321625\,n^{47}$$

$$- 12548071159343852775\,n^{46} + 18969161382750939975\,n^{45} + 745870467512663518425\,n^{44}$$

$$- 1128290281960370747625\,n^{43} - 40741955724050405185815\,n^{42}$$

$$+ 61677078727205579315255\,n^{41} + 2037403517295994456062105\,n^{40}$$

$$- 3086943815307519580669425\,n^{39} - 92892743972918792931979215\,n^{38}$$

$$+ 1408825878670319491883035355\,n^{37} + 38440801319935575716281075055\,n^{36}$$

$$- 583656149192385233203631305\,n^{35} - 143661455685653899104194862575\,n^{34}$$

$$+ 218410464274442774822310450375\,n^{33} + 482188437763210009019591874092575\,n^{32}$$

$$- 73420317985537152270503336575\,n^{31} - 144454346203723200094520563130945\,n^{30}$$

$$+ 220352535204877485903133361364705\,n^{29}$$

$$+ 3835777363441145214490489821776601575\,n^{28}$$

$$- 586383671921961696030891400733137557535\,n^{27}$$

$$- 89566980000013517422479553047230275\,n^{26}$$

$$+ 13728238835965683609352638996075022557525\,n^{25}$$

$$+ 18225568624004205673958553352593533375\,n^{24}$$

$$- 28024764877804592691405461978694051757579\,n^{23}$$

$$- 3198140056705128049926283522504023789757\,n^{22}$$

$$+ 49373333909446715038346452593649505943357\,n^{21}$$

$$+ 4780287270950647322878219556608765277993\,n^{20}$$

$$- 741729760189830673623465196459562321415\,n^{19}$$

$$- 59976577920690730359407068479195406475775\,n^{18}$$

$$+ 93673515681985248907227928701090921320735\,n^{17}$$

$$+ 62049008836806766449775171965413604528767977\,n^{16}$$

$$- 97757189039309409191988975916258613997555\,n^{15}$$

$$- 517661643416249993136904075369262668897357\, n^{14}$$

$$+ 825371059644029694301350601012023310344813\, n^{13}$$

$$+ 3383923653824771257657447184430197585513459\, n^{12}$$

$$- 54885710105591717333636846077151308033442595\, n^{11}$$

$$- 16668277160738443687528775415549500981079517\, n^{10}$$

$$+ 27746701246387251398111586161899905488340573\, n^{9}$$

$$+ 58425736474280464435875412475021674456957299\, n^{8}$$

$$- 101511955334614322352868911793482464429606235\, n^{7}$$

$$- 132484090200503821871685923715024073271007077\, n^{6}$$

$$+ 249482112968062893983963341469277342121313733\, n^{5}$$

$$+ 158066409312975824664747708325851379405662139\, n^{4}$$

$$- 361840670453495183989103233223415740169150075\, n^{3}$$

$$- 32087320293617409741163441912368739834034285\, n^{2}$$

$$+ 229051315667173706606296779480260979835626465\, n$$

$$- 76350438555724568868765593160086993278542155\, )\, n\, /\, 364650$$

$$s_{n,\,65} := (\, 510\, n^{62} + 15810\, n^{61} + 150195\, n^{60} - 316200\, n^{59} - 11770035\, n^{58} + 23856270\, n^{57}$$

$$+ 1067342535\, n^{56} - 2158541340\, n^{55} - 94390985895\, n^{54} + 190940513130\, n^{53}$$

$$+ 7864031179035\, n^{52} - 15919002871200\, n^{51} - 611510259166475\, n^{50}$$

$$+ 1238939521204150\, n^{49} + 44194742427118575\, n^{48} - 89628424375441300\, n^{47}$$

$$- 2958994336851030787\, n^{46} + 6007617098077502874\, n^{45} + 182956303642475525539\, n^{44}$$

$$- 371920224383028553952\, n^{43} - 10411669571372408734515\, n^{42}$$

$$+ 21195259367127846022982\, n^{41} + 543352466525631155742551\, n^{40}$$

$$- 1107900192418390157508084\, n^{39} - 25900254349469443719284183\, n^{38}$$

$$+ 52908408891357277596076450\, n^{37} + 1122770729109048718210503243\, n^{36}$$

$$- 2298449867109454714017082936\, n^{35} - 44050373476219475758524320011\, n^{34}$$

$$+ 90399196819548406231065722958\, n^{33} + 1555808428570768430922155202495\, n^{32}$$

$$- 3202016053961085268075376127948\, n^{31} - 49171700454572788718352633687309\, n^{30}$$

$$+ 101545416963106662704780643502566\, n^{29}$$

$$+ 1381359090572199503450918088525027\, n^{28}$$

$$- 2864263598107505669606616820552620\, n^{27}$$

$$- 34230771199019208608341235330154619\, n^{26}$$

$$+ 71325805996145922886289087480861858\, n^{25}$$

$$+ 741734997038108853835057750031909703\, n^{24}$$

$$- 1554795800072363630556404587544681264\, n^{23}$$

$$- 13912375797288039647772584046118548239\, n^{22}$$

$$+ 29379547394648442926101572679781777742\, n^{21}$$

$$+ 2232022987346110626022631021631044491755 \, n^{20}$$
$$- 475784144863870568130627777005990761252 \, n^{19}$$
$$- 30195820465375717519220950016666388896051 \, n^{18}$$
$$+ 65149482379390140719748177803392685533354 \, n^{17}$$
$$+ 338504099465015527974917927621920496686543 \, n^{16}$$
$$- 742157681309421196669584033047233679264 40 \, n^{15}$$
$$- 3076008220558244632475827724304670253333563 \, n^{14}$$
$$+ 6894174122425910461621239481656574185933566 \, n^{13}$$
$$+ 22011280637576604728291006622968531005 82931 \, n^{12}$$
$$- 50916735397579119918203252727593636197 59428 \, n^{11}$$
$$- 119113568553831070541787252074021523511 95515 \, n^{10}$$
$$+ 289143872505241261001777756875636683221 50458 \, n^{9}$$
$$+ 457942321178269745431015778423715979058 22999 \, n^{8}$$
$$- 1205028514861780751863809313723068641337 96456 \, n^{7}$$
$$- 1114355151038129119900016172990421888625 55423 \, n^{6}$$
$$+ 3433738816938038991663841659703912418589 07302 \, n^{5}$$
$$+ 1237590522134324957145079233666427390729 34619 \, n^{4}$$
$$- 5908919861206688905954000127036767200047 76540 \, n^{3}$$
$$+ 6639467739316073869140322687157738016676 1805 \, n^{2}$$
$$+ 4581026313343474132125935589605215961712 52930 \, n$$
$$- 2290513156671737066062967794802609798356 26465 \,) (n + 1)^2 \, n^2 \, / \, 33660$$

$$s_{n,\,66} := (2\,n + 1)\,(n + 1)\,(41055\,n^{64} + 1313760\,n^{63} + 13137600\,n^{62} - 20363280\,n^{61}$$
$$- 1024951760\,n^{60} + 1547609280\,n^{59} + 95741354800\,n^{58} - 144385836840\,n^{57}$$
$$- 8757143717000\,n^{56} + 13207908493920\,n^{55} + 755886528360320\,n^{54}$$
$$- 1140433746787440\,n^{53} - 60985187455306128\,n^{52} + 9204799805635291 2\,n^{51}$$
$$+ 4579790470474799496\,n^{50} - 6915709704740375700\,n^{49} - 31912796648721028 3620\,n^{48}$$
$$+ 482149804583185613280\,n^{47} + 2057143424324327031264 0\,n^{46}$$
$$- 3109822626271564982756 00\,n^{45} - 12227972813362299467661 60\,n^{44}$$
$$+ 18497450351379231692870 40\,n^{43} + 6679343555831684415362 4720\,n^{42}$$
$$- 10111502585504422781508 0600\,n^{41} - 3340174636700653189435 621560\,n^{40}$$
$$+ 50608194679785018980609 72640\,n^{39} + 1522909047585011535134 13343680\,n^{38}$$
$$- 23096676687174098121915 0501840\,n^{37} - 6302090279842146188057 079492640\,n^{36}$$
$$+ 95686188031990897726951 94489880\,n^{35} + 2355225279147148289130 81670125700\,n^{34}$$
$$- 35806810127367178825597 0102433490\,n^{33}$$
$$- 79051294210845675379864 53338556026\,n^{32}$$

$$+ 12036728182263687201107665059050784 \, n^{31}$$

$$+ 23682241490322103346753053459974912 \, n^{30}$$

$$- 36125198644596339380184963442918776 0 \, n^{29}$$

$$- 62884724419663624345881682893990192 48 \, n^{28}$$

$$+ 96133346561725253487831772513131227 52 \, n^{27}$$

$$+ 146838562214088395666695039257384063696 \, n^{26}$$

$$- 22506451064921885617443414751173265 6920 \, n^{25}$$

$$- 29879496800046451069018438673902694 07160 \, n^{24}$$

$$+ 45944567753315770884399828748412704 39200 \, n^{23}$$

$$+ 52431184761345936046747038579864017 959040 \, n^{22}$$

$$- 80944005529684692614340549307216662 158160 \, n^{21}$$

$$- 78369339888709550720062770317849550 0576880 \, n^{20}$$

$$+ 12160121010954856071081118294213515 81944400 \, n^{19}$$

$$+ 98327245916624087031130282154531009 95636040 \, n^{18}$$

$$- 15357092938041355858223598237890327 284426260 \, n^{17}$$

$$- 10172484596982736858317090597662220 148794980 \, n^{16}$$

$$+ 16026581542376173080386815808387849 3865405600 \, n^{15}$$

$$+ 84866868831866415471100367279527571 4642539840 \, n^{14}$$

$$- 13531359401898770974684395882348528 18896512560 \, n^{13}$$

$$- 55476971979412126173045123676464071 087278095856 \, n^{12}$$

$$+ 89981137670067574746909883453135330 40365400064 \, n^{11}$$

$$+ 27326430486875450676932933052566120 197700203632 \, n^{10}$$

$$- 45488702613816554752744893751505946 816733005480 \, n^{9}$$

$$- 95784753937832715996198804946536831 22514111752 \, n^{8}$$

$$+ 16642148221365735137067065429555649 8092137670368 \, n^{7}$$

$$+ 21719804911862987373386128326306797 1477685357024 \, n^{6}$$

$$- 40900781478477348628612725204238020 6262596870720 \, n^{5}$$

$$- 25913840433222538073549664188578080 9999810921925880 \, n^{4}$$

$$+ 59321151389072481424630858884986145 3128481324180 \, n^{3}$$

$$+ 52604832464568789773391065705413471 38295739110 \, n^{2}$$

$$- 37551300564221559178316295428074274 7271684270755 \, n$$

$$+ 12517100188073853059438765142691424 9090561423585) \, n \, / \, 5501370$$

$$s_{n,\,67} := (\, 4830 \, n^{64} + 154560 \, n^{63} + 1519840 \, n^{62} - 3194240 \, n^{61} - 126247345 \, n^{60}$$

$$+ 255688930 \, n^{59} + 12202004045 \, n^{58} - 24659697020 \, n^{57} - 1152996182653 \, n^{56}$$

$$+ 2330652062326 \, n^{55} + 102880632643201 \, n^{54} - 208091917348728 \, n^{53}$$

$$- 8589253199561049 \, n^{52} + 17386598316470826 \, n^{51} + 668241223060200997 \, n^{50}$$

$$-\,1353869044436872820\,n^{49} - 48300754797122370345\,n^{48} + 9795537863868161 3510\,n^{47}$$

$$+\,3234005192219670117845\,n^{46} - 656596576307802184 9200\,n^{45}$$

$$-\,19996189093528746958 9985\,n^{44} + 40648974763365296102 9170\,n^{43}$$

$$+\,1137944404268634918825 7245\,n^{42} - 2316537783300635133754 3660\,n^{41}$$

$$-\,5938578625046847572552 48125\,n^{40} + 1210881102842375865848 039910\,n^{39}$$

$$+\,2830772358865204345713 7642945\,n^{38} - 5782632828014646278012 3325800\,n^{37}$$

$$-\,12271340625824605345227 10706057\,n^{36} + 25120944534450675318255 44737914\,n^{35}$$

$$+\,48144926414553881360997 078517429\,n^{34} - 98801947282552830253819 701772772\,n^{33}$$

$$-\,1700423322746718844389772 709892643\,n^{32}$$

$$+\,3499648592775990519033365 121558058\,n^{31}$$

$$+\,5374228907038868615931099 6706558427\,n^{30}$$

$$-\,1109842267335533628376553 58534674912\,n^{29}$$

$$-\,1509758639149206212623294 648072936961\,n^{28}$$

$$+\,3130501505031965788084244 654680548834\,n^{27}$$

$$+\,3741257642304900862914166 4269286302813\,n^{26}$$

$$-\,7795565435112998304636757 3193253154460\,n^{25}$$

$$-\,8106804576840545263450201 09033973239805\,n^{24}$$

$$+\,1699316569719239035736407 791261199634070\,n^{23}$$

$$+\,1520555349804280080132023 2513721047483745\,n^{22}$$

$$-\,3211042356580484063837687 2818703294601560\,n^{21}$$

$$-\,2439493113000829788908214 89492777286558825\,n^{20}$$

$$+\,5200090461659707984200198 51804257867719210\,n^{19}$$

$$+\,3300257053099810240296193 461519127230489205\,n^{18}$$

$$-\,7120523152365591279012406 774842512328697620\,n^{17}$$

$$-\,3699685998079190688985129 560212656680414 0505\,n^{16}$$

$$+\,8111424311394940505871499 7979095645936978630\,n^{15}$$

$$+\,3361928130726233486058526 13836282983981850245\,n^{14}$$

$$-\,7534998692591961022704202 2565166161390067912 0\,n^{13}$$

$$-\,2405726456587579237140847 119242315985331289393\,n^{12}$$

$$+\,5564952782434354576552114 46413629358456325790 6\,n^{11}$$

$$+\,1301853662795600184421785 2126286519795981491181\,n^{10}$$

$$-\,3160202603834635826498781 8716709333176526240268\,n^{9}$$

$$-\,5005087962800983383880866 7391208822754380227 49\,n^{8}$$

$$+\,1317037852943660259426051 53499127097727402285766\,n^{7}$$

$$+\,1217936254154374161870126 50179462588464341296577\,n^{6}$$

$$- 375291036125240858316630453858052274656084878920\, n^5$$
$$- 1352626551150263728515086207811752628053460922185\, n^4$$
$$+ 645816346355293604019647695420402800266777063290\, n^3$$
$$- 72566169416169740821048544856372901952265684475\, n^2$$
$$- 500684007522954122377550605707656996362245694340\, n$$
$$+ 250342003761477061188775302853828498181122847170\,)\,(n + 1)^2\, n^2\, /\, 328440$$

$$s_{n,68} := (2\,n + 1)\,(n + 1)\,(\,105\, n^{66} + 3465\, n^{65} + 35805\, n^{64} - 55440\, n^{63} - 2960496\, n^{62}$$
$$+ 4468464\, n^{61} + 294471232\, n^{60} - 443941080\, n^{59} - 28746411848\, n^{58} + 43341588312\, n^{57}$$
$$+ 2654164397656\, n^{56} - 4002917390640\, n^{55} - 229600506155376\, n^{54} + 346402217928384\, n^{53}$$
$$+ 18534065103282792\, n^{52} - 27974298763888380\, n^{51} - 1392029688073806420\, n^{50}$$
$$+ 2102031681492653820\, n^{49} + 97002237754148245260\, n^{48} - 146554372471968694800\, n^{47}$$
$$- 6525949680655565259120\, n^{46} + 9452701707219332236080\, n^{45}$$
$$+ 371685559950297751532640\, n^{44} - 562254690779056293417000\, n^{43}$$
$$- 20302766548073601460376760\, n^{42} + 30735277167499930337273640\, n^{41}$$
$$+ 1015291302354139594697945640\, n^{40} - 1538304592114959357215555280\, n^{39}$$
$$- 46290883722301024052714166048\, n^{38} + 70205477879509015757679026712\, n^{37}$$
$$+ 1915605734998100676471586900156\, n^{36} - 2908511341436905522586219863590\, n^{35}$$
$$- 71590263864692501116643559664994\, n^{34} + 108839651467757204436258449429286\, n^{33}$$
$$+ 2402871209010451087603417858822078\, n^{32}$$
$$- 3658726639249555233623256012947760\, n^{31}$$
$$- 71985382170445631426971114384340688\, n^{30}$$
$$+ 109807436575293224757268299582984912\, n^{29}$$
$$+ 1911466413349560429836624668974521856\, n^{28}$$
$$- 2922103338311987257133571153253275240\, n^{27}$$
$$- 44633570783513634698629003422772515960\, n^{26}$$
$$+ 68411407844426445676510290710785411560\, n^{25}$$
$$+ 9082277947249629458594395093790677802 00\, n^{24}$$
$$- 1396547396009657641627414409423994376080\, n^{23}$$
$$- 15937169099366986572463898297541312022672\, n^{22}$$
$$+ 24604027347055308679509554651023965222048\, n^{21}$$
$$+ 23821422836377458623944465583204209006664\, n^{20}$$
$$- 369623356219189533698996475700318296121020\, n^{19}$$
$$- 298878988216905150739838935426014962539412\, n^{18}$$
$$+ 4667996501363172027947082269243181591869628\, n^{17}$$
$$+ 3092064539849215601441788555160947528 1754764\, n^{16}$$

$$-\ 4871496634841982003560036946203580371856696 0\ n^{15}$$

$$-\ 257964348061911767931087943586517758164567888\ n^{14}$$

$$+\ 4113040052670775619144321001107945391061353 12\ n^{13}$$

$$+\ 16862977397540527171120192151587461684007527 36\ n^{12}$$

$$-\ 27350986122646178566252448727935165221541967 60\ n^{11}$$

$$-\ 83062388449508551726065791481976785508739672 40\ n^{10}$$

$$+\ 13826907573558591687222491158693276087388049240\ n^{9}$$

$$+\ 29115073931613853814571501443310002641224161400\ n^{8}$$

$$-\ 50586064684200076565468497744311642005530266720\ n^{7}$$

$$-\ 6602029026451859094204248382482374225451243952 8\ n^{6}$$

$$+\ 1243234677388779246957979746093914343845337926 52\ n^{5}$$

$$+\ 78768629562383357220196756827594170073871748626\ n^{4}$$

$$-\ 18031467821301399817819412254608697230307451926 5\ n^{3}$$

$$-\ 15989951671850958659019530620024145531014363211\ n^{2}$$

$$+\ 11414226661428343707762635720307970444805880444 9\ n$$

$$-\ 3804742220476114569254211906769323481601960148 3\ )\ n\ /\ 14490$$

$$s_{n,\,69} := (\ 30\ n^{66} + 990\ n^{65} + 10065\ n^{64} - 21120\ n^{63} - 884720\ n^{62} + 1790560\ n^{61} + 90957875\ n^{60}$$

$$-\ 183706310\ n^{59} - 9163896175\ n^{58} + 18511498660\ n^{57} + 873742090255\ n^{56}$$

$$-\ 1765995679170\ n^{55} - 78127199022435\ n^{54} + 158020393724040\ n^{53}$$

$$+\ 6525967903394955\ n^{52} - 13209956200513950\ n^{51} - 507781649453477275\ n^{50}$$

$$+\ 1028773255107468500\ n^{49} + 36703800690142510575\ n^{48} - 74436374635392489650\ n^{47}$$

$$-\ 2457543099613935202751\ n^{46} + 4989522573863262895152\ n^{45}$$

$$+\ 151952721981386979102947\ n^{44} - 308894966536637221101046\ n^{43}$$

$$-\ 8647339225545954721568655\ n^{42} + 17603573417628546664238356\ n^{41}$$

$$+\ 451277844586358265699304543\ n^{40} - 920159262590345078062847442\ n^{39}$$

$$-\ 21511290252403950629187193539\ n^{38} + 43942739767398246336437234520\ n^{37}$$

$$+\ 932510065695714301997186903099\ n^{36} - 1908962871158826850330811040718\ n^{35}$$

$$-\ 36585756959874337392012983774453\ n^{34} + 75080476790913694434356778589624\ n^{33}$$

$$+\ 1292166777977905958089679289447155\ n^{32}$$

$$-\ 2659414032746725610613715357483934\ n^{31}$$

$$-\ 40839242550610771639824188533509457\ n^{30}$$

$$+\ 84337899133968268888578553064502848\ n^{29}$$

$$+\ 1147278992487222640898726288447521811\ n^{28}$$

$$-\ 2378895884108413550686031129959546470\ n^{27}$$

$$-\ 28430148947141619119134029240638540799\ n^{26}$$

$$+\ 592391937783916517889540896112366280 68\ n^{25}$$

$$+ 616043276461676741987576166970922170063 \, n^{24}$$
$$- 1291325746701745135764106423553080968194 \, n^{23}$$
$$- 1155483508768817997967113267139381384 6515 \, n^{22}$$
$$+ 2440099592207810509510637176634070866 1224 \, n^{21}$$
$$+ 18537924727242701643323076281452031388 1067 \, n^{20}$$
$$- 39515949046693213796156789739538133642 3358 \, n^{19}$$
$$- 25078946320803341176256739089518397339 59451 \, n^{18}$$
$$+ 54109487546276003732129157152990608043 42260 \, n^{17}$$
$$+ 28114242332277521673838915675669695831 3564431 \, n^{16}$$
$$- 61639433419182643720891229228692977431 471122 \, n^{15}$$
$$- 25547590314424047865125620385905181023 0112127 \, n^{14}$$
$$+ 57259123970766360102340363694679659789 1695376 \, n^{13}$$
$$+ 18281328907585547466271031718800574869 18811075 \, n^{12}$$
$$- 42288570212247730942776099807069115717 29317526 \, n^{11}$$
$$- 98929015532672308423410787101927775576 10013183 \, n^{10}$$
$$+ 24014660127759234778959767401092466686 949343892 \, n^{9}$$
$$+ 38034107746875950882924791586759640011 826744799 \, n^{8}$$
$$- 10008287562151113654480935057461174671 0602833490 \, n^{7}$$
$$- 92552057153877235559225049876808554475 012595635 \, n^{6}$$
$$+ 28518698992926560766325945032822885566 0628024760 \, n^{5}$$
$$+ 10278729239144839221488734129352446946 2970907 15 \, n^{4}$$
$$- 49076157471216239209303413291527779458 5222206190 \, n^{3}$$
$$+ 55143676332275467583806471119172723212 513095680 \, n^{2}$$
$$+ 38047422204761145692542119067693234816 0196014830 \, n$$
$$- 1902371110238057284627105953384661740 80098007415 \,) \, (n+1)^2 \, n^2 \, / \, 2100$$

$$s_{n,70} := (2\,n+1)\,(n+1)\,(165\,n^{68} + 5610\,n^{67} + 59840\,n^{66} - 92565\,n^{65} - 5235065\,n^{64}$$
$$+ 7898880\,n^{63} + 553415280\,n^{62} - 834072360\,n^{61} - 57540790120\,n^{60} + 86728221360\,n^{59}$$
$$+ 5670403722720\,n^{58} - 8548969694760\,n^{57} - 524701875411560\,n^{56} + 791327297964720\,n^{55}$$
$$+ 45413881439108600\,n^{54} - 6851648580764 5260\,n^{53} - 3666428929060968972\,n^{52}$$
$$+ 5533901636495276088\,n^{51} + 275381745927297011904\,n^{50} - 415839569709193155900\,n^{49}$$
$$- 19189863965850497609580\,n^{48} + 2899271573363034 2992320\,n^{47}$$
$$+ 1237017786617573243788560\,n^{46} - 18700230377931750371790 00\,n^{45}$$
$$- 735303977091100320046954 80\,n^{44} + 1112306080825616355256327 20\,n^{43}$$
$$+ 4016488246343270393415588000\,n^{42} - 6080347673556186407886198360\,n^{41}$$
$$- 20085468143460082823392743 8440\,n^{40} + 304322195988679335554834256840\,n^{39}$$

$$+\ 9157707499635518547320355425260\ n^{38} - 13888722347447617488757950266310\ n^{37}$$

$$-\ 37896353726846697998820356259687 0\ n^{36} + 57538966707642427872668431902846 0\ n^{35}$$

$$+\ 1416267405399915745275465196339592 0\ n^{34}$$

$$-\ 2153170591453694831849532010460811 0\ n^{33}$$

$$-\ 475359076692003491970860817080822406\ n^{32}$$

$$+\ 7238044679952737121155388856735376 64\ n^{31}$$

$$+\ 142408401567735806610821169594191516 32\ n^{30}$$

$$-\ 2172316246915800784768094488196549628 0\ n^{29}$$

$$-\ 3781446571095276110698801610060838703 60\ n^{28}$$

$$+\ 578078566898870420528660713950108553680\ n^{27}$$

$$+\ 8829841948390935606717804139579657164960\ n^{26}$$

$$-\ 13533802206035838620341036566344540024280\ n^{25}$$

$$-\ 179674351385808838219039789137758034687320\ n^{24}$$

$$+\ 276278428181731176638730201989809322043120\ n^{23}$$

$$+\ 3152843964351648282727149068400284703395320\ n^{22}$$

$$-\ 4867405160618338012410088703595331716114540\ n^{21}$$

$$-\ 47125828146559696373341205671520886398710700\ n^{20}$$

$$+\ 73122444800148713566216852859078995456123320\ n^{19}$$

$$+\ 591271139934540474607852572389667440183125120\ n^{18}$$

$$-\ 923467932301885068694887285014040658002749340\ n^{17}$$

$$-\ 6117019252959346815614090459783182472423797580\ n^{16}$$

$$+\ 9637262845589962757768579332181794037637071040\ n^{15}$$

$$+\ 51032986644863846567851028082657186044482823280\ n^{14}$$

$$-\ 81368111390090751230660831790076676085542770440\ n^{13}$$

$$-\ 333599625989708183238408989199094587235628016712\ n^{12}$$

$$+\ 541083494679607650472943899693680218896213410288\ n^{11}$$

$$+\ 1643220000081917041926366969563071293802653221824\ n^{10}$$

$$-\ 2735371747462679388126022404191447050152086537880\ n^{9}$$

$$-\ 5759823751922775739464272938588535551569095686560\ n^{8}$$

$$+\ 10007421501615503303259420609978526852429686798780\ n^{7}$$

$$+\ 13060768345207908054736533943723622060860582385810\ n^{6}$$

$$-\ 24594863268619613733734511220574696517505716978105\ n^{5}$$

$$-\ 15582767350186627851519317863810295372963882189625\ n^{4}$$

$$+\ 35671582659589748644146232406002791318198681773490\ n^{3}$$

$$+\ 31632859201370871566486882051109747979048154344520\ n^{2}$$

$$-\ 22580720210000505057046148510667858127671572403525\ n$$

$$+ 7526906736668350190153828368892860425571908011175) \, n \, / \, 23430$$

$$
\begin{aligned}
s_{n,\,71} := (\, & 330 \, n^{68} + 11220 \, n^{67} + 117810 \, n^{66} - 246840 \, n^{65} - 10940820 \, n^{64} + 22128480 \, n^{63} \\
& + 1194275280 \, n^{62} - 2410679040 \, n^{61} - 128032096995 \, n^{60} + 258474873030 \, n^{59} \\
& + 13016380657335 \, n^{58} - 26291236187700 \, n^{57} - 1243690215703551 \, n^{56} \\
& + 251367166759 4802 \, n^{55} + 111263710476417147 \, n^{54} - 225041092620429096 \, n^{53} \\
& - 9295030226228685543 \, n^{52} + 1881510154507 7800182 \, n^{51} + 72326318424943159597 9 \, n^{50} \\
& - 146534147004394099 2140 \, n^{49} - 5227977779094793408 4715 \, n^{48} \\
& + 10602489705193980916 1570 \, n^{47} + 350045633700413762304 1655 \, n^{46} \\
& - 71069375710602150552448 80 \, n^{45} - 21643735131235673132826099 5 \, n^{44} \\
& + 439981640195773677117668 70 \, n^{43} + 1231703777057751215027523301 5 \, n^{42} \\
& - 2507405718135079797826223290 0 \, n^{41} - 64278806137001091759063821369 5 \, n^{40} \\
& + 1310650179921372633159538660290 \, n^{39} + 306401052543315072622509086347 15 \, n^{38} \\
& - 625908606885843871576613559297 20 \, n^{37} - 132824234522051706195857731379661 3 \, n^{36} \\
& + 2719075551129618511074815983522946 \, n^{35} \\
& + 5211177169509435095671639374496612 1 \, n^{34} \\
& - 1069426189413183204245076034734551 88 \, n^{33} \\
& - 1840527727959695052996607071822202 189 \, n^{32} \\
& + 3787999807486070842641772174711785956 6 \, n^{31} \\
& + 5817032258237211601165761247852832885 7 \, n^{30} \\
& - 12012864323960494044973294670417451728 0 \, n^{29} \\
& - 16341534494101124814110640583921307196 35 \, n^{28} \\
& + 33884355420598299032718610634884359565 50 \, n^{27} \\
& + 40495142222175944348351446555423793495 975 \, n^{26} \\
& - 84378719986411718599974541743360229485 00 \, n^{25} \\
& - 87747553281278345292882021092524611632 5455 \, n^{24} \\
& + 183932978561197862445761517602482825559 9410 \, n^{23} \\
& + 164583974252085589161252946535494831571 55595 \, n^{22} \\
& - 347561246360290964567082044831237945699 10600 \, n^{21} \\
& - 264049231585018641506265470843584343744 443495 \, n^{20} \\
& + 562854587806066379469239146170292482058 797590 \, n^{19} \\
& + 357217790146623371597474371638134804254 4370315 \, n^{18} \\
& - 770721039073853381141872657893298856714 7538220 \, n^{17} \\
& - 400451733063924761988752552865031210314 71078891 \, n^{16} \\
& + 877975570035234862091692371519392306300 89696002 \, n^{15} \\
& + 363893029593498442154111595790930946689 103935287 \, n^{14}
\end{aligned}
$$

$$- 8155836161905203705173924287338011240082975665576\, n^{13}$$
$$- 2603943494984169182323599243909257773747815648723\, n^{12}$$
$$+ 6023470606158858735164590916552316671503928864022\, n^{11}$$
$$+ 14091183839190188665307794462846223238468453875079\, n^{10}$$
$$- 34205838284539236065780179842244763148440836614180\, n^{9}$$
$$- 5417476374702175214202891149948674270033756723635\, n^{8}$$
$$+ 1425553657785827403498380028412182485491159710 86850\, n^{7}$$
$$+ 13182866978190103345681971462436912731968175 7328715\, n^{6}$$
$$- 4062127053423848072634774320899565031884794857 44280\, n^{5}$$
$$- 146407464546349118575415569455045645080981782189950\, n^{4}$$
$$+ 69902763443508304441430857100004779335044305012 4180\, n^{3}$$
$$- 78545174697535461522600503372009599143162656219790\, n^{2}$$
$$- 5419372850400121213691075642560285950641177376 84600\, n$$
$$+ 2709686425200060606845537821280142975320588684 42300\,) \,(n + 1)^2\, n^2\, /\, 23760$$

$s_{n,\,72} := (2\,n + 1)\,(n + 1)\,(10555545\, n^{70} + 369444075\, n^{69} + 4063884825\, n^{68} - 6280549275\, n^{67}$

$$- 375576846645\, n^{66} + 566505544605\, n^{65} + 42120772879185\, n^{64} - 6464412091080\, n^{63}$$
$$- 4655495981832696\, n^{62} + 7014976178794584\, n^{61} + 488655536457160072\, n^{60}$$
$$- 736490792775137400\, n^{59} - 48261265104409004648\, n^{58} + 72760143053001075672\, n^{57}$$
$$+ 4468164353281281068416\, n^{56} - 6738626601448422140460\, n^{55}$$
$$- 386778359952101603522036\, n^{54} + 583536853228876616353284\, n^{53}$$
$$+ 31227048664421286175711212\, n^{52} - 47132341423246367571743460\, n^{51}$$
$$- 2345450792123956868089999260\, n^{50} + 3541742358897558485920870620\, n^{49}$$
$$+ 163442104176881913842495896860\, n^{48} - 24693402744477165000670428 0600\, n^{47}$$
$$- 10535816784604267440452540023080\, n^{46} + 15927192190628786985682162174920\, n^{45}$$
$$+ 626266576635551421334157480407320\, n^{44} - 9473634610486415254940773 01698440\, n^{43}$$
$$- 34208878232724495182853794840962088\, n^{42}$$
$$+ 51786999079611063537027730912292352\, n^{41}$$
$$+ 1710701730165606578405423351684094916\, n^{40}$$
$$- 2591946094788215399376648892982288550\, n^{39}$$
$$- 77997216603490037187210210462149615594\, n^{38}$$
$$+ 118291797952629163480503640139715567666\, n^{37}$$
$$+ 3227674733682987850821751880094186920998\, n^{36}$$
$$- 4900657999500796357972879640211138165330\, n^{35}$$
$$- 120625075278431227242156073524776732073678\, n^{34}$$
$$+ 183387941917397239042220550107270667193182\, n^{33}$$

$$+ 40486862998660750450038935109314732524 0326\ n^{32}$$
$$- 616472342075860987627169430169335632145 7080\ n^{31}$$
$$- 1212908247023627154798559334037846224064 92040\ n^{30}$$
$$+ 1850185987639233781579197472565236117704 66600\ n^{29}$$
$$+ 32207002404890236420261015657066867846312 68600\ n^{28}$$
$$- 49235596601154971521181122221882919828321 36200\ n^{27}$$
$$- 75204749166766246373732105996558449181676 993240\ n^{26}$$
$$+ 115268903580207118136657215105931819763931 557960\ n^{25}$$
$$+ 153030650000875016334873787164616165107250 1382640\ n^{24}$$
$$- 235309420180322880409143541502220838649071 7852940\ n^{23}$$
$$- 2685312385962441929116919292970593969413300 5575188\ n^{22}$$
$$+ 4145623289033824333879950710207001373444448 67289252\ n^{21}$$
$$+ 4013759369367277509221169843689541824614003 39856716\ n^{20}$$
$$- 62279202185026074805257523010446628055932294 3429700\ n^{19}$$
$$- 50359222767780409412343233060915849256345564 56128700\ n^{18}$$
$$+ 78652794260921917858777725741896105287314961 55907900\ n^{17}$$
$$+ 520993355550358052296901664450355819536280291 67444700\ n^{16}$$
$$- 820816430455998037374741359546481781948077918 29121000\ n^{15}$$
$$- 4346536418534454910338370417902205464065099785 30173528\ n^{14}$$
$$+ 6930212843029681384194926306626549087071688637 09820792\ n^{13}$$
$$+ 28413052398133035632325168056000870980910850348 14711496\ n^{12}$$
$$- 46084685018714394140585215237314581014902119840 76977640\ n^{11}$$
$$- 139954881020844040713633659806837112608368852007 65651552\ n^{10}$$
$$+ 232974664040623258140743097328912959420043379318 6966148\ n^{9}$$
$$+ 490570616144641318557365100577155748932655349738 03709114\ n^{8}$$
$$- 852343256237273606906419199530190103108985193572 99046745\ n^{7}$$
$$- 111240021403296840706239351696568404470388850810 088661887\ n^{6}$$
$$+ 209477194916808941404679987521362111861032535893 782516203\ n^{5}$$
$$+ 132720168350077423098646850314034207124811006299 695465209\ n^{4}$$
$$- 303818849983520605350310269231732366617732777396 434455915\ n^{3}$$
$$- 269420591622319019795386763048593598303038168417 62244133\ n^{2}$$
$$+ 192322513735108155644463149073155223054322113960 860594157\ n$$
$$- 641075045783693852148210496910517410181073713202 86864719\ )\ n\ /\ 1541109570$$

$$s_{n,73} := (570570\ n^{70} + 19969950\ n^{69} + 216341125\ n^{68} - 452652200\ n^{67} - 2194792619\ n^{66}$$
$$+ 42842237438\ n^{65} + 2452142245553\ n^{64} - 4947126728544\ n^{63} - 279202265366024\ n^{62}$$
$$+ 563351657460592\ n^{61} + 30205639746347065\ n^{60} - 60974631150154722\ n^{59}$$

$$-\,3077396604887529925\,n^{58} + 6215767840925214572\,n^{57} + 29419078994360897 9421\,n^{56}$$

$$-\,59459734772814317 3414\,n^{55} - 263239041708274306 0077\,n^{54}$$

$$+\,5323937818189362929356 8\,n^{53} + 2199055402931049077221521\,n^{52}$$

$$-\,4451350184043991783736610\,n^{51} - 1711138202145923822771378 45\,n^{50}$$

$$+\,34667899061322875633801230 0\,n^{49} + 123686783250778481878574365 65\,n^{48}$$

$$-\,250840356407689252132052885430\,n^{47} - 8281603489381787146070595453 49\,n^{46}$$

$$+\,16814047335171263543461719761 28\,n^{45} + 51206138926003298346952586731273\,n^{44}$$

$$-\,10409368258552372304825134543867 4\,n^{43}$$

$$-\,291404400297058035315721790048923 7\,n^{42}$$

$$+\,5932181688526684429362687146417148\,n^{41}$$

$$+\,152074934232159113868935302756023341\,n^{40}$$

$$-\,3100820501528449121672332926584638 30\,n^{39}$$

$$-\,7249033203506427098779466062751571963\,n^{38}$$

$$+\,14808148457165699109726165418161607756\,n^{37}$$

$$+\,314244118520193686624217358733696902921\,n^{36}$$

$$-\,6432963854975530723581608828855555413598\,n^{35}$$

$$-\,12328938180245185210718897387177990903217\,n^{34}$$

$$+\,253011727459879234937959556572415372200 32\,n^{33}$$

$$+\,435443889921990199272309446016213245351293\,n^{32}$$

$$-\,896188952589968322038414847689668027922618\,n^{31}$$

$$-\,1376230912398526844249115977190488698080 6939\,n^{30}$$

$$+\,2842080720056050520702073439149944198953 6496\,n^{29}$$

$$+\,3866618535508330176968831343475103489387199257\,n^{28}$$

$$-\,8016578782172208591446834213417064207639350 10\,n^{27}$$

$$-\,9580601250633219162798960973723363149673937669\,n^{26}$$

$$+\,19962860379483659184742605368788432720111810348\,n^{25}$$

$$+\,20759880631961986041715055613711509144195892025 3\,n^{24}$$

$$-\,435160473018723380019043717643018615604029650854\,n^{23}$$

$$-\,389383353910138654893782883106895039074664295135 7\,n^{22}$$

$$+\,8222827551221496477894701379780919397097315553568\,n^{21}$$

$$+\,6247046582706221130727036749770815046588553484472 1\,n^{20}$$

$$-\,1331637592053459190924354363751972203288683852430 10\,n^{19}$$

$$-\,8451288276136554981864227054861877562159379827850 45\,n^{18}$$

$$+\,182342141443265691546528084734757273276074435081310 0\,n^{17}$$

$$+\,9474144709905356154833793741445883347476048315807285\,n^{16}$$

$$- 20771710834243369225132868330293394277128409824276670 \, n^{15}$$

$$- 8609215385126912748794653763645054287469651117098142 9 \, n^{14}$$

$$+ 192956018536781624201025943603140425177105863243390528 \, n^{13}$$

$$+ 6160577031129641354221039581049094895076306198007893 53 \, n^{12}$$

$$- 14250714247627098950452338598129594041923671029259692 34 \, n^{11}$$

$$- 333378292072609046271303334233066349941622803094135438 9 \, n^{10}$$

$$+ 809263726621489082047130054447428640302482316480867801 2 \, n^{9}$$

$$+ 1281701411146811500215919614228564508531582938905777588 5 \, n^{8}$$

$$- 3372666548915112082478969282904557657365648194292422978 2 \, n^{7}$$

$$- 3118887474582819967831340112309014493404441044686346730 0 \, n^{6}$$

$$+ 961044149808075201841649507522586644174530283665116438 2 \, n^{5}$$

$$+ 3463801979603436674175398884653499839113649714122359282 1 \, n^{4}$$

$$- 16538045457287625366492447276829586322401829711909835002 4 \, n^{3}$$

$$+ 1858272270806874161764118669309619059390177723926231029 3 \, n^{2}$$

$$+ 128215009156738770429642099382103482036214742640573729438 \, n$$

$$- 6410750457836938521482104969105174101810737132028686471 9 \, \big) \, (n+1)^2 \, n^2 \, /$$

$$42222180$$

$$s_{n,74} := (2 n + 1)(n + 1)(57057 \, n^{72} + 2054052 \, n^{71} + 23279256 \, n^{70} - 35945910 \, n^{69}$$

$$- 2269385118 \, n^{68} + 3422050632 \, n^{67} + 269548565286 \, n^{66} - 406033873245 \, n^{65}$$

$$- 31612780552353 \, n^{64} + 47622187765152 \, n^{63} + 3527427815847936 \, n^{62}$$

$$- 5314952817654480 \, n^{61} - 371069883186381808 \, n^{60} + 559262301188399952 \, n^{59}$$

$$+ 36667739320208772776 \, n^{58} - 55281240130907359140 \, n^{57}$$

$$- 3395250270810618825748 \, n^{56} + 5120516026281381918192 \, n^{55}$$

$$+ 293913136474635025055136 \, n^{54} - 443429962725093228541800 \, n^{53}$$

$$- 23729647785024584501631624 \, n^{52} + 35816186658899423366718336 \, n^{51}$$

$$+ 1782327688162141456953991608 \, n^{50} - 2691399625572661897114346580 \, n^{49}$$

$$- 124201085272580712689174281764 \, n^{48} + 1876473277216573999823 18595936 \, n^{47}$$

$$+ 8006260374130238798366701768128 \, n^{46} - 12103214225056186897541211950160 \, n^{45}$$

$$- 475905526277506025367781320932544 \, n^{44} + 719909896528787131500442587373896 \, n^{43}$$

$$+ 2599563018605255604451713663746682 8 \, n^{42}$$

$$- 3935340022734322763252592624988719 0 \, n^{41}$$

$$- 12999774303764226823720157275449159 34 \, n^{40}$$

$$+ 19696428456783056373742865544423174 96 \, n^{39}$$

$$+ 59270777293808312228702057941975179 248 \, n^{38}$$

$$- 89890987363551621161740230190183927 620 \, n^{37}$$

$$- 24527386831246724299123179518040587 80404 \, n^{36}$$

$$
\begin{aligned}
&+ 37240535183687844554493470428011801344416\, n^{35} \\
&+ 9166406553030392866337317758944046099 0228\, n^{34} \\
&- 139358125054640285222784439905561281552550\, n^{33} \\
&- 3076632660735309012346584219208094692016478\, n^{32} \\
&+ 46846280536302836611312685487649226788 00992\, n^{31} \\
&+ 92169974428437210591835855760592455688139776\, n^{30} \\
&- 140597275669470957718319417915271144871610160\, n^{29} \\
&- 244743870392597927455380085091021350291115 6688\, n^{28} \\
&+ 374145669372370439068986098532295582680254 0112\, n^{27} \\
&+ 571487565082620777815447311632384127696591 65576\, n^{26} \\
&- 875938631092549688676620272375190706789001 8420\, n^{25} \\
&- 1162893494373363432411158688674749597036534 369828\, n^{24} \\
&+ 1788137173114672633050569046630883944088746 563952\, n^{23} \\
&+ 20405927204635818080866308202251360543185204 239136\, n^{22} \\
&- 31502959393511063437824746826692482786822179 640680\, n^{21} \\
&- 30500913761986168075239925475419805818085581 6928008\, n^{20} \\
&+ 473265186126548052847511255544643328664694815 212352\, n^{19} \\
&+ 382684204360475679238276586215114869922967105 5162776\, n^{18} \\
&- 597689565847040921499790442099904471317685399 0350340\, n^{17} \\
&- 39590747590616702646080768624109940009841453813120148\, n^{16} \\
&+ 62374569215160258576620105146664378371350607714855392\, n^{15} \\
&+ 330297161002828958900555509389737706648106213251301536\, n^{14} \\
&- 526633026111823567639143316657938749157834623734380000\, n^{13} \\
&- 2159133074902952186682261331633797138733966979520332520\, n^{12} \\
&+ 3502016125410340063842963655779665082679867781147688780\, n^{11} \\
&+ 10635295651165854178570202810360990017625593945577350890\, n^{10} \\
&- 17703951539453951299776786043431317567778324808939870725\, n^{9} \\
&- 37278896615944495361711006866553282063537272220961245025\, n^{8} \\
&+ 64770320693643718692454903321545581879150707359118 02900\, n^{7} \\
&+ 84532279773576758307344633458720901634616069796649 63400\, n^{6} \\
&- 159183580012858373092329146679580926184789945837453346550\, n^{5} \\
&- 100855234128301348602359693793165231657758370295485458590\, n^{4} \\
&+ 2308746411988812094497041140295383105790325283619548 61160\, n^{3} \\
&+ 2047350993059420601978017950165227043266622904994794230\, n^{2} \\
&- 1461475854953319137545223262672475609385156082884696 21925\, n \\
&+ 48715861831777304584840775422415853646171869429489873975\, )\, n\, /\, 8558550
\end{aligned}
$$

$$
\begin{aligned}
s_{n,75} := \ & ( \, 6006\, n^{72} + 216216\, n^{71} + 2414412\, n^{70} - 5045040\, n^{69} - 249175927\, n^{68} \\
& + 503396894\, n^{67} + 30504890559\, n^{66} - 61513178012\, n^{65} - 3682426426250\, n^{64} \\
& + 7426366030512\, n^{63} + 423139355846376\, n^{62} - 853705077723264\, n^{61} \\
& - 45875554013517293\, n^{60} + 92604813104757850\, n^{59} + 4676291454306992793\, n^{58} \\
& - 9445187721718743436\, n^{57} - 447097684190452902997\, n^{56} + 90364055610262454 9430\, n^{55} \\
& + 40004821451799243599937\, n^{54} - 80913283459701111749304\, n^{53} \\
& - 3342154848943097996089437\, n^{52} + 6765222981345897103928178\, n^{51} \\
& + 260061656183729570418344681\, n^{50} - 526888535348805037940617540\, n^{49} \\
& - 18798133949618217157055387421\, n^{48} + 38123156434585239352051392382\, n^{47} \\
& + 1258652752772677265172247406041\, n^{46} - 2555428661979937696965 46204464\, n^{45} \\
& - 77823998333483662111597452992733\, n^{44} + 158203425328947263992891452189930\, n^{43} \\
& + 4428815802928977608749103219580073\, n^{42} \\
& - 9015835031186902481491097891350076\, n^{41} \\
& - 231126184944993580972775279909791431\, n^{40} \\
& + 4712682049211740644270416577 10932938\, n^{39} \\
& + 11017209365033609136609640518748153655\, n^{38} \\
& - 22505686934988392337646322695207240248\, n^{37} \\
& - 477593790599675965569432448752003089601\, n^{36} \\
& + 977693268134340323476511220199213419450\, n^{35} \\
& + 18737739141971999346319786677351513472101\, n^{34} \\
& - 38453171552078339016116084574902240363652\, n^{33} \\
& - 661795355049396982613881904254595118044839\, n^{32} \\
& + 1362043881650872304243879893084092476453330\, n^{31} \\
& + 20916201751398236274681290700236008367569479\, n^{30} \\
& - 431944473844473448536064612935561092115922 88\, n^{29} \\
& - 587589714536495323026308880588297880380104 13\, n^{28} \\
& + 12183738764574379909062242242471215685287613114\, n^{27} \\
& + 1456076788079085793596122635985112901970475 6825\, n^{26} \\
& - 303399096380391538628286769421734737246971 26764\, n^{25} \\
& - 315512351685602506445126547889621049071321878917\, n^{24} \\
& + 66136461300924416675308177272141557186734084 8598\, n^{23} \\
& + 591791734631977426962714221101919349327511 9155921\, n^{22} \\
& - 1249719930564879270600736619475980255841757 9196440\, n^{21} \\
& - 94943723104292311147180150119203289072938788852861\, n^{20} \\
& + 202384645514233415000367666433166380704295156902162\, n^{19}
\end{aligned}
$$

$$+ 1284441797161152037974747769285966768804815882145737\, n^{18}$$
$$- 2771268239836537490949863205005099918313926921193636\, n^{17}$$
$$- 1439897334009605228502951828908892562417372174818 4173\, n^{16}$$
$$+ 315692149200286420610088997831829511666613704175 61982\, n^{15}$$
$$+ 1308443839579325439152096221883582776020560466539 60009\, n^{14}$$
$$- 2932579828358937298914281441598995063707734637254 82000\, n^{13}$$
$$- 936295435070781613477702609831477469433842023451121885\, n^{12}$$
$$+ 216584885297745695684683336382285444523845751062772 5770\, n^{11}$$
$$+ 506674247301866094606604595545067662025887678936065 1945\, n^{10}$$
$$- 1229933379901477884897892527472420768575621108934902 9660\, n^{9}$$
$$- 1947952560801743339616891842368210351967442958505084 3070\, n^{8}$$
$$+ 5125838501504964564131676212208841472510507025945071 5800\, n^{7}$$
$$+ 4740140558579894706438444609543633162985275064459104 7660\, n^{6}$$
$$- 1460611961866475397700856543129610779848105715486328 11120\, n^{5}$$
$$- 526434774714139378496993196091147515134440931091584 22135\, n^{4}$$
$$+ 251348151129475415469484293531190581011698757766949 655390\, n^{3}$$
$$- 282423519011830985650605959207635832135056400244950 79745\, n^{2}$$
$$- 19486344732710921833936310168966341458468747717195 94955900\, n$$
$$+ 974317236635546091696815508448317072923473885897974 7950)\,(n+1)^2\,n^2 \,/\, 456456$$

$$s_{n,76} := (2\,n+1)\,(n+1)\,(195\,n^{74} + 7215\,n^{73} + 84175\,n^{72} - 129870\,n^{71} - 8643570\,n^{70}$$
$$+ 13030290\,n^{69} + 1085568900\,n^{68} - 1634868495\,n^{67} - 134863956045\,n^{66}$$
$$+ 203113368315\,n^{65} + 15968346309175\,n^{64} - 24054076147920\,n^{63} - 1785742474821840\,n^{62}$$
$$+ 2690640750306720\,n^{61} + 18793161163656840\,n^{60} - 283275062120638620\,n^{59}$$
$$- 1857528399476084 5780\,n^{58} + 2800456352320158 7980\,n^{57}$$
$$+ 172003524077249498 1420\,n^{56} - 2594055142920343 266120\,n^{55}$$
$$- 14889773290647507 3643640\,n^{54} + 22464362693117278 2098520\,n^{53}$$
$$+ 120215719395782177559 51840\,n^{52} - 181446797228329130249 77020\,n^{51}$$
$$- 902937604316299063044 202740\,n^{50} + 136347874633586505107 8792620\,n^{49}$$
$$+ 629210022361950150293791 47900\,n^{48} - 950632427274604550696081 18160\,n^{47}$$
$$- 405601884644888970173774158 0160\,n^{46} + 613155989103706478014141642 9320\,n^{45}$$
$$+ 241096555504402277181895505240 100\,n^{44} - 364710613202121948162913966 074810\,n^{43}$$
$$- 131695400915289297900692641631 14030\,n^{42}$$
$$+ 199366654438944556591853532277 08450\,n^{41}$$
$$+ 658576259666606716551079192603 961090\,n^{40}$$
$$- 997832722221857302656211465519 795860\,n^{39}$$

$$- 30026926552879733325517774448246031660 \, n^{38}$$
$$+ 45539306190430528639604767405128945420 \, n^{37}$$
$$+ 124257193268573022932536306105219480420 \, n^{36}$$
$$- 188662755212381060830784697528085679010 \, n^{35}$$
$$- 46437555808488563983056538990051346679570 \, n^{34}$$
$$+ 705996509033903650500020083384106285331330 \, n^{33}$$
$$+ 1558640314090970408094220295453326376486130 \, n^{32}$$
$$- 23732602965881507946663314473491948789957660 \, n^{31}$$
$$- 466938545268596057134749245278927192101756000 \, n^{30}$$
$$+ 7122741193858348396754555251551367625476128000 \, n^{29}$$
$$+ 1239886931869202850675367998345896055081169480 \, n^{28}$$
$$- 18954441037730960179968247737660092074913486000 \, n^{27}$$
$$- 289518982655240915565168154086905352240902769000 \, n^{26}$$
$$+ 443755569450172685343773635499924103296509982780 \, n^{25}$$
$$+ 5891287264994159896781701033963483337726949890780 \, n^{24}$$
$$- 9058808744742103271891419728444845523072973275600 \, n^{23}$$
$$- 10337763488465404733012461222018174080752842630320 \, n^{22}$$
$$+ 159595856669935212263112790169724968388266575058328 \, n^{21}$$
$$+ 154519434226865379279261121723154351897755411185024 \, n^{20}$$
$$- 2397589441752656750504480776695940120407965993067000 \, n^{19}$$
$$- 19387008273514836090823509360016098116968309321593320 \, n^{18}$$
$$+ 302793071311485825114875044283721172356564469789234800 \, n^{17}$$
$$+ 200569070358315828162651423650127839403953518062016440 \, n^{16}$$
$$- 315993259103048033499720887689377817723758500582486400 \, n^{15}$$
$$- 16733049653255219363984344920011710233712403379339343200 \, n^{14}$$
$$+ 26679540775398069213475121818464454439187397571921446800 \, n^{13}$$
$$+ 109382959395243726633146324286724307479774268451203861400 \, n^{12}$$
$$- 177414209480564624556457047339318688439255101462766515500 \, n^{11}$$
$$- 538790371881165890867670726089765703214234770000456973700 \, n^{10}$$
$$+ 896892662562031148579734612804307899040979705732068718300 \, n^{9}$$
$$+ 18885709650039181532653482601095814021453855895485262319000 \, n^{8}$$
$$- 3281302778768928041878896965665260527385682371888237070000 \, n^{7}$$
$$- 4282455321513566268281011311305017617887409482375277464200 \, n^{6}$$
$$+ 8064334371663795804515461815240789453200398342157328049800 \, n^{5}$$
$$+ 5109385849202296374660479773237854175444813742246139832400 \, n^{4}$$

$$- 11696245959635342464248450567477175989767419784447873735 \, n^3$$

$$- 103720012974048079466282645816609835759812002452297543557 \, n^2$$

$$+ 74039231744283924241184649709877355312808899290084000403 \, n$$

$$- 2467974391476130808039488323662578510426963309669466801 \,) \, n \, / \, 30030$$

$$s_{n,\,77} := (\, 30 \, n^{74} + 1110 \, n^{73} + 12765 \, n^{72} - 26640 \, n^{71} - 1385910 \, n^{70} + 2798460 \, n^{69}$$

$$+ 179254415 \, n^{68} - 361307290 \, n^{67} - 22903521150 \, n^{66} + 46168349590 \, n^{65}$$

$$+ 2790375830845 \, n^{64} - 5626920011280 \, n^{63} - 321325520515400 \, n^{62} + 648277961042080 \, n^{61}$$

$$+ 34855316670626115 \, n^{60} - 70358911302294310 \, n^{59} - 3553405882587323035 \, n^{58}$$

$$+ 7177170676476940380 \, n^{57} + 339749953341616830375 \, n^{56} - 686677077359710601130 \, n^{55}$$

$$- 30399940212512404121655 \, n^{54} + 61486557502384518844440 \, n^{53}$$

$$+ 2539731506676687521920875 \, n^{52} - 5140949570855759562686190 \, n^{51}$$

$$- 197623129682911684775146795 \, n^{50} + 400387208936679129112979780 \, n^{49}$$

$$+ 14284868301088167982197318495 \, n^{48} - 28970123811113015093507616770 \, n^{47}$$

$$- 956461390940138028233773488983 \, n^{46} + 1941892905691389071561054594736 \, n^{45}$$

$$+ 59139147058990717260421354073011 \, n^{44} - 120220187023672823592403762740758 \, n^{43}$$

$$- 3365496442641187219742401835445285 \, n^{42}$$

$$+ 6851213072306047263077207433631328 \, n^{41}$$

$$+ 175634839694397438642175966158355379 \, n^{40}$$

$$- 358120892461100924547429139750342086 \, n^{39}$$

$$- 8372075199451971567040618272435884037 \, n^{38}$$

$$+ 17102271291365044058628665684622110160 \, n^{37}$$

$$+ 362927761225202730136386684383921225667 \, n^{36}$$

$$- 742957793741770504331402034452464561494 \, n^{35}$$

$$- 14238974314813816656746575931134356653889 \, n^{34}$$

$$+ 292209064233694038178245538967211177869272 \, n^{33}$$

$$+ 5029041652692793876057443391412538628777995 \, n^{32}$$

$$- 10350292369619281790293132321792289036252622 \, n^{31}$$

$$- 158944073906500159829920513634380366466366121 \, n^{30}$$

$$+ 3282384401826196014501341595905530219635750442 \, n^{29}$$

$$+ 44651464029674419835669604376603501095048484963 \, n^{28}$$

$$- 925853124611750356858405503491125324097327430 \, n^{27}$$

$$- 1106485677317289714790963246485677979382518788333 \, n^{26}$$

$$+ 23055566670957544652677670433204684911747703196 \, n^{25}$$

$$+ 239760636949220510092085212573542790444078949991 \, n^{24}$$

$$- 5025768405693985648368480955802902657999056033178 \, n^{23}$$

$$- 4497077926699813639936900112076209983298059809319 \, n^{22}$$

$$+ 9496732693969025844710648319732710232396025221816 \, n^{21}$$

$$+ 721485780325295096042277204309312374065748901 86187 \, n^{20}$$

$$- 153793888759028045053166089181595185045545805594190 \, n^{19}$$

$$- 97605872406044583696247326478288104212046767 5648827 \, n^{18}$$

$$+ 210591133687991971897811261874735726928648115 6891844 \, n^{17}$$

$$+ 109419076653974834320745035640488444754847473 94041439 \, n^{16}$$

$$- 23989726667674886583127119746845046220255975944974722 \, n^{15}$$

$$- 99429808916776781770266271517703674051966217489202375 \, n^{14}$$

$$+ 222849344501228450123659662782252394324188410923379472 \, n^{13}$$

$$+ 71149921290216903730819820511131239056633119038066131 \, n^{12}$$

$$- 164584777030556652474005607300487717545685464 8999511734 \, n^{11}$$

$$- 385026258432867282542051293014768588650888343 2201941088 \, n^{10}$$

$$+ 934637293896291217558108193330024894847462151 3403393910 \, n^{9}$$

$$+ 148026644354666677812694972262587919805570180 36859161493 \, n^{8}$$

$$- 389517018098962477381200763858178329095886575 87121716896 \, n^{7}$$

$$- 360207489019775290801914051678233181124505313 83226726606 \, n^{6}$$

$$+ 110993199613851305898502886721464446913448972 0353575170108 \, n^{5}$$

$$+ 40004245863393021492583134331592323037996688390493784015 \, n^{4}$$

$$- 19100169134063734888366915538464911521048309713 4562738138 \, n^{3}$$

$$+ 21461613926034750200649927982447202292432649277 197368666 \, n^{2}$$

$$+ 148078463488567848482369299419754710625617798580168000806 \, n$$

$$- 74039231744283924241184649709877355312808899290084000403 \,) \, (n + 1)^2 \, n^2 \, / \, 2340$$

$$s_{n,78} := (2\,n + 1)\,(n + 1)\,(\,105\,n^{76} + 3990\,n^{75} + 47880\,n^{74} - 73815\,n^{73} - 5171971\,n^{72}$$

$$+ 7794864\,n^{71} + 685800402\,n^{70} - 1032598035\,n^{69} - 90104745171\,n^{68} + 135673416774\,n^{67}$$

$$+ 11301501536612\,n^{66} - 17020089013305\,n^{65} - 1341115457121741\,n^{64}$$

$$+ 2020183230189264\,n^{63} + 150057917736536712\,n^{62} - 226096968219899700\,n^{61}$$

$$- 15796005878381838676\,n^{60} + 23807057301682707864\,n^{59}$$

$$+ 1561159888404919428192\,n^{58} - 2353643361258220496220\,n^{57}$$

$$- 144561556096224326401356\,n^{56} + 218019155824965599850144\,n^{55}$$

$$+ 12514238788629313921187912\,n^{54} - 18880367760856453681706940\,n^{53}$$

$$- 1010363926087797463323168636\,n^{52} + 1524986073012124421825606424\,n^{51}$$

$$+ 75888220144762257846785595312\,n^{50} - 114594823253649448981091196180\,n^{49}$$

$$- 5288253580628039969300197635028\,n^{48} + 7989677782568884678440842050632\,n^{47}$$

$$+ 340891842036120352441950656606316\,n^{46} - 515332601945464971002146405934790\,n^{45}$$

$$- 20263182244750254302331944738445078\,n^{44}$$

$+ 3065243966809113938998990310635012 \, n^{43}$

$+ 110684613691124679301565807014955369 6 \, n^{42}$

$- 1675595425200919246492986600379648050 \, n^{41}$

$- 5535064883772303656809154941282334773 8 \, n^{40}$

$+ 8386377096918501447538381741942484563 2 \, n^{39}$

$+ 2523640721873497511420868570611950597 596 \, n^{38}$

$- 3827392968294838774368994764627638319 210 \, n^{37}$

$- 1044331035240445121224385485942288331 63818 \, n^{36}$

$+ 1585633517702141875708423202736570689 05332 \, n^{35}$

$+ 3902887376837322884938870924494327592 485336 \, n^{34}$

$- 5933612741141091421193727546878319923 180670 \, n^{33}$

$- 1309973620012133086428938187023060845 78306998 \, n^{32}$

$+ 1994628493723905086749375918268982860 29050832 \, n^{31}$

$+ 3924428047567290751506663943570673657 793291176 \, n^{30}$

$- 5986373496037131381597464711269459630 104462180 \, n^{29}$

$- 1042074401555522785580050434165644775 99568256388 \, n^{28}$

$+ 1593043469813469835278062974804814462 14404615672 \, n^{27}$

$+ 2433288978492537980262993063405279817 788233931936 \, n^{26}$

$- 3729585641229480462158392743848160449 789553205740 \, n^{25}$

$- 4951386689595450780640754877031513981 9819917015196 \, n^{24}$

$+ 7613559316454650194069051952739678995 4624652125664 \, n^{23}$

$+ 8688468620622287900422400865308653862 0014802193832 \, n^{22}$

$- 1341338089675616436033681272743328202 907334529353580 \, n^{21}$

$- 1298672829045079293353388315306964915 336546462458156 \, n^{20}$

$+ 2015076148051399761831766536597613783 1501839958364024 \, n^{19}$

$+ 1629398981898966143080224747615862111 266036510506 70992 \, n^{18}$

$- 2544852280251019202711925448253673856 0565656396555188500 \, n^{17}$

$- 1685701240911536914693502475889629971 244221503265343276 \, n^{16}$

$+ 2655794475379856332175849986247128649 66916045317560916 4 \, n^{15}$

$+ 1406344583157070699288426879723402696 2735937102392180042 \, n^{14}$

$- 2242306598504598865541432818897460476 938485880176074645 \, n^{13}$

$- 9193191655011014288627437743810539630 5010625230034739021 \, n^{12}$

$+ 1491094078176882086571187302516453968 41985180785140145854 \, n^{11}$

$+ 4528313347859186075642731292638737185 44202534906522402592 \, n^{10}$

$- 7538017060877220156749690590216332762 37296392752353676815 \, n^9$

$$- 1587266876976183384663529320390382678959882504320088976731\, n^8$$
$$+ 275780116850813608483277851009639065655712195285631030350 4\, n^7$$
$$+ 3599229051980460433293531280689714387231178074879624637922\, n^6$$
$$- 6777744162224758692356686176082766909125328088747592108635\, n^5$$
$$- 4294230436881045758038893690203425681003802740075079938011\, n^4$$
$$+ 9830217736433947983236683623346521976068368154486415961334\, n^3$$
$$+ 8717244102760256218728958143799134095658710355618661 81092\, n^2$$
$$- 6222695483631012424427685533243131102382990630586007252305\, n$$
$$+ 2074231827877004141475895177747710367460996878620024174 35\, )\, n \,/\, 16590$$

$$s_{n,\,79} := (42\, n^{76} + 1596\, n^{75} + 18886\, n^{74} - 39368\, n^{73} - 2154362\, n^{72} + 4348092\, n^{71}$$
$$+ 293958378\, n^{70} - 592264848\, n^{69} - 39691980692\, n^{68} + 79976226232\, n^{67}$$
$$+ 5118578060428\, n^{66} - 10317132347088\, n^{65} - 624952473306672\, n^{64}$$
$$+ 1260222078960432\, n^{63} + 72003967638178408\, n^{62} - 145268157355317248\, n^{61}$$
$$- 7811528838849253847\, n^{60} + 15768325835053824942\, n^{59} + 79638970112654938236 3\, n^{58}$$
$$- 1608547728088152589668\, n^{57} - 7614539621503514949 9875\, n^{56}$$
$$+ 153899340158158451589418\, n^{55} + 6813304919660732126465839\, n^{54}$$
$$- 13780509179479622704521096\, n^{53} - 569210777037711479223215127\, n^{52}$$
$$+ 1152202063254902581150951350\, n^{51} + 44291779860149563056395339467\, n^{50}$$
$$- 89735761783554028693941630284\, n^{49} - 3201559767330661821882937698371\, n^{48}$$
$$+ 6492855296444877672459817027026\, n^{47} + 214364477603023088540938868997759\, n^{46}$$
$$- 4352218105024910547543375550 22544\, n^{45}$$
$$- 132544110076983020341083897716487 41\, n^{44}$$
$$+ 2694404382589909512297111709832002 6\, n^{43}$$
$$+ 754283335732261041717203083050684289\, n^{42}$$
$$- 15355107152904211785573772831996886 04\, n^{41}$$
$$- 393637120193832674419284123554280656 01\, n^{40}$$
$$+ 802629347540569560624142019940558198 06\, n^{39}$$
$$+ 1876370074060536795202736977536593524789\, n^{38}$$
$$- 3833003082875130546467888157067242869384\, n^{37}$$
$$- 81340262000750292924892245828093560061089\, n^{36}$$
$$+ 166513527084375716396252379813254362991562\, n^{35}$$
$$+ 3191273925920943626645158122502804858719165\, n^{34}$$
$$- 6549061378926262969686568624818864080429892\, n^{33}$$
$$- 112712117767629559315907553815421877290601901\, n^{32}$$
$$+ 231973296914185381601501676255662618661633694\, n^{31}$$
$$+ 3562293656292453958651843892560443035965182913\, n^{30}$$

$$- 7356560609499093298905189461376548690591999520\, n^{29}$$
$$- 10007395881311859761548937348149067587060948 7943\, n^{28}$$
$$+ 207504478235736288529883936424357900431810975406\, n^{27}$$
$$+ 2479882899820880491249801487745009273204381379531\, n^{26}$$
$$- 5167270277877497271029486911914376446840573734468\, n^{25}$$
$$- 53735743336697268245777846354641674191649096605267\, n^{24}$$
$$+ 11263875695127203376258517962119772483013876694500 2\, n^{23}$$
$$+ 1007896159724702479551445657646429798244653805022623\, n^{22}$$
$$- 2128431076400676992865476494914057321319446376990248\, n^{21}$$
$$- 16170116665500819417238216744962446073494351524897607\, n^{20}$$
$$+ 34468664407402315827341909984838949468308149426785462\, n^{19}$$
$$+ 218756680600430317607982805842943591932767051590581883\, n^{18}$$
$$- 471982025608262951043307521670726133333842252607949228\, n^{17}$$
$$- 2452327243550690796235415851696512035808129382491148275\, n^{16}$$
$$+ 5376636512709644543514139225063750204950101017590245778\, n^{15}$$
$$+ 22284453194459798213742768581408842651005824290462518319\, n^{14}$$
$$- 49945542901629240970999676387881435506961749598515282416\, n^{13}$$
$$- 159462952614988861017650886451546551591183096091637475422\, n^{12}$$
$$+ 368871448131606963006301449290974538689327941781790233260\, n^{11}$$
$$+ 862930315180108304392051498549115452851784340457691870302\, n^{10}$$
$$- 2094732078491823571790404446389205444392896622697173973864\, n^{9}$$
$$- 3317609541435946772528166190546426521513526631050966030994\, n^{8}$$
$$+ 8729951161363717116846736827482058487419949884799106035852\, n^{7}$$
$$+ 8073058790723145717109682569627458450953210067865699894690\, n^{6}$$
$$- 24876068742810008551066101966736975389326370020530505825232\, n^{5}$$
$$- 8965849922014942934686107892084383076605293369831311372236\, n^{4}$$
$$+ 42807768586839894420438317750905741542536956760193128569704\, n^{3}$$
$$- 4810029670403914078411997453471187831580503365200544945372\, n^{2}$$
$$- 3318770924603206626361432284396336587937595002979203 8678960\, n$$
$$+ 16593854623016033131807161421981682939687975014896019339480 \big)\,(n+1)^2\, n^2\, /$$
$$3360$$

$$s_{n,\,80} := (2\,n+1)\,(n+1)\,\big(268345\, n^{78} + 10465455\, n^{77} + 129073945\, n^{76} - 198843645\, n^{75}$$
$$- 14648148515\, n^{74} + 22071644595\, n^{73} + 2047758137425\, n^{72} - 3082673028435\, n^{71}$$
$$- 284103756821885\, n^{70} + 427696971747045\, n^{69} + 37686658734599875\, n^{68}$$
$$- 56743836587773335\, n^{67} - 4737447622351813985\, n^{66} + 7134543351821607645\, n^{65}$$
$$+ 562482325759263396025\, n^{64} - 847290760314805897860\, n^{63}$$

$$- 6294472563919462484870\,0\,n^{62} + 9484073383894934022198\,0\,n^{61}$$

$$+ 662616164272315069886938\,0\,n^{60} - 998666283100420071841506\,0\,n^{59}$$

$$- 654885963735728477574894380\,n^{58} + 98732227701909481672154910\,0\,n^{57}$$

$$+ 6064178807492272727577895585540\,n^{56} - 914563432508936877520415286\,0\,n^{55}$$

$$- 5249571342106531499425763865380\,n^{54} + 792008518478524406852624787450\,0\,n^{53}$$

$$+ 42383546230981037873700436755382\,0\,n^{52} - 639713236057108190139769675267980\,n^{51}$$

$$- 3183419266199314592487733979640698\,0\,n^{50}$$

$$+ 480711456110182729823858945322444\,60\,n^{49}$$

$$+ 221835858489350685598305187752386920\,0\,n^{48}$$

$$- 335157345014576942046577076355192603\,0\,n^{47}$$

$$- 143000015836060443213122143760893267730\,n^{46}$$

$$+ 216175810479163549529916101023115864610\,n^{45}$$

$$+ 850016053767897105389041060756982231995\,0\,n^{44}$$

$$- 128583287117580383556005739618662914122\,30\,n^{43}$$

$$- 464308603739925804627171591488837585400890\,n^{42}$$

$$+ 702892069965767726118557674214189523807450\,n^{41}$$

$$+ 232189295527807475610012041303977207247502\,70\,n^{40}$$

$$- 351798403641540052045610850327036758490291\,30\,n^{39}$$

$$- 10586368429019834833662466365074612029308772\,70\,n^{38}$$

$$+ 16055451845350522276516504972775436423208304\,70\,n^{37}$$

$$+ 438084272658583087316827076140091623226604859\,30\,n^{36}$$

$$- 665154134910549892113498866696525153051511441\,30\,n^{35}$$

$$- 1637214178315843452130047134539400182115970845870\,n^{34}$$

$$+ 2489078974219292672800746451439265308265318408\,70\,n^{33}$$

$$+ 549518132814220239598655741386455968011838019070\,30\,n^{32}$$

$$- 8367225940924268227619873403054035846718896878098\,0\,n^{31}$$

$$- 16462502298657120477841505112579860868095358263598\,20\,n^{30}$$

$$+ 25112114745031894128143251339022493094478982239302\,20\,n^{29}$$

$$+ 437137642047328177549859749783334037999139126677185\,80\,n^{28}$$

$$- 668262520443508213388861250344512303545948181135429\,80\,n^{27}$$

$$- 102073537637063585082038639926062589221895806628090310\,0\,n^{26}$$

$$+ 156451619057812918690002266140816445350573450847812614\,0\,n^{25}$$

$$+ 207704699311622274470502498412320299150002146498495851\,40\,n^{24}$$

$$- 319379629920324057640253860925521270992531892290134407\,80\,n^{23}$$

$$- 3644707786844788729160650101303407640267402460216395373\,32\,n^{22}$$

$$+\ 5626751495227345122561102082417872095897369636469660263 88\ n^{21}$$

$$+\ 5447775873126584677331724750603253828949654435939946783084\ n^{20}$$

$$-\ 8453001384451244272125642230025774348219350135733403187820\ n^{19}$$

$$-\ 6835132192465455472633863485123582756700360945456430571 9292\ n^{18}$$

$$+\ 10675348357920745422557077339186662852461508924971316017 2848\ n^{17}$$

$$+\ 7071313377896693897099995356151158105590121377349392175 94374\ n^{16}$$

$$-\ 11140737484741078116777846901186070301008257512272654064 77985\ n^{15}$$

$$-\ 58994458943588253082656109077674757653809801631394431005 86447\ n^{14}$$

$$+\ 94062057157752918682373087067105171631218831203227973541 18663\ n^{13}$$

$$+\ 385643301184822185956825440175273516780299952003498226897 07329\ n^{12}$$

$$-\ 625495980356109738276424703796462860986059343606861327116 20325\ n^{11}$$

$$-\ 1889957282932944394132513842277916828746732419860926801708 560827\ n^{10}$$

$$+\ 3162107234172220781125919986066983861694015969717332689186 51403\ n^{9}$$

$$+\ 6658393094249975343522615168556499781933575053765765713720 34809\ n^{8}$$

$$-\ 11568643258461073405846882745868241603747370565507314915173 77915\ n^{7}$$

$$-\ 15098331737373093592181438003367539114706112721445794448283 45557\ n^{6}$$

$$+\ 28431819235290177091195598377985429473932854364922349130012 07293\ n^{5}$$

$$+\ 18013778716605724861928056252978411502973175411980191956695 03339\ n^{4}$$

$$-\ 41236577692553675838489883568460331991426190300431462500048 58655\ n^{3}$$

$$-\ 36567787545144579592531788600718864284135646588560512829534 2841\ n^{2}$$

$$+\ 26103456978048524858124710074337995638333442138499808174454 43589\ n$$

$$-\ 87011523260161749527082366914459985461111473794999360581514 7863\ )\ n\ /$$

$$43471890$$

$$s_{n,81} := (\ 13090\ n^{78} + 510510\ n^{77} + 6211205\ n^{76} - 12932920\ n^{75} - 743518545\ n^{74}$$

$$+\ 1499970010\ n^{73} + 106877343265\ n^{72} - 215254656540\ n^{71} - 15227929505635\ n^{70}$$

$$+\ 30671113667810\ n^{69} + 2075338608492765\ n^{68} - 4181348330653340\ n^{67}$$

$$-\ 268210930396681745\ n^{66} + 540603209124016830\ n^{65} + 3276438198131 4524585\ n^{64}$$

$$-\ 66069367171753066000\ n^{63} - 3775439018352184271220\ n^{62}$$

$$+\ 7616947403876121608440\ n^{61} + 4096009430019735283529 65\ n^{60}$$

$$-\ 8268188334078231783143 70\ n^{59} - 41759374705920082223535985\ n^{58}$$

$$+\ 84345568245247987625386340\ n^{57} + 39927568054543244779187 94105\ n^{56}$$

$$-\ 80698591791538969434629745 50\ n^{55} - 357262346236485192232346107205\ n^{54}$$

$$+\ 722594551652124281408155188960\ n^{53} + 2984713070121330 6564421811260205\ n^{52}$$

$$-\ 6041685595407873741025177770937 0\ n^{51} - 2322483402861493092537880635510505\ n^{50}$$

$$+\ 4705383661677064922486013048730380\ n^{49}$$

$$+\ 16787696277020257036250370449374 4545\ n^{48}$$

$$- \, 3404593092020822056474934220 36219470 \; n^{47}$$
$$- \, 11240414081872625672333858777 818043379 \; n^{46}$$
$$+ \, 22821287472947333550315210977 672306228 \; n^{45}$$
$$+ \, 69500819292595180495326374409 1144221173 \; n^{44}$$
$$- \, 14128376733248509434568426991 59960748574 \; n^{43}$$
$$- \, 39551595153991326147050643313 410536755265 \; n^{42}$$
$$+ \, 80516027981307503237558129325 981034259104 \; n^{41}$$
$$+ \, 20640752988149614331142484615 68229033488057 \; n^{40}$$
$$- \, 42086666256112303694660550524 62439101235218 \; n^{39}$$
$$- \, 98389326683004321398238106165 340289530260081 \; n^{38}$$
$$+ \, 20098731999161987316594226738 3143018161755380 \; n^{37}$$
$$+ \, 42651573488207541082140595755 66228971674378621 \; n^{36}$$
$$- \, 87313020176331280895940614185 1560096151 0512622 \; n^{35}$$
$$- \, 16733761488402784335862915892 136390851 3113071177 \; n^{34}$$
$$+ \, 34340653178568881480685237926 124341798 7736654976 \; n^{33}$$
$$+ \, 59101717350452103518562209602 979644621 33819121825 \; n^{32}$$
$$- \, 12163750001876109518519294299 857172342 255374898626 \; n^{31}$$
$$- \, 18679240259468880634643023564 639364025 826290333723 \; n^{30}$$
$$+ \, 38574855519125387078780534142 913590039 07955566072 \; n^{29}$$
$$+ \, 52474773299065243703484557686 510906533 23134390505829 \; n^{28}$$
$$- \, 10880703215004302611484716878 731317207 040176736577730 \; n^{27}$$
$$- \, 13003512054453734234087536858 760659435 2695225866331473 \; n^{26}$$
$$+ \, 27095094430407898729323545405 394450591 2430628469240676 \; n^{25}$$
$$+ \, 28176870217712512680387031224 388552466 92025023131610521 \; n^{24}$$
$$- \, 59063249878465815233706416989 316549992 96480674732461718 \; n^{23}$$
$$- \, 52850035231762144811847857781 990744355 02128679856672613 \; n^{22}$$
$$+ \, 11160639545137087114706639325 532980387 030073038034445806944 \; n^{21}$$
$$+ \, 84789611234042823739402468225 100017691 5048945677561067725 \; n^{20}$$
$$- \, 18073986201322273459351157577 533015770 0398629389567942394 \; n^{19}$$
$$- \, 11470723611124724235680381497 822314802 397489403297413748297 \; n^{18}$$
$$+ \, 24748845842381675817295878753 401957624 95377435984395438988 \; n^{17}$$
$$+ \, 12859021236559203133812589964 647796318 98933279077773983869121 \; n^{16}$$
$$- \, 28192927057356573849354767804 635788614 22282033251532363177230 \; n^{15}$$
$$- \, 11685074152573832945570945824 732796142 94820977207569613404936 \; n^{14}$$
$$+ \, 26189441010883323276077368429 929171147 3192398766 6671589987102 \; n^{13}$$

$+ 836159814932215457531804630241479276983655545169030072607875 7\, n^{12}$

$- 1934214039973264147824382944782250265440503489104727304214461 6\, n^{11}$

$- 4524860732903394178344075193209052519120543532724749447925742 5\, n^{10}$

$+ 1098393550578005250451253333120035530368159054554226200065946 6\, n^{9}$

$+ 1739621482415523765114876165169964668089381878707339951098489 93\, n^{8}$

$- 4577636515409052780681005663459964866546922812870102522203574 52\, n^{7}$

$- 4233188482773304615429046653406695586739688670096688492364920 99\, n^{6}$

$+ 1304401348095566201153909897027335604002630015306347950693341 650\, n^{5}$

$+ 4701332371289202443666216121996376584946788663290138692283795 49\, n^{4}$

$- 2244667822353406689887153121426610920991987747964375689150100 748\, n^{3}$

$+ 2522186785750858496727528915687056058848791360321942387599025 11\, n^{2}$

$+ 1740230465203234990541647338289199709222229475899987211630295 726\, n$

$- 8701152326016174952708236691445998546111147379499936058151478 63\, ) \, (n+1)^2\, n^2$

$/ \ 1073380$

$s_{n,\,82} := (2\,n+1)\,(n+1)\,(19635\,n^{80} + 785400\,n^{79} + 9948400\,n^{78} - 15315300\,n^{77}$

$- 1184723540\,n^{76} + 1784742960\,n^{75} + 174371789020\,n^{74} - 262450055010\,n^{73}$

$- 25509102861242\,n^{72} + 38394879319368\,n^{71} + 3573256572686944\,n^{70}$

$- 5379082298690100\,n^{69} - 475056552894114092\,n^{68} + 715274370490516188\,n^{67}$

$+ 59749884762901886014\,n^{66} - 89982464329598087115\,n^{65}$

$- 7095114151169737348007\,n^{64} + 10687662458919405065568\,n^{63}$

$+ 794006757995298247516864\,n^{62} - 1196353968222407073808080\,n^{61}$

$- 83585404537804554448769552\,n^{60} + 125976283790818035210058368\,n^{59}$

$+ 8261045106439119424127015824\,n^{58} - 12454555801554088153795552920\,n^{57}$

$- 764964940033475942300835960632\,n^{56} + 1153674687950990957528151717408\,n^{55}$

$+ 66220648570651535866581331695424\,n^{54} - 99907810199952799278636073401840\,n^{53}$

$- 5346467027467586368448942306086688\,n^{52}$

$+ 8069654446301355952312731495830952\,n^{51}$

$+ 401572020013933708044618813899561116\,n^{50}$

$- 606392857244051240043084586597257150\,n^{49}$

$- 27983456318642265393990223090110637046\,n^{48}$

$+ 42278380906585423711006876928464584144\,n^{47}$

$+ 1803871892449478586406081903396226062432\,n^{46}$

$- 2726947029127510591464626293558571385720\,n^{45}$

$- 107225167685736175331790731053298438865176\,n^{44}$

$+ 162201225043168018293418409726726943990624\,n^{43}$

$+ 5857015014570034874082836697045318367307752\,n^{42}$

$$- 886662313437663632027096425043134102295694 0\, n^{41}$$

$$- 29289489343615661359977376248875399341795535 6\, n^{40}$$

$$+ 44377565172142323855979612585834 6660638411504\, n^{39}$$

$$+ 1335416107728382629961603165633035096792745427 2\, n^{38}$$

$$- 2025312944178645106870394554742469978221038716 0\, n^{37}$$

$$- 552620852158627146956639306782981862793346476296\, n^{36}$$

$$+ 8390578429588339459693109329481851440811249080 24\, n^{35}$$

$$+ 20652617563660082578709661948764791053240508319212\, n^{34}$$

$$- 31398455266969540841049148389621279151901324932830\, n^{33}$$

$$- 6931889542383125363225885713872344813388054385657 66\, n^{32}$$

$$+ 105548265899095357490440743127566236158415882042006 4\, n^{31}$$

$$+ 2076660272175229355112387152602612444725067780406291 2\, n^{30}$$

$$- 3167764541212391711413801100467701785166809611630440 0\, n^{29}$$

$$- 5514267261703550513492031498612897352888596268613860 64\, n^{28}$$

$$+ 8429789119615945355808737302942731185912348835023129 6\, n^{27}$$

$$+ 1287605350667525301835140075238277448221580617567666228 8\, n^{26}$$

$$- 197355697159936767953175379937212982792532710076901090 80\, n^{25}$$

$$- 2620088304800013722615153851884290466308079632706810079 2\, n^{24}$$

$$+ 4028810305779988967899318467795042190858398299944472057 28\, n^{23}$$

$$+ 459761203206981335862988169548259253158298315264351136662 4\, n^{22}$$

$$- 709785856339371948633978846661364090691739464396249065280 0\, n^{21}$$

$$- 687208999105777822726093933130796012889205594479863172095 2\, n^{20}$$

$$+ 10663027914756353315208399323026870646796781239179192907828\, n^{19}$$

$$+ 8622168866950967824920425583490468488611356597025565952773 34\, n^{18}$$

$$- 1346640496164264031410583413870465361510188017342448936991 5\, n^{17}$$

$$- 892009932486216488370407162509698850779695283156251600086480 7\, n^{16}$$

$$+ 1405346922210146079571316035471483508850298018743048624598216 8\, n^{15}$$

$$+ 744184856858705838382272931062382857054495952524061903629108 64\, n^{14}$$

$$- 118654463139856606155197519836714846102425882972324528667357380\, n^{13}$$

$$- 486469254960320065052583392358208996859571255326225742831402932\, n^{12}$$

$$+ 7890311140104084006564738484556709183405698244755008785807830 88\, n^{11}$$

$$+ 2396213745156922497686045741611726923809306904219407135596285084\, n^{10}$$

$$- 39888361747405879468573055366454258448842452685668611426848191 70\, n^{9}$$

$$- 8399221554843923494499652165719581523143258057737861139374193082\, n^{8}$$

$$+ 14593250419636179215178131016902085207157009720890222280403699208\, n^{7}$$

$$+ 19045771491058683152135554473478673375762033194926829894818443104 \, n^6$$

$$- 35865282446406114335792397218669052667221554652835355982429514260 \, n^5$$

$$- 22723458398898671082950467339602526674441538610947400567169802988 \, n^4$$

$$+ 52017828821551063792321899618738316345273085242838778841969461612 \, n^3$$

$$+ 46128389389832044065744593268561675452949028532816909201438269 66 \, n^2$$

$$- 32928172819250338506022638799653409490578896901341925801200471255 \, n$$

$$+ 10976057606416779502007546266551136496859623004473086004001570 85 \, ) \, n \, /$$

$$3259410$$

$$s_{n,\,83} := ( \, 1870 \, n^{80} + 74800 \, n^{79} + 935000 \, n^{78} - 1944800 \, n^{77} - 117317629 \, n^{76} + 236580058 \, n^{75}$$

$$+ 17742264353 \, n^{74} - 35721108764 \, n^{73} - 2663730788851 \, n^{72} + 5363182686466 \, n^{71}$$

$$+ 383082699495419 \, n^{70} - 771528581677304 \, n^{69} - 52322086431230341 \, n^{68}$$

$$+ 105415701444137986 \, n^{67} + 6765521316936166169 \, n^{66} - 13636458335316470324 \, n^{65}$$

$$- 826576475856030700306 \, n^{64} + 1666789410047377870936 \, n^{63}$$

$$+ 95249445487045212708884 \, n^{62} - 192165680384137803288704 \, n^{61}$$

$$- 10333784598611334098096881 \, n^{60} + 20859734877606805999482466 \, n^{59}$$

$$+ 1053545518695842134965934349 \, n^{58} - 2127950772269291075931351164 \, n^{57}$$

$$- 100733144131835595292365336301 \, n^{56} + 203594239035940481660662023766 \, n^{55}$$

$$+ 9013362339516683505312164225489 \, n^{54} - 18230318918069307492284990474744 \, n^{53}$$

$$- 753012498572262325582685861896249 \, n^{52}$$

$$+ 1524255316062593958657656714267242 \, n^{51}$$

$$+ 58593874928617733488098993172660165 \, n^{50}$$

$$- 1187120051732980609348556430595875 72 \, n^{49}$$

$$- 4235363648351463980854163178292983823 \, n^{48}$$

$$+ 858943930187622602264318199964555 5218 \, n^{47}$$

$$+ 283584122748604706219663336338020692607 \, n^{46}$$

$$- 5757576847990856384619698546756869 40432 \, n^{45}$$

$$- 17534344133684217083593161883976343963833 \, n^{44}$$

$$+ 356444459521675198056482936226283 74868098 \, n^{43}$$

$$+ 997846194540987232056134652222851692290437 \, n^{42}$$

$$- 2031336835034141983917917598068331759448972 \, n^{41}$$

$$- 52074503548188439744790988640791341949526473 \, n^{40}$$

$$+ 1061803439314110214734998948796510156585019 18 \, n^{39}$$

$$+ 2482261836284963130256046239619297684380350037 \, n^{38}$$

$$- 5070704016501337281985592374118246384419201992 \, n^{37}$$

$$- 10760554695981527064826230334657355231208161203 3 \, n^{36}$$

$$+ 220281797936131878578510199067265351008582426058 \, n^{35}$$

$$
\begin{aligned}
&+ 4221756456775766234372359525169006605656079379517\, n^{34} \\
&- 866379471148766434732322924940527856232074118 5092\, n^{33} \\
&- 14910757333534195004968311502575489493435619075 9423\, n^{32} \\
&+ 30687894138217156444668945930091506843103312270 3938\, n^{31} \\
&+ 471258080424619120829023087577592493429758719602 8047\, n^{30} \\
&- 973204054987455398102715121085276493702620751476 0032\, n^{29} \\
&- 1323884729364132783074479861053934450209292801999 60785\, n^{28} \\
&+ 2745089864227011105959231234216396549788847679146 81602\, n^{27} \\
&+ 3280652769833002336061764345367057943692509502976 281901\, n^{26} \\
&- 6835814526088705782719451814155755542363903773867 245404\, n^{25} \\
&- 7108735466069966532642742980357543451855921616974 6057741\, n^{24} \\
&+ 1490105238474880364355743114213066245794823361133 59360886\, n^{23} \\
&+ 1333352203180126982511900890978695827481544851472 765220529\, n^{22} \\
&- 2815714930207742001459376093378698279542572039058 889801944\, n^{21} \\
&- 2139154958930912792052786243205034100955528621816 1856288921\, n^{20} \\
&+ 4559881410882599784251510095747938029865314447538 2602379786\, n^{19} \\
&+ 2893945960848037471750845145275628566911161985724 60676539749\, n^{18} \\
&- 6243880062784334921926841300126050936808855416203 03955459284\, n^{17} \\
&- 3244199217903639523961399893354561494844960244120 402540533586\, n^{16} \\
&+ 7112786442085712540115483916721728083370806029861 109036526456\, n^{15} \\
&+ 2948024404777296670592820887136090269890583223674 0700792570324\, n^{14} \\
&- 6607327453763164595197190165944353348118247050334 2510621667104\, n^{13} \\
&- 2109545483860945159851174183325548973167835323093 10399874727901\, n^{12} \\
&+ 4879823713098206779222067383245533281147495351219 63310371122906\, n^{11} \\
&+ 1141575970733524626508965520355218427389278982194 09826019201089\, n^{10} \\
&- 2771134312776869930940145842395597013592605331560 782962409525084\, n^{9} \\
&- 4388886641430680875257069600498899034597032240750 143424935670451\, n^{8} \\
&+ 1154890759563823168145428504339339508278666981306 1069812280865986\, n^{7} \\
&+ 1067990052462700969447083944648971720896442761254 2589996963454779\, n^{6} \\
&- 3290870864489225107039596393637282950071552503814 6249806207775544\, n^{5} \\
&- 1186097955782100856425019750461082719492623152898 7109977952756517\, n^{4} \\
&+ 5663066776053426819889635894559448389056798809612 0469762113288578\, n^{3} \\
&- 6363218667433575095433086939694968951564729447165 617680256330119\, n^{2} \\
&- 4390423042566711800803018506620454598743852920178 9234401600628340\, n \\
&+ 2195211521283355900401509253310227299371926460089 4617200800314170 \,) \,(n+1)^{2} \\
&n^{2}\,/\,157080
\end{aligned}
$$

$$s_{n,84} := (2n+1)(n+1)(18723705\, n^{82} + 767671905\, n^{81} + 9979734765\, n^{80}$$

$$- 15353438100\, n^{79} - 1245675611180\, n^{78} + 1876190135820\, n^{77} + 192781863033760\, n^{76}$$

$$- 290110889618550\, n^{75} - 29696863069123250\, n^{74} + 4469035004849494150\, n^{73}$$

$$+ 438641281399576 2550\, n^{72} - 660196439601789 0900\, n^{71} - 615815995175170489700\, n^{70}$$

$$+ 927024974960764680000\, n^{69} + 8191577449840997 6474350\, n^{68}$$

$$- 12333717423509534 7051525\, n^{67} - 10304271305599414647679175\, n^{66}$$

$$+ 15518075545516669645044525\, n^{65} + 12236408914260106 06301177425\, n^{64}$$

$$- 184322037491177424427 4288400\, n^{63} - 13693749258299763712288186480\, n^{62}$$

$$+ 2063278490619523427905 69423920\, n^{61} + 144154934274085506185234 01379360\, n^{60}$$

$$- 2172640406564380209918038 6781000\, n^{59} - 1424735610048787480671450441097880\, n^{58}$$

$$+ 21479666171060031220567658 55037320\, n^{57}$$

$$+ 13192918272650729335726717908 5049160\, n^{56}$$

$$- 19896775739831386059629151555092400\, n^{55}$$

$$- 11420701610794043202285632338366449440\, n^{54}$$

$$+ 17230536294890221733726913083327220360\, n^{53}$$

$$+ 9220750082234171272371960138123 78467380\, n^{52}$$

$$- 13917277804825708017226574772602313 11250\, n^{51}$$

$$- 692568611165206370403414174254278147 50950\, n^{50}$$

$$+ 1045811555650222409613734548767718378 82050\, n^{49}$$

$$+ 4826148865407411393883450124955494172559450\, n^{48}$$

$$- 72915138758936282113058619148716271777 30200\, n^{47}$$

$$- 311103610255113084972642687458805999747396200\, n^{46}$$

$$+ 470301172320616441564616962145644813209959400\, n^{45}$$

$$+ 1849251985015159222907719274040687590110 2220000\, n^{44}$$

$$- 279739303613876965643980975916831362582583 09700\, n^{43}$$

$$- 101012634214602602270222813275147652273888 9586860\, n^{42}$$

$$+ 152917647839973288233554124792305635223746 3535140\, n^{41}$$

$$+ 505139301511909320348203150882942259704171 23305220\, n^{40}$$

$$- 76535483465986264493398243256402867131744 416725400\, n^{39}$$

$$- 230311683475293799751814587242810630783340 8706639960\, n^{38}$$

$$+ 349294299386240012852391793027036089531598 5268322640\, n^{37}$$

$$+ 953074012269639345524013853023193533623116 97717057020\, n^{36}$$

$$- 144470757333737710189286403691861421049112 5539209746850\, n^{35}$$

$$- 356184045686688108776644817234273372054761 1600985589910\, n^{34}$$

$$+ 541511447196901018259610427697340768606698 0171083258290\, n^{33}$$

$$+ 1195503888961627409517311957691741676163557761399430197770\, n^{32}$$
$$- 1820331405802286165188948457922479552675671542954561588100\, n^{31}$$
$$- 3581498834131864370006878026036939182813315912834088092720\, n^{30}$$
$$+ 5463264821487910863269764461951532751853757446398860218480\, n^{29}$$
$$+ 9510145705342566533226984733290666447294039014745518264940\, n^{28}$$
$$- 14538381799088245343003965323033576308533746394438220408340\, n^{27}$$
$$- 22206603188913751305475500350875920906984140050673182506518\!00\, n^{26}$$
$$+ 34036823873325039225363448792465560175902897395731684780194\!00\, n^{25}$$
$$+ 45187185090869668762984619984251839937595037684160591387801\!000\, n^{24}$$
$$- 69482618829970755105745102416001037915187701396027471320711\!200\, n^{23}$$
$$- 79292421361731725701560920927598121399446201983414386769632\!8968\, n^{22}$$
$$+ 12241276298409612630762863651219723399492868800449229537204\!849052\, n^{21}$$
$$+ 11851905976533218393730167698002721994344065929187167964814\!312906\, n^{20}$$
$$- 18389922779720308222133394729565069161490742334005366715823\!893885\, n^{19}$$
$$- 14870168297834744051705040042823665043102001159556348491143\!3809263\, n^{18}$$
$$+ 23224748585738131488664229800713751022727538856034791072506\!2660837\, n^{17}$$
$$+ 15383992153358132391617527213216575820771838240655513467716\!22208001\, n^{16}$$
$$- 24237225659324105161859502309860551282294134303785009755199\!64642420\, n^{15}$$
$$- 12834536457069305708363003172084066080286175081943224753506\!966178508\, n^{14}$$
$$+ 20463665968570163820637479873619126684543969338104087618020\!431588972\, n^{13}$$
$$+ 83898608396663591517761120108574471060638904636040687811004\!941687296\, n^{12}$$
$$- 13607974557928046918696042009967126993323034162311307552551\!7628325430\, n^{11}$$
$$- 41326146840671663310100416570162397191319216501589264968738\!1153185778\, n^{10}$$
$$+ 68793207539971518424498645860227159283640341833539551229383\!0543941382\, n^{9}$$
$$+ 14485663644337507936721338280509077505980103593100316856002\!7856936566\, n^{8}$$
$$- 25168155843504837826306939713774974590079032630642025089869\!57057375540\, n^{7}$$
$$- 32847167783933425222038574682733479860033613671473801311335\!80412518308\, n^{6}$$
$$+ 61854829597652556746211331880987707085089936822531714511938\!49147465232\, n^{5}$$
$$+ 39189866948170876579514191440441911455074177227039770979076\!088026209206\, n^{4}$$
$$- 89712215221082593242376953101156720725156234251825525875610\!56613046425\, n^{3}$$
$$- 79555031236290653221965016603704629784260437981291241132860\!8082716515\, n^{2}$$
$$+ 56789362295984894604483229041134054830217182823106449107344\!0430597985\, n$$
$$- 18929787431994964868161076347044684943405727607702149702578\!13476865995\, )\, n\,/$$
$$3183029850$$

$$s_{n,\,85} := (\,870870\, n^{82} + 35705670\, n^{81} + 458222765\, n^{80} - 952151200\, n^{79} - 60198598460\, n^{78}$$
$$+ 121349348120\, n^{77} + 9566165423815\, n^{76} - 19253680195750\, n^{75} - 1511347957444325\, n^{74}$$

$$+\ 3041949595084400\ n^{73} + 229037970216430375\ n^{72} - 461117890027945150\ n^{71}$$

$$-\ 33010895993522672225\ n^{70} + 66482909877073289600\ n^{69}$$

$$+\ 4511067518642905995775\ n^{68} - 9088617947162885281150\ n^{67}$$

$$-\ 583380641893051574291720\ n^{66} + 1175849901733266033864590\ n^{65}$$

$$+\ 71276731961455102062145165\ n^{64} - 143729313826434701581549 20\ n^{63}$$

$$-\ 8213545540654731384506382620\ n^{62} + 16570820395134106239170920160\ n^{61}$$

$$+\ 891104183169760288540734006175\ n^{60} - 1798779186734654683320638932510\ n^{59}$$

$$-\ 90849511579606832116154896178555\ n^{58} + 183497802345948318915630431289620\ n^{57}$$

$$+\ 8686438230030500214286000876736395\ n^{56}$$

$$-\ 17556374262406948747487632184762410\ n^{55}$$

$$-\ 777241872845096628478096424518009055\ n^{54}$$

$$+\ 1572040119952600205703680481220780520\ n^{53}$$

$$+\ 64933908908679570376830990266946881815\ n^{52}$$

$$-\ 131439857937311740959365661015114544150\ n^{51}$$

$$-\ 5052677547340188286925270455970937977525\ n^{50}$$

$$+\ 10236794952617688314809906572956990499200\ n^{49}$$

$$+\ 365224639041849768334639272583744497121975\ n^{48}$$

$$-\ 740686073036317224984088451740445984743150\ n^{47}$$

$$-\ 24455407702386920104577924279059370800185 2617\ n^{46}$$

$$+\ 49648840120774719316542574032927861988448384\ n^{45}$$

$$+\ 1512024713741066458157510498049399660778577599\ n^{44}$$

$$-\ 3073698267602907635631563570131727183545603582\ n^{43}$$

$$-\ 86046452331842361746397573163768361743750495775\ n^{42}$$

$$+\ 175166602931287631128426709897668450671046595132\ n^{41}$$

$$+\ 4490497946261768767396136234623256209256554925011\ n^{40}$$

$$-\ 9156162495454825165920699179144180869184156445154\ n^{39}$$

$$-\ 214050848657858968744205507378102846427632696327603\ n^{38}$$

$$+\ 437257859811172762654331713935349873724449549100360\ n^{37}$$

$$+\ 9279060859073620018078361005658487254115630724748983\ n^{36}$$

$$-\ 18995379577958412798811053725252324381955710998598326\ n^{35}$$

$$-\ 364051261309416244244618730930890192132640425668781981\ n^{34}$$

$$+\ 747097902196790901288048515587032708647236562336162288\ n^{33}$$

$$+\ 12857871054213925072973202599586037250969545 9589449247655\ n^{32}$$

$$-\ 26462840010624641047234453714759107210586335741234657598\ n^{31}$$

$$-\ 4063761146275676928972519065630941795246779082490153231 69\ n^{30}$$

$+\ 839215069265760026841738266840947466259942152239265303936\ n^{29}$

$+\ 1141614658466999756948515724948313599816908964885003833164\,7\ n^{28}$

$-\ 2367150823860575516581205276580721946259812144993934196723\,0\ n^{27}$

$-\ 282897839087589950974724330553421017484966180870326549790811\ n^{26}$

$+\ 589467186413785657115260713872649254432530483190592441548852\ n^{25}$

$+\ 613001753946326679124639751503353363985294311320152523519370\,7\ n^{24}$

$-\ 1284950226534031923960805574393971653413841670959364291193626\,6\ n^{23}$

$-\ 1149778667498347137416817557052859237867216563226457438703204\,47\ n^{22}$

$+\ 2428052357650097467229715671545115641075817293548851306525771\,60\ n^{21}$

$+\ 1844663994763414883026853734904211117833524042701988720274031\,27\ n^{20}$

$-\ 3932085131033307407260046356962933787674286537895282874707383\,414\ n^{19}$

$-\ 2495512680550595104363670080344796047754097789793016966856945332\,4\ n^{18}$

$+\ 5384233874204520949453344796385885474275624233375562221184629006\,2\ n^{17}$

$+\ 2797543698479578303012079984993793554325458734520772193255767321\,25\ n^{16}$

$-\ 6133510784379608700969494449626175656078479892379100608629997543\,12\ n^{15}$

$-\ 2542145701483150504120876254774659775681672853971412773910933231\,976\ n^{14}$

$+\ 5697642481404261878338701954511937116971119369718073560868684866218\,264\ n^{13}$

$+\ 1819107051892067284607743192899106755356722823516624327186172729082\,3\ n^{12}$

$-\ 4207978351924560757049356581249407222410565016751322215240832079991\,0\ n^{11}$

$-\ 984405842168024694543943676745296777575611350808660615585179784996\,53\ n^{10}$

$+\ 238960951952850546479282301161553427739227920329245345269444277799216\ n^{9}$

$+\ 378463261421089794325533576385923881731708654731092075071990454986471\ n^{8}$

$-\ 995887474795030135130349453933401191202645229791429495413425187772158\ n^{7}$

$-\ 920951100911906643790367273044777614120544492659506008378740907487489\ n^{6}$

$+\ 2837789676618843422711084000022956419443734215110441512170907002747136\ n^{5}$

$+\ 1022798120308369752758794369026701382867689843693645493636962672567167\ n^{4}$

$-\ 4883385917235582928228672738076359185179113902497732499444832347881470\ n^{3}$

$+\ 548714215418294977298228734333711098248984190478651279464602697074740\ n^{2}$

$+\ 37859574863989929736322152694089369886811455215404299405156269537319\,90\ n$

$-\ 18929787431994964868161076347044684943405727607702149702578134768659\,95\ )$

$(n+1)^{2}\ n^{2}\ /\ 74894820$

$s_{n,\,86} := (2\,n+1)\,(n+1)\,(15015\,n^{84} + 630630\,n^{83} + 8408400\,n^{82} - 12927915\,n^{81}$

$-\ 1098872775\,n^{80} + 1654773120\,n^{79} + 178598119700\,n^{78} - 268724566110\,n^{77}$

$-\ 28932208775470\,n^{76} + 4353267544626\,0\,n^{75} + 4500021405162920\,n^{74}$

$-\ 6771798445467510\,n^{73} - 666174189101622838\,n^{72} + 1002647182875168012\,n^{71}$

$+\ 93576152807981824526\,n^{70} - 140865552803410320795\,n^{69}$

$$- 1244915475287636363575963\, n^{68} + 1874416490571625052432\, n^{67}$$

$$+ 156604422399674603547065656\, n^{66} - 23584384184477917784681555\, n^{65}$$

$$- 18597060559647334507303042\, 3\, n^{64} + 280135127603934006198779712\, n^{63}$$

$$+ 2081198903645480727018994865\, 6\, n^{62} - 3135805111848417790843123284\, 0\, n^{61}$$

$$- 2190891912362081641397648029768\, n^{60} + 3302016894102364551050664201072\, n^{59}$$

$$+ 216533838397227248584313979593696\, n^{58} - 32645176604289205515199630149108\, 0\, n^{57}$$

$$- 20050831214368513484584137656478808\, n^{56}$$

$$+ 302394727045742162544522046354637\, 52\, n^{55}$$

$$+ 1735738503280314357617462572529528876\, n^{54}$$

$$- 261872749127275864455419961112025190\, n^{53}$$

$$- 14013859676853462462471962282825727712\, 6\, n^{52}$$

$$+ 2115172588984383162593561442229419282\, 84\, n^{51}$$

$$+ 10525780713802976536338433574523702731\, 312\, n^{50}$$

$$- 158944297001536839626373284338970250611\, 10\, n^{49}$$

$$- 73348667318608230724811853065085323948625\, 4\, n^{48}$$

$$+ 1108177224629200302853496460193228371759936\, n^{47}$$

$$+ 47282079039202710977246357533953798656601848\, n^{46}$$

$$- 71477207171118666617296284531027312170782740\, n^{45}$$

$$- 281052599961315028831278055428902369011669464\, 4\, n^{44}$$

$$+ 4251527603005284765777818973699049191260433336\, n^{43}$$

$$+ 15352079492171537025904273119175676301999761796\, 8\, n^{42}$$

$$- 2324069561840756977714530062744846691256266436\, 20\, n^{41}$$

$$- 767719679001637845156481055973419621352194678400\, 4\, n^{40}$$

$$+ 1163199866311660552623294234273853665484573349781\, 6\, n^{39}$$

$$+ 350031785645629862935560513454213591483159341192428\, n^{38}$$

$$- 530863677800003097166457241352689655552161878537550\, n^{37}$$

$$- 1448498805328663527895511275530331668327626827360641\, 4\, n^{36}$$

$$+ 2199291391882995446701589775363131985269048334967839\, 6\, n^{35}$$

$$+ 5413348365523750331540083408213619411084692586445548\, 8\, n^{34}$$

$$- 8229987117879775269645204601088585715890491296410724\, 30\, n^{33}$$

$$- 18169480361791494794076549682805665542951276512376271558\, n^{32}$$

$$+ 27665719898581230954597084754262927600221439333384943552\, n^{31}$$

$$+ 544322551631839806886569096646460111878611481305181339376\, n^{30}$$

$$- 8303166873970503258071521873468216316180279416244644804\, 40\, n^{29}$$

$$- 14453688292146060402000134275474884080060328382475016306184\, n^{28}$$

$+ 2209569078191761576590377750688573693589950654452475669969 \, n^{27}$

$+ 3374998979664680628194992946132059641137518301699184887156448 \, n^{26}$

$- 517297692340660902112200830673251846385774985271401114233204 \, n^{25}$

$- 68676286183084696465554189389707154387182067022182215471672484 \, n^{24}$

$+ 105600917736330349208892288237926990653965988025909023764625324 \, n^{23}$

$+ 1205100297935646990403117926442732080939718575466621605766240464 \, n^{22}$

$- 18604509057716356602091230337830616167365608572128869205316733354 \, n^{21}$

$- 1801273712940064197959760666633130640304297359086806947379657783 \, n^{20}$

$+ 279493311469867807995009715163884904129327408146355476709607034254 \, n^{19}$

$+ 225994579346075056810031986469037064682062185888599232880999901654 \, n^{18}$

$- 352973852475404648921255283728549804908775698296076587676303369546 \, n^{17}$

$- 2338086441184023948416697999678402581237905055802062428475253970034 \, n^{16}$

$+ 3683616588013738247085674644138187877431124543284839747209669612352 \, n^{15}$

$+ 195061563799649746906098013759195731395885044747075247549048120603654 \, n^{14}$

$- 311010428639543311594575393845702990965383794284854858684055661523044 \, n^{13}$

$- 1275105946304240023554226589816310925771918902753517961663169496103854 \, n^{12}$

$+ 2068164133776131691128627581647317884140570251272704371836782074917246 \, n^{11}$

$+ 6280821169910865757734600511296806245743604385722125931372885082813654 \, n^{10}$

$- 10455313821754364482166214557768868310685691704219541082977718661679046 \, n^{9}$

$- 220155687943358468791479684807085931313086232389582003066780576585597454 \, n^{8}$

$+ 38251010102380952559805059999947323852305780710547071001505945818738564 \, n^{7}$

$+ 49921708787498993839709645840897254247772536907908868811080706277442784 \, n^{6}$

$- 940080682324389670394669987613195432978116957171368387173740323255095 \, n^{5}$

$- 595614555896160704978040899097330098037734485387456383106528616387575 \, n^{4}$

$+ 136346217500643589266439634245259286354566020666668768246663086208910 \, n^{3}$

$+ 120909148943461844816391427403910728081587674530213896254682617182240 \, n^{2}$

$- 8630948109184107135567853123321625238952116151287552285053554688781544 \, n$

$+ 287698270306136904518928437444054174631737205042918409501785156292605 \, ) \, n \, / \\ 2612610$

$s_{n,87} := ( 2730 \, n^{84} + 114660 \, n^{83} + 1509690 \, n^{82} - 3134040 \, n^{81} - 207443600 \, n^{80} \\
\quad + 418021240 \, n^{79} + 34596931460 \, n^{78} - 69611884160 \, n^{77} - 5744572476097 \, n^{76} \\
\quad + 11558756836354 \, n^{75} + 916135636022789 \, n^{74} - 1843830028881932 \, n^{73} \\
\quad - 139140449487152113 \, n^{72} + 280124729003186158 \, n^{71} + 20064736718560512797 \, n^{70} \\
\quad - 4040598166124211752 \, n^{69} - 274228220831464278158 \, n^{68} \\
\quad + 5524974014795409788068 \, n^{67} + 35464910994532547566 2122 \, n^{66} \\
\quad - 714823193905446361112312 \, n^{65} - 43330941807345476947083928 \, n^{64}$

$$
\begin{aligned}
&+ 8737670680859640025528016 8\, n^{63} + 4993233631235600338727821292\, n^{62} \\
&- 10073843969279797077710922752\, n^{61} - 5417262943230711323395883 56963\, n^{60} \\
&+ 109352643261542206175688763 6678\, n^{59} + 55229878018047002386939658730167\, n^{58} \\
&- 111553282468709426835636205097012\, n^{57} \\
&- 528072116487642493352099394 6914431\, n^{56} \\
&+ 1067299561222155929387762409 8925874\, n^{55} \\
&+ 472506398003879216761448176274437883\, n^{54} \\
&- 9556857916199799928167739766 47801640\, n^{53} \\
&- 394750830301363779204824777 89738991881\, n^{52} \\
&+ 79905851851892735833781729556125785402\, n^{51} \\
&+ 3071659617479348446524062123697928009277\, n^{50} \\
&- 6223225086810589628881905976951981803956\, n^{49} \\
&- 220299956329830890703808036789620710774369\, n^{48} \\
&+ 450283137746472371036497979556193403352694\, n^{47} \\
&+ 14866296173493996637559625824531062574308541\, n^{46} \\
&- 3018287548473446564615574962861831855196977 6\, n^{45} \\
&- 919200802154100709595282060673490345151891599\, n^{44} \\
&+ 1868584479792935884836719870975599008855752974\, n^{43} \\
&+ 52309970390887088857679431666604318401120358051\, n^{42} \\
&- 106488525261567113600195583204184235811096469076\, n^{41} \\
&- 2729895402351983520963640521094615787033613943699\, n^{40} \\
&+ 5566279329965534155527476625393415809878324356474\, n^{39} \\
&+ 130127312073962470257195362199242754970578667343951\, n^{38} \\
&- 265820903477890474669918201023878925751035659044376\, n^{37} \\
&- 5640992575983711518212166304393582651813679917581809\, n^{36} \\
&+ 11547806055445313511094250809811044229378395494207994\, n^{35} \\
&+ 221316628214135883323415088294942973159737667679426821\, n^{34} \\
&- 454181062483717080157924427399696990548853730853061636\, n^{33} \\
&- 7816648286001045530983711212410336738631192684869215753\, n^{32} \\
&+ 16087477634485808142125346852220370467811239100591493142\, n^{31} \\
&+ 247047053628236914386227817007476264277111303683248873669\, n^{30} \\
&- 5101815848909596369145809808671728990220338464670892 40480\, n^{29} \\
&- 6940184907559175634190854859359073702007075887878631571587\, n^{28} \\
&+ 14390551400009310905296290699585320303036185622224352383654\, n^{27} \\
&+ 171981263437285330499520883452613567750639198401012921797399\, n^{26}
\end{aligned}
$$

$$- 3583530782745799719043380576048124558043145824242501959784 52\, n^{25}$$

$$- 37266037970802550828702565040566973296386940535090682269 51583\, n^{24}$$

$$+ 78115606724350901376448510657182071150817026894423866498 81618\, n^{23}$$

$$+ 69898161310585469637345432163682990159603074766587004691 947867\, n^{22}$$

$$- 14760788329360602941233571539308418743428785222616396033 777352\, n^{21}$$

$$- 11214066173293915424599894849631501344695577429061601494 05124338\, n^{20}$$

$$+ 23904211179523891143323146853193844563734033380349366948 44026028\, n^{19}$$

$$+ 15170897915270966469464876666202256777074303286143211127 545636382\, n^{18}$$

$$- 32732216948494322053262068017723898010522009910321358949 935298792\, n^{17}$$

$$- 17007026329266849587958708670657649558538752298098400962 2869972 28\, n^{16}$$

$$+ 37287274353383131381243624143087688918129705587228937820 2509293248\, n^{15}$$

$$+ 15454392687933232093240447672362173848171418642121323089 3466567232\, n^{14}$$

$$- 34637512811204777324605257759033116588155807843147158399 89442427712\, n^{13}$$

$$- 11058844780189668919333399952482269532406501013664407810 061870992673\, n^{12}$$

$$+ 25581440841499815571127325680867850723628582811643531460 113184413058\, n^{11}$$

$$+ 59844699067737961007369361034856579030750532332850421275 09484419557\, n^{10}$$

$$- 14527083897697573758586604775058104652977868927821361571 5132153252172\, n^{9}$$

$$- 23007807367394596214937046635654107116972245604254379690 3354419378001\, n^{8}$$

$$+ 60542698632486766188460698046366318886922360136330120952 1840992008174\, n^{7}$$

$$+ 55987113372665062621764853103663778800505594983947380419 681974181373\, n^{6}$$

$$- 17251692537781689143199040425369387664702347913311959703 61204940370920\, n^{5}$$

$$- 62178669706081328032083574858803420839213048553148467932 0743233205690\, n^{4}$$

$$+ 29687426478997954749615755397130071832544957623941653290 02691406782300\, n^{3}$$

$$- 33357824272535011940507402008028689310029906102540902649 4205078220730\, n^{2}$$

$$- 23015861624490952361514274995524333970538976403433472760 14281250340840\, n$$

$$+ 11507930812245476180757137497762166985269488201716736380 07140625170420\,)$$

$$(n + 1)^{2}\, n^{2}\, /\, 240240$$

$$s_{n,\,88} := (2\,n + 1)\,(n + 1)\,(31395\, n^{86} + 1349985\, n^{85} + 18449795\, n^{84} - 28349685\, n^{83}$$

$$- 2521861979\, n^{82} + 3796967811\, n^{81} + 429942952803\, n^{80} - 646812913110\, n^{79}$$

$$- 73154079227866\, n^{78} + 110054525298354\, n^{77} + 11965748780578502\, n^{76}$$

$$- 18003650433516930\, n^{75} - 1865300844139392790\, n^{74} + 2806953091425847650\, n^{73}$$

$$+ 276285728807145681260\, n^{72} - 415832069756431445715\, n^{71}$$

$$- 38814504754627096497925\, n^{70} + 5842967316681886046 9745\, n^{69}$$

$$+ 5163963587656899219272195\, n^{68} - 7775160218068758259143165\, n^{67}$$

$$- 649607337897355801093343395\, n^{66} + 9782985869550680807 69586675\, n^{65}$$

$$+ 771422143178745187996183 09535\, n^{64} - 1162024707702893122 39812257640\, n^{63}$$

$$- 8632997342899032700593244878040\ n^{62} + 1300759724973369370700977344 5880\ n^{61}$$

$$+ 9088014821110586995913616609885 20\ n^{60}$$

$$- 1369706021791454896240547378205720\ n^{59}$$

$$- 89820167959401101255406201989254168\ n^{58}$$

$$+ 135415104949997379331229576672984112\ n^{57}$$

$$+ 831726373609600769028257005196348171 6\ n^{56}$$

$$- 125436031566190102250946986628171463 0\ n^{55}$$

$$- 719999823596868812517233573903384098 282\ n^{54}$$

$$+ 1086271536973612723888395095788217004 738\ n^{53}$$

$$+ 5813074078360334227438084174305883821 5174\ n^{52}$$

$$- 8773924694389181977351546016248236582 5130\ n^{51}$$

$$- 436618779125507866733826458167295068004 4790\ n^{50}$$

$$+ 659315131035456391089415460259066720297 9750\ n^{49}$$

$$+ 3042568190131261363186737997968503499079 09430\ n^{48}$$

$$- 4596818041748664864334577769965708546335 4020\ n^{47}$$

$$- 1961302841729329235667057332741768161774 9503100\ n^{46}$$

$$+ 2964938352802737177822258887962480785585 5931660\ n^{45}$$

$$+ 1165831270920953527358339545845489792565 760040260\ n^{44}$$

$$- 1763571598145443976266206132080470927765 68026220\ n^{43}$$

$$- 6368179603430031257162371738307811578713 2953694340\ n^{42}$$

$$+ 9640447985052319084589888638122119722708 7714554620\ n^{41}$$

$$+ 3184569753865383546549473105625868718266 406814616640\ n^{40}$$

$$- 4825056870723336915247159101629413676013 154079202270\ n^{39}$$

$$- 1451963090106250023084829080669977661727 05767988072850\ n^{38}$$

$$+ 2202069919512991719203479416513113560970 65229021710410\ n^{37}$$

$$+ 6008502335069268989423338056043980790540 698041766397310\ n^{36}$$

$$- 9122856998579553070095181054891626863859 579677160451170\ n^{35}$$

$$- 2245505220655861396858819889570165841160 63082130817282462\ n^{34}$$

$$+ 3413872115976689860638705739629706896060 2441303480614 9278\ n^{33}$$

$$+ 7536862631795457049040647341755340120027 394922820409884134\ n^{32}$$

$$- 1147598755349202006659290629961449552484 4104590748017900840\ n^{31}$$

$$- 2257898529483903092271819656057578328883 45178018988355866008\ n^{30}$$

$$+ 3444227731993314738740694015584439970949 39819323856542749432\ n^{29}$$

$$+ 5995518914771771457602705955537700834320 796710374876956565576\ n^{28}$$

$$- 9165489758757322923341093634085773250028 66497522424 3706223080\ n^{27}$$

$$-13999797014378937527619067499578952941535381794150744381526536 0 \, n^{26}$$

$$+214579700095062724375956559310727180748045093998732875760009580 \, n^{25}$$

$$+2848753650171790105083951328870792166772305211171826510700335670 \, n^{24}$$

$$-4380420325305216519813905272961551840532480346457676409838508295 \, n^{23}$$

$$-4998863601061845258246221365298898647775646553952135369239553969 7 \, n^{22}$$

$$+7717316417858028713360023115964255636900938482510868743512563693 \, n^{21}$$

$$+7471844139936015110181241837111357199667089880515250809930002675 9 \, n^{20}$$

$$-11593632030796924100939864121246524858134568174489841965207063219 85 \, n^{19}$$

$$-937465923844079778771510795979257104050896615866230848469617812340 7 \, n^{18}$$

$$+146416704592010428866196551457511828036701776467179548253046203461 03 \, n^{17}$$

$$+96985912516927929226268040852379639855950039852935851748996148709799 \, n^{16}$$

$$-1527997040049924152827118888514450511857601486027627550361465332377 50 \, n^{15}$$

$$-8091327775079034447797292045239289101289349708427551448636459556138 82 \, n^{14}$$

$$+12900990182643513748109497512116158907862825306555140948135422000396 98 \, n^{13}$$

$$+528925327908114462192026943121431498679484183705417497172879773347709 4 \, n^{12}$$

$$-8578929427753892620285879022427280425585404020864019504999967700235490 \, n^{11}$$

$$-26053406828319212150308695492903426134624745173983335680822355673090198 \, n^{10}$$

$$+43369574956355764535605982750568779414729819771407013273733517359753042 \, n^{9}$$

$$+9132254443153026177644292709416195650070921153871788825891724388618675 6 \, n^{8}$$

$$-158668604125473274932467382016527324458428727193780339025242624509156655 \, n^{7}$$

$$-20707970397826935672145991188195309049788914757095094411988609757387636 1 \, n^{6}$$

$$+38995385803014067254842355883119329797604808495331658569245045861539286 9 \, n^{5}$$

$$+2470662341410298263419341150599721778197940441726256964539397319855052767 \, n^{4}$$

$$-565576280226615075787112952005554915717715108735596837535321209090275585 \, n^{3}$$

$$-5015419419646616119281209516966311691635484724224360075026568372105135 9 \, n^{2}$$

$$+3580194314080067796827746187572721332333898252311638198930591301267148 31 \, n$$

$$-1193398104693355932275915395857573777444632750770546066310197100422382 77 \, )$$

$$n \, / \, 5588310$$

$$s_{n,89} := (\, 20930 \, n^{86} + 899990 \, n^{85} + 12149865 \, n^{84} - 25199720 \, n^{83} - 1744421315 \, n^{82}$$

$$+ 3514042350 \, n^{81} + 304985960565 \, n^{80} - 613485963480 \, n^{79} - 53158009488090 \, n^{78}$$

$$+ 106929504939660 \, n^{77} + 8909989146292845 \, n^{76} - 17926907797525350 \, n^{75}$$

$$- 1424074208473841325\, n^{74} + 2866075324745208000\, n^{73} + 216399956208100338225\, n^{72}$$

$$- 43566598774094588 4450\, n^{71} - 321002203061121194 1650\, n^{70}$$

$$+ 6285571004896336976 7750\, n^{69} + 4265665972790862773 846025\, n^{68}$$

$$- 8594187655630688917 459800\, n^{67} - 5516670876020690903 972431035\, n^{66}$$

$$+ 1111928362859770496 862321870\, n^{65} + 6740268553484586098 6917803045\, n^{64}$$

$$- 1359172994325514924 70697927960\, n^{63} - 7767140761120224157 457480228700\, n^{62}$$

$$+ 1567019882167299980 7385658385360\, n^{61} + 8426733487248458406 06667131006405\, n^{60}$$

$$- 1701016896271364681 020719920398170\, n^{59}$$

$$- 8591192302492349155 9182655250921505\, n^{58}$$

$$+ 1735248629461183477 99386030422241180\, n^{57}$$

$$+ 8214338506673990004 612121910384769945\, n^{56}$$

$$- 1660220187629409835 7023629851191781070\, n^{55}$$

$$- 7349995173027919905 42006062563337341635\, n^{54}$$

$$+ 1486660123648187807 9441035754977866464340\, n^{53}$$

$$+ 6140481294855188613 946866769523 9181205605\, n^{52}$$

$$- 1242962271335856503 58378371145456228875550\, n^{51}$$

$$- 4778069348991508058 6019867598989932067537725\, n^{50}$$

$$+ 9680434925116601767 56235189094332036395 1000\, n^{49}$$

$$+ 3453750288253318924 28876397681742916769278025\, n^{48}$$

$$- 7004304925757803866 2531514725442915 3902507050\, n^{47}$$

$$- 2312502130053260594 0971953814789819885552959223\, n^{46}$$

$$+ 4695047309364099226 8569222776834068925008425496\, n^{45}$$

$$+ 1429847615118685544 9630061863876443 86461372590981\, n^{44}$$

$$- 2906645703331012082 194581595552122841847753607458\, n^{43}$$

$$- 8136990985548837606 6622932517626702636548989 15465\, n^{42}$$

$$+ 1656464654143086872 9551916809907746336915755 1438388\, n^{41}$$

$$+ 4246443673059623521 9235915613801915474528 01490949689\, n^{40}$$

$$- 8658533811533555731 14270229085946055827476053 3337766\, n^{39}$$

$$- 2024173895353886776 4260805864627405470320845709371207\, n^{38}$$

$$+ 4134933128823109110 1635881958340757949891645 1952080180\, n^{37}$$

$$+ 8774752766507124724 8343181912106323128724024 43762600597\, n^{36}$$

$$- 1796299884589656036 06849952020046722052437213 39477281374\, n^{35}$$

$$- 3442654230682732168 3974815648930449 50751459548955 7659789\, n^{34}$$

$$+ 7064938449824429940 401813081806136712202729 12318592600952\, n^{33}$$

$$+ 1215905804670199965 924857255383022359436430 2106364389902585\, n^{32}$$

$$- 250246099383864423125373264158410608599488771250473724061 22 \ n^{31}$$

$$- 384289961045336075273969804517974995019638849227505252928931 \ n^{30}$$

$$+ 79360453202905859286047693545179105089922657558005787826 3984 \ n^{29}$$

$$+ 1079568992466829422792842669937530773106937416747582873281 9813 \ n^{28}$$

$$- 2238498438136564704871733033420240651303797491053171534390 3610 \ n^{27}$$

$$- 26752261181101594799529569096926917245146625844472513983 7646977 \ n^{26}$$

$$+ 55743020800339754303930871227274075141597049179998199501919 7564 \ n^{25}$$

$$+ 579685694275241265622337290286447661216244504191808126442 0789049 \ n^{24}$$

$$- 12151144093508222855486054518001693975740860575636144523860 775662 \ n^{23}$$

$$- 10872892954071377227519225440069782286491695493746568849009 1086120 \ n^{22}$$

$$+ 229609003174935767405870563319397339705574770450567521504 042947902 \ n^{21}$$

$$+ 17443855288312490260316848441583167989866275984342259814029 32946441 \ n^{20}$$

$$- 37183800608374338194692402516360309376788299673190194843099 08840784 \ n^{19}$$

$$- 23598839505114786071305465757698721506221020331594829912219 90518523 \ n^{18}$$

$$+ 50916059071860391033730333403175775238923034033637985466753 889877830 \ n^{17}$$

$$+ 2645499871213179805266897885427487072963785036969284012915291 81264613 \ n^{16}$$

$$- 5800160333144963520871099104886731898316800414274947880498122 52407056 \ n^{15}$$

$$- 24039825113487522865546800254817578842316856669100433303364636 13502526 \ n^{14}$$

$$+ 53879810560120009253964699614521889582950513752475814487227394 79412108 \ n^{13}$$

$$+ 172024032157887002924123505617140338473033259694688874635808148 63664685 \ n^{12}$$

$$- 3979278748758940151022117108488025665290170331418535637588436920 6741478 \ n^{11}$$

$$- 9309043251379337771387306258168433173980890628432222093663168369 3202669 \ n^{10}$$

$$+$$

$$2259736525151761569379672962482489201325195158828297982491477365 93146816 \ n^9$$

$$+$$

$$3578941447429625763614262374530035542522277175628463260838679209 53265137 \ n^8$$

$$-$$

$$9417619420011013096608197711542560286369749510085224504168835784 99677090 \ n^7$$

$$-$$

$$8708982884450447864314238233024700905893436110993764764198492383 21018750 \ n^6$$

$$+$$

$$2683558518891190882523667417759196209815662173207275403256582055 141714590 \ n^5$$

$$+$$

$$9672100196409591974128852180274695174699312483132639419426598204 71618745 \ n^4$$

$$-$$

$$4617978558173109277349437853814135244755524669833803287141901696 084952080 \ n^3$$

$$+ 5188921220465207402608458331207069562108132087610825441056551974089018855\, n^2$$
$$+ 3580194314080067796827746187572721332333898252311638198930591301267148310\, n$$
$$- 179009715704003389841387309378636066616694912615581909946529565063357 4155\,)$$
$$(n+1)^2\, n^2\, /\, 1883700$$

$s_{n,90} := (2\,n+1)\,(n+1)\,(15646785\, n^{88} + 688458540\, n^{87} + 9638419560\, n^{86}$
$$- 14801858610\, n^{85} - 1376572850730\, n^{84} + 2072260205400\, n^{83} + 245903277087930\, n^{82}$$
$$- 369891045734595\, n^{81} - 43893838573205775\, n^{80} + 66025703382675960\, n^{79}$$
$$+ 7541112257223012160\, n^{78} - 11344681237525856220\, n^{77} - 1236272919168107037140\, n^{76}$$
$$+ 186008171937092 3483820\, n^{75} + 192823560726022460839990\, n^{74}$$
$$- 29016538194871915 3001895\, n^{73} - 285645889279286236598 18267\, n^{72}$$
$$+ 42991966082867295066228348\, n^{71} + 4013082100898518722077196424\, n^{70}$$
$$- 604111913438921173064890 8810\, n^{69} - 533913335461069588759992021842\, n^{68}$$
$$+ 80389056275879898900 05312487168\, n^{67} + 6716444286097072760811 5836504294\, n^{66}$$
$$- 1011486095728354909066676411000025\, n^{65}$$
$$- 79759205497059935652621748825939 97\, n^{64}$$
$$+ 120144551293454080933466005293910 08\, n^{63}$$
$$+ 892586638444623784719912620983142464\, n^{62}$$
$$- 134488718523160838112654223173940920 0\, n^{61}$$
$$- 9396320347370412751890195589124027392\, n^{60}$$
$$+ 14161724880317199545839856449955574568 8\, n^{59}$$
$$+ 9286726476072095490780996348919783397764\, n^{58}$$
$$- 14000898338509729233900693805629452969490\, n^{57}$$
$$- 85994220735734278009652113779776722947156 2\, n^{56}$$
$$+ 12969137602052690347617320535994655570692088\, n^{55}$$
$$+ 744425399450488046802753854942304721607 25424\, n^{54}$$
$$- 1123122667976758415377939442681454410264341 80\, n^{53}$$
$$- 601027924109060731160889181913127442463724626 0\, n^{52}$$
$$+ 9071574995034748888182234700830984357469086480\, n^{51}$$
$$+ 45143081769727880104825466128732751129245180186 0\, n^{50}$$
$$- 68168201404343557601647310928140675911741224 6030\, n^{49}$$
$$- 3145785549426344885040559489745575554635619 8645126\, n^{48}$$
$$+ 475276242484168910636166289008243366990300 4090704\, n^{47}$$
$$+ 20278389018806234313454663787221500731477 51715701632\, n^{46}$$
$$- 30655221649451435925500078825336372780711 74075597800\, n^{45}$$

$$- 1205381419892837757432878112958405468780096787267650 16 \, n^{44}$$

$$+ 1823399740663982354112067208850276389560501051279464 24 \, n^{43}$$

$$+ 6584216399043958005676602316287830656989002023314748 972 \, n^{42}$$

$$- 9967494585599136126220506834874259804961528087536096 670 \, n^{41}$$

$$- 3292604433770478898490110033370844961316814474283185 72806 \, n^{40}$$

$$+ 4988744123583714028366267584230638741000029351862459 07544 \, n^{39}$$

$$+ 1501220094913030971053305826883847964810754323827406 4717712 \, n^{38}$$

$$- 2276773862987465026721790078246925140921131632500422 0030340 \, n^{37}$$

$$- 6212337288186720632915158503886571687074461959402679 94768916 \, n^{36}$$

$$+ 9432344625429454200708827259742203787657749520729041 02168544 \, n^{35}$$

$$+ 2321682681502618288370062702403435399069476170089361 0026678092 \, n^{34}$$

$$- 3529685745381074703558638189903864117542503002737686 7091101410 \, n^{33}$$

$$- 7792546320597254795065620019411079474384812379846521 55335254314 \, n^{32}$$

$$+ 1186530376816493592777636193861181241745434371990066 666544382176 \, n^{31}$$

$$+ 2334496426136965916064773600379490804932512959069182 3563727163168 \, n^{30}$$

$$- 3561071158046273553736042210262295269486041157203306 8678862935840 \, n^{29}$$

$$- 6198913412893944773635910591250594732970065818660436 20627505333512 \, n^{28}$$

$$+ 9476423677243230838140667997389006862929400785850819 65280689468188 \, n^{27}$$

$$+ 1447473199966207515333267059787631141123404352802486 1991277058635074 \, n^{26}$$

$$- 2218591918335527427190603929668391745999753533132983 3969555932686705 \, n^{25}$$

$$- 2945395963737480375648524843361954687581054846929300 55357256880301837 \, n^{24}$$

$$+ 4529023541522996934832317461526351618671569947050600 00002066328 6796108 \, n^{23}$$

$$+ 5168447146334986623473684847322737902741225726868739 857848943701368984 \, n^{22}$$

$$- 7979121896578629781952143144060424435045417087655639 786783747195451530 \, n^{21}$$

$$- 7725322114151903065317951351487573395722550612675962 0847236132103949922 \, n^{20}$$

$$+ 1198693926605678608707453418443438131533609677339672 5116424607175365064 8 \, n^{19}$$

$$+ 9692689110078963604964227384621259314762315021561488 52350192555361918114 \, n^{18}$$

$$- 1513838062842128471180006778615360803791027737101206 90410741186891970249 5 \, n^{17}$$

$$- 1002761032880219238386526431367753283045310221997311 2250684071748416180827 \, n^{16} +$$

$$1579833452462435281138789985982397964757516719851027 182807981355708412248 8 \, n^{15} +$$

$$8365821371938423031403959143944614959718430319620571 140053489225636064774 4$$

$$n^{14} -$$

$$1333864878413885218767533370890812142195640383935637030148422451630830328 60$$

$$n^{13} -$$

$$546868812557778371190052711752787401496936842491361318021428443112828643 764$$

$$n^{12} +$$

$$886996462757361817723455736173721709355187282933823828539563787250784482 076$$

$$n^{11} +$$

$$269372535280874020794958206470418889611715429171989069551634542497334392 9918$$

$$n^{10} -$$

$$448408626059179122078610096514314419885332507904674795754430003108540813 5915$$

$$n^{9} -$$

$$944206087285383280016821940373039676196622593989090831540103210558193790 6895$$

$$n^{8} + 16405134439576644810645379588167167242376001449359736451873698173915610 9\backslash$$

$$28300 \ n^{7} + 2141047627030739359602170622465253925322823482815594297183122882 95\backslash$$

$$77640089960 \ n^{6} - 40318281625249412799355249131062392501030352966913782683 6836\backslash$$

$$92331324265599090 \ n^{5} - 25544781268500528848536847759734474567915405965189 8956\backslash$$

$$279066699683333795482650 \ n^{4} + 58476312715375499672482896205132908102388285 4312\backslash$$

$$417347837018506906628023023520 \ n^{3} +$$

$$518556461145276311074706967326843566928665761950513588149229113496447664 9710$$

$$n^{2} - 370165032748668945023620520526124691075551241291448785712140893620477 781279\backslash$$

$$86325 \ n + 1233883442495563150078735087082303585170804304829285707136312068 259\backslash$$

$$2709328775 \ ) \ n \ / \ 2847714870$$

$$s_{n,\,91} := (\ 1360590 \ n^{88} + 59865960 \ n^{87} + 828145780 \ n^{86} - 1716157520 \ n^{85} - 124119142455 \ n^{84}$$

$$+ 249954442430 \ n^{83} + 22724059650215 \ n^{82} - 45698073742860 \ n^{81} - 4152825388665800 \ n^{80}$$

$$+ 8351348851074460 \ n^{79} + 730689531882655730 \ n^{78} - 1469730412616385920 \ n^{77}$$

$$- 12274316719622831542 9 \ n^{76} + 246956064805073016778 \ n^{75}$$

$$+ 1962837430379067638927 3 \ n^{74} - 39503704672386425795324 \ n^{73}$$

$$- 2983086777642037946788216 \ n^{72} + 6005677259956462319371756 \ n^{71}$$

$$+ 43024608421003780010272495 4 \ n^{70} - 866497845680032062523921664 \ n^{69}$$

$$- 5880485851642821223004629683 1 \ n^{68} + 118476214878536456522616515326 \ n^{67}$$

$$+ 76050892330311714474412486060 79 \ n^{66} - 1532865468094087935140511372748 4 \ n^{65}$$

$$- 929190279004650758370293560073686 \ n^{64}$$

$$+ 1873709212690242396091992233874856 \ n^{63}$$

$$+ 10707514892802282052101072537546724 4 \ n^{62}$$

$$- 216024007068735883438113442984809344 \ n^{61}$$

$$- 1161680725133170457754729041768604480 1 \ n^{60}$$

$$+\ 23449638509732145038532694278356898946\ n^{59}$$

$$+\ 11843524653983192106569527531996245505309\ n^{58}$$

$$-\ 23921545693063705663524382006776059909564\ n^{57}$$

$$-\ 11324006887706175490488296706274461498 3103\ n^{56}$$

$$+\ 228872292323429880376118372326166835875770\ n^{55}$$

$$+\ 1013245265411682776551398732664914017 1097463\ n^{54}$$

$$-\ 204937776005570854114040930256244471780 70696\ n^{53}$$

$$-\ 8465058076408021814534836645305970101 14461045\ n^{52}$$

$$+\ 171350539288216144831837142208681846 7406992786\ n^{51}$$

$$+\ 658688343654743490218128248839137527 59202224873\ n^{50}$$

$$-\ 133451174123830859491944021189914323 985811442532\ n^{49}$$

$$-\ 4761222348653191524485028971847808478539285664923\ n^{48}$$

$$+\ 9655895871430213908462001964885531281064382772378\ n^{47}$$

$$+\ 31879365628677149067007118898202137984124619722587\ n^{46}$$

$$-\ 64724320844497319524860437992892980724 9313622217552\ n^{45}$$

$$-\ 19711391537023688336246210750837438165 75924705672633\ n^{44}$$

$$+\ 40070026282492349867741025881603897440 401163033562818\ n^{43}$$

$$+\ 11217378240414152444132456358515343571 06916614576509797\ n^{42}$$

$$-\ 22835456743653228369423229758467261165 4234392186582412\ n^{41}$$

$$-\ 58540024121842517755298137846364365590 27649120 2560680123\ n^{40}$$

$$+\ 11936359391805035834929050799031340379 2207216797307942658\ n^{39}$$

$$+\ 27904570926627594905545252942094581696 61425279958282755307\ n^{38}$$

$$-\ 57002777792435693394583410964092297431 15057776713873453272\ n^{37}$$

$$-\ 12096574879195939757905307758697935988 2873189488257994339877\ n^{36}$$

$$+\ 24763177536316236449756449627036794950 8861436753229862133026\ n^{35}$$

$$+\ 47459257021564272964073331772149499818 1035053304 1341887452025\ n^{34}$$

$$-\ 97394831796760169573122308507002679131 2956250283591363 7037076\ n^{33}$$

$$-\ 16762062708345517307469204064862524353 472264823298521482579 2847\ n^{32}$$

$$+\ 34498073734658636310669631214795075498 2574858968806343288622770\ n^{31}$$

$$+\ 52976903313466440492285721913212239543 281707019860778190432 43607\ n^{30}$$

$$-\ 10940361400039874461563840694790398663 63891626294096198137 5109984\ n^{29}$$

$$-\ 14882569916361783289842784687101453870 5708996783290981822 277062561\ n^{28}$$

$$+\ 30859175972727554025841953443681947607 50569098295229256259 29235106\ n^{27}$$

$$+\ 36879754812034300964807821827766109742 955819384494352982729 09304509\ n^{26}$$

$$-\ 76845427221341357332199838999900414246 662207867283935221717 47844124\ n^{25}$$

$$
\begin{aligned}
&-799134926509208476003394452348774393371016449152533269156300564775 06\, n^{24}\\
&+1675115280239758309338988743697449200988695106172350473534318607991 36\, n^{23}\\
&+1498899937949603691860574380921108494881796001089972410732246393775 304\, n^{22}\\
&-3165311403923183214655047636211961909862461512797179868817924648349 744\, n^{21}\\
&-2404750393450973126734077437759885189463121260709851620592549077922 3391\, n^{20}\\
&+5126031927294264574933659639140966569912488672699421228066890620679 6526\, n^{19}\\
&+\\
&3253255524533475803053770341297627147609851484962809711620217136166 91839\, n^{18}\\
&-\\
&7019114241796378063600906646509350952210951837195561546047123334401 80204\, n^{17}\\
&-\\
&3646995891118646731740020949775696952479406010875907934220634356382 241636\, n^{16}\\
&+\\
&7995903206416931269840132564202329000179907205471372023045981046204 663476\, n^{15}\\
&+3314048296358241234281033763130141790915092621503741411804537302743 0132934\, n^{14}-\\
&7427686913358175595546080782680516481848175963554620025913672710106 4929344\, n^{13}-\\
&2371464634265034198559701828110341622860327811857470220214841011899 36708837\, n^{12}+\\
&5485697959865885956674011734488734893905473220070402443021049294809 38347018\, n^{11}+\\
&1283312951833836942201018871963788760305776305671635550168689604681 263632601\, n^{10}-\\
&3115195699654262480069438917376451010002099933350311344639484138843 465612220\, n^{9}-\\
&4933380661075007833364386091660623248176000071403840279686259428680 9917838700\, n^{8}+129828089211544191473571607505889159735221013614271169383646727124 633012\backslash\\
&89620\, n^{7}+120059067630375955684540721394588739233630972169842437136862568920\backslash\\
&11385282630\, n^{6}-36994622447229610284265305029506663820248295795395604365737 1\backslash\\
&86496486071854880\, n^{5}-13333627439799326249482330424447339975713323627675633 1\backslash\\
&49728477664570615409175\, n^{4}+636618773268282627832299658784013437716749430507\backslash\\
&46870665194141825627302673230\, n^{3}-\\
&7153269813502868390040281197554600182421385428787721189870829547628232679065\\
&n^{2}-49355337699822526003149403483292143406832172193171428285452482730370837 3\backslash\\
&15100\, n+24677668849911263001547017416460717034160860965857141427262413651 8\backslash\\
&5418657550\,\big)\,(n+1)^2\, n^2\, /\, 125174280
\end{aligned}
$$

$$s_{n,92} := (2\,n + 1)\,(n + 1)\,(1031415\,n^{90} + 46413675\,n^{89} + 665262675\,n^{88} - 1021100850\,n^{87}$$
$$- 99182929230\,n^{86} + 149284944270\,n^{85} + 18544507521540\,n^{84} - 27891403754445\,n^{83}$$
$$- 3468794819815479\,n^{82} + 5217137931600441\,n^{81} + 625214449498384053\,n^{80}$$
$$- 940430243213376300\,n^{79} - 10765426841599243916\,n^{78} + 16195485538400554 4024\,n^{77}$$
$$+ 176585896046159797 03222\,n^{76} - 265686861834615972331845\,n^{75}$$
$$- 275461028509133102616 2495\,n^{74} + 4145199858554304525409665\,n^{73}$$
$$+ 40807739618988437431 1543185\,n^{72} - 61418869421410371373001961 0\,n^{71}$$
$$- 5733188405009249645391 5899670\,n^{70} + 86304920422245796537738859 310\,n^{69}$$
$$+ 762763383887284710908342 9822320\,n^{68} - 114846032185203935618940141 63135\,n^{67}$$
$$- 959530335119363729829225461696045\,n^{66}$$
$$+ 1445037804288305791524785199625635\,n^{65}$$
$$+ 11394628433791983729874864908009 0175\,n^{64}$$
$$- 17164194540902390884388536621994808 0\,n^{63}$$
$$- 127517486581500598005393418716524896 64\,n^{62}$$
$$+ 1921344395992960165523095549058870853 6\,n^{61}$$
$$+ 1342385278694495732640140618959863566 428\,n^{60}$$
$$- 2023184640021708399787826406185089703910\,n^{59}$$
$$- 1326728386513551847141568191447672851278 90\,n^{58}$$
$$+ 200020850297043631271129141920243472543790\,n^{57}$$
$$+ 1228538108416836836377198720260439827382 5310\,n^{56}$$
$$- 18528082051401074361293545374866719147009860\,n^{55}$$
$$- 106350748269726093190327619061779811836578 5596\,n^{54}$$
$$+ 160452526507159193503556105861413053712218332 4\,n^{53}$$
$$+ 85864573547807602862638731534310501148319132952\,n^{52}$$
$$- 1295991229542472002614758778307728169910397910 90\,n^{51}$$
$$- 6449270174205631341765913951827269569417699802870\,n^{50}$$
$$+ 97387048227855706127796088666562907626220695998 50\,n^{49}$$
$$+ 4494159486472185823175364173576184170873217067167 70\,n^{48}$$
$$- 678993275382220658782694430469755771012293594875080\,n^{47}$$
$$- 2897028832620321418345760527074538271861450876292644 0\,n^{46}$$
$$+ 4379492912699593160457755121352951963427909941827200\,n^{45}$$
$$+ 17220424780766720934035843852492002657087991337566662 0\,n^{44}$$
$$- 26049611816785061059076654553447687448803382503441353 0\,n^{43}$$
$$- 94064004446086982382801886338374961269585657169162704750\,n^{42}$$
$$+ 142398487259969726627156662235279680341622502666261263890\,n^{41}$$

$+\ 4703909156788431857287770046333298265278573194950278438290\ n^{40}$

$-\ 712706297881263264924523340061758723808867104375854828 9380\ n^{39}$

$-\ 21446860966319371237191686495233723838674170742700282 4678332\ n^{38}$

$+\ 325266445984196884882497914128814786199156896662383511 162188\ n^{37}$

$+\ 887512327786558050992760492020330058044942051648140883 8548784\ n^{36}$

$-\ 134753181397904692073326563373693582637737092230533050 13404270\ n^{35}$

$-\ 331682248637784382885593159725114904954178032005771486 55255978\ n^{34}$

$+\ 504261032026571808932455302127451914875013559412392375 489586102\ n^{33}$

$+\ 111326552367482771342411818658030568890346497883158349 53979708526\ n^{32}$

$-\ 169511133711357016058280004497683112908984814621799486 18714355840\ n^{31}$

$-\ 333512856958056680490167996406243755494258418816241672 476247858200\ n^{30}$

$+\ 508744842122652871538165994834249788886882368955452480 23728965220\ n^{29}$

$+\ 885594554449997212193761355553033089879322716156932620 2767872262150\ n^{28}$

$-\ 135382907378112846186755033307126212426332819268317155 45663672875835\ n^{27}$

$-\ 206790174054705165077991878508306724635645468993757015 494953782738625\ n^{26}$

$+\ 316954406450963389926325569427816397574784844454051381 015262510545855\ n^{25}$

$+\ 420787717531156426365048181666847777150533261993979540 2369225501556935\ n^{24}$

$-\ 647029296619282809043888550971662485604539135213671879 4061469507608330\ n^{23}$

$-\ 738379187267928530896755453026260031606022170428850074 2670452673022134\ n^{22}$

$+\ 113992024573285693679732760708797317168926021240395870 6036413763337366\ n^{21}$

$+\ 110366116806801192515834132578366524691162408379083942 19826371887656715268\ n^{20}$

$-\ 171248853330682172421498626710937736225882426748278926 5042576038366741585\ n^{19}$

$-\ 138472537396766137477032834077218893428519360035500664 8353487571151808 6659\ n^{18}$

$+\ 216271248761683314836624182451375226954073161390664943 5782360158646050 0781\ n^{17}$

$+\ 143257318014190166734037804736159487810983893286062114 926565056729949023033\ n^{16}$

$-\ 225699539459369415842887916226807993064179497998626419 5687593858881537 84940\ n^{15}$

$-\ 119516524219868558427314249795437873771958145753021152 1859008676397853754028\ n^{14}$

$+\ 190559763302771308433115770504497210311146193529463049 2572892707540857523512$

$$n^{13} +$$

$$781272474934617351993513456084708862515653721777301440080181748041769759540 6$$
$$n^{12} - 1267188594053311682206828069379311898929053679430683684748917257439697 5 \backslash$$
$$154865 \; n^{11} - 3848333320270772135659678097030093311885930275529171273578926862 \backslash$$
$$2228676284867 \; n^{10} + 6406094277432814044592931180234795917293422253009098752 74 \backslash$$
$$28489220541502004733 \; n^{9} + 13489199049613829226681165524018747088181328080240 9 \backslash$$
$$8544354325428552138066474 69 \; n^{8} - 234368457131371508623182138761455185909187 03 \backslash$$
$$24686602754168630588930914609735 70 \; n^{7} - 30587620652557600132017708817259167 38 \backslash$$
$$3775949937649486640463500091796547358110 2 \; n^{6} + 5759985383540975629185670163 9 \backslash$$
$$6151037112327652990724373153840308234939408584 38 \; n^{5} + 364940073835281984843 60 \backslash$$
$$0787270223565731697741978873389065044716931374809317824 \; n^{4} - 8354093799299478 \backslash$$
$$5541132953172514290045316299561784630225525909080880918440595 5 \; n^{3} - 740824637 \backslash$$
$$44356078570983944604546127812746353994149818164132428621374918039161 \; n^{2} + 528 \backslash$$
$$8283855815080455621406827693906419457010288001478783738281883346696926171 9$$
$$n - 1762761285271693485207135609231302139819003429333826261246093961121556 56 \backslash$$
$$420573 ) \; n \; / \; 191843190$$

$s_{n,\,93} := ( 43890 \; n^{90} + 1975050 \; n^{89} + 27979875 \; n^{88} - 57934800 \; n^{87} - 4373530238 \; n^{86}$

$\qquad + 8804995276 \; n^{85} + 837620061201 \; n^{84} - 1684045117678 \; n^{83} - 160323468043816 \; n^{82}$

$\qquad + 322330981205310 \; n^{81} + 2957840504960671 \; n^{80} - 5947801991126652 \; n^{79}$

$\qquad - 5215717669287103286 \; n^{78} + 10490913350565333224 \; n^{77} + 876618574681715857453 \; n^{76}$

$\qquad - 1763728062713997048130 \; n^{75} - 1402022026051997214949 70 \; n^{74}$

$\qquad + 2821681332731134400380 70 \; n^{73} + 21308386506331086819659695 \; n^{72}$

$\qquad - 42898941145935287079357460 \; n^{71} - 307330122817763401738284518 0 \; n^{70}$

$\qquad + 618950139750120332184504782 0 \; n^{69} + 42005125829324020268830831406 5 \; n^{68}$

$\qquad - 846292017983981608698461675950 \; n^{67} - 543242342524800139348642288495 70 \; n^{66}$

$\qquad + 10949476052294400947842691937500 0 \; n^{65}$

$\qquad + 663733917654733400210278769296295 5 \; n^{64}$

$\qquad - 1338417311361761201368400230530100 0 \; n^{63}$

$\qquad - 7648531404701991248326617889273528 76 \; n^{62}$

$\qquad + 15430904540540158616790075801600067 52 \; n^{61}$

$\qquad + 829805203764443614778944829573523633 97 \; n^{60}$

$\qquad - 167504131206942738817467973494864733546 \; n^{59}$

$\qquad - 8459999546481300783682861016262138244727 \; n^{58}$

$\qquad + 17087503224169544306183190006019141223000 \; n^{57}$

$\qquad + 8088900557840259600594233147317158451371 97 \; n^{56}$

$\qquad - 163486761479222146442502981946945083149739 4 \; n^{55}$

$$-7237756277325326298040457161226150369589 3287\, n^{54}$$

$$+14638999316129874742523417304399245822 3283968\, n^{53}$$

$$+60467124124553957890542530131616192940 22147861\, n^{52}$$

$$-12239814818072090326433740199367231046 267579690\, n^{51}$$

$$-47051053254394812926454850615563260492 9642538335\, n^{50}$$

$$+95326087990596834885553075251063244090 5552656360\, n^{49}$$

$$+34010094218446943936616784976648922988 237457319085\, n^{48}$$

$$-68973449316799856222089100705808478417 380467294530\, n^{47}$$

$$-22771888167810821019302688801867285394 80619586880479\, n^{46}$$

$$+46233510828789640600826286610792655573 78619641055488\, n^{45}$$

$$+14080129728469518595832519976601628094 6116257297014853\, n^{44}$$

$$-28622594565226933597673302639311182744 96111134235085194\, n^{43}$$

$$-80127341868120178613928236973191036215 67091201872626355\, n^{42}$$

$$+16311694319276305058762380421031319070 58379537980337904\, n^{41}$$

$$+41815978968056331283815025344771376594 9146242340920039797\, n^{40}$$

$$-85263127368040293073506288731645885096 88876278219820417498\, n^{39}$$

$$-19932635294988180824429055910252611312 7004820092503755194635\, n^{38}$$

$$+40717901863656764579593174707821681476 36984029672057080 6768\, n^{37}$$

$$+86407569576870254977178204238664465791 5439302191957502876389\, n^{36}$$

$$-17688692934010618641231572594811109973 072488444680635576559546\, n^{35}$$

$$-33900828078286294236461730472736409779 18027470513767600 8292875\, n^{34}$$

$$+69570525449973650337046618204953930556 6779785495598759 3145296\, n^{33}$$

$$+11973381838971437331398057341748262171 2334037075457402 7809386693\, n^{32}$$

$$-24642468932442611166166580865546063648 0234853936410404 3211918682\, n^{31}$$

$$-37842161973454044362469949107888518662 3686804820679751 09306671411\, n^{30}$$

$$+78148570840152349841556556302331643689 5397095035000542 61825261504\, n^{29}$$

$$+10630833176937647308043291552005481943 3034194967436587 7968823871013\, n^{28}$$

$$-22043152062276818114502148667034280323 5022236088522318 10199473003530\, n^{27}$$

$$-26343737890461348567242076318675090936 3230651723992256 85845258653460\, n^{26}$$

$$+54891790987150378945934367504053609904 9963539536506831 81889990310450\, n^{25}$$

$$+57083354133911155724694108867571772906 1999504256684259 14157046033475\, n^{24}$$

$$-11965588736653734934398165448554890680 28996254804987535 01204082377400\, n^{23}$$

$$-10706857269150874277060943792636539832 83749017329572227 321344798645 6546\, n^{22}$$

$$+22610273411967122047561704130128568733 70397660139643208 1437100055290492\, n^{21}$$

$$+$$

$$17177477014132519110829997324046864291 36509971784271177 6854461536743937\, n^{20}$$

$$- 366159813694617504264161650611065854561005970958250655635146023128778366\, n^{19}$$

$$- 232384708807877867498616870464881496137019396177466275890393583805006746 0\, n^{18}$$

$$+ 5013853989852174854236499059908695777301393894507576173443017699228913286\, n^{17}$$

$$+ 2605101480009530038364912866483238884330118607133687505632450469036838706 3\, n^{16}$$

$$- 5711588359004277562153475638957347346390376603718132628609202707996568741 2\, n^{15}$$

$$- 2367272237045974435117919425360212160066844890825564655394438342782392166 66\, n^{14}$$

$$+ 5305703309992376626451186414616159054772727442022942573649796956364441207 44\, n^{13}$$

$$+ 1693971206153222206795379626972128757793123848665431545492397493887865226 093\, n^{12}$$

$$- 3918512743305682076235877895405873421063520441533157348349774683412174572 930\, n^{11}$$

$$- 9166888502065946577104090902036686855728037062132105650684882053810741857 222\, n^{10} + 2225228974743757523044405969947924713251959456579736864971953879103365 8\backslash$$
$$287374\, n^9 + 3524288835928301020201838232259784028461716900383752813747268526 3\backslash$$
$$380665786519\, n^8 - 9273806646600359563448082434467492770175393257347242492466 4\backslash$$
$$909317794989860412\, n^7 - 8575991421710146618167452111112027178536266692820481 9\backslash$$
$$2130281937597889907527 84\, n^6 + 2642578949002065279978298665669154712724792664 2\backslash$$
$$98820633507212968373729713659 80\, n^5 + 9524401346847271949666343579896452143023\backslash$$
$$77041880579984294872314388593599282 89\, n^4 - 4547459218371519669911567381648445\backslash$$
$$14132954674805998060209695759715092051222558\, n^3 + 5109683239140663497486480 81\backslash$$
$$59292043084576994469616403980238483745390369190706\, n^2 + 35255225705433869704 1\backslash$$
$$42712184626042796380068586676525224921879222431131284114 6\, n - 176276128527169\backslash$$
$$3485207135609231302139819003429333826261246093961121556564205 73\,)\,(n+1)^2\, n^2\, /$$
$$4125660$$

$$s_{n,\,94} := (2\, n + 1)\, (n + 1)\, (1155\, n^{92} + 53130\, n^{91} + 779240\, n^{90} - 1195425\, n^{89} - 121162965\, n^{88}$$
$$+ 182342160\, n^{87} + 23687532330\, n^{86} - 35622469575\, n^{85} - 4638120930135\, n^{84}$$
$$+ 6974992629990\, n^{83} + 876038306794420\, n^{82} - 1317544956506625\, n^{81}$$
$$- 158253745476601005\, n^{80} + 238039390693154820\, n^{79} + 27264804938351636810\, n^{78}$$
$$- 41016227102874032625\, n^{77} - 4472773889534674887625\, n^{76}$$

$$+ 6729668947853449347750 \, n^{75} + 69774320073961424063780 0 \, n^{74}$$

$$- 10499796355833480856305 75 \, n^{73} - 103366926279100568482603051 \, n^{72}$$

$$+ 1555753792364425267667198 64 \, n^{71} + 145223256889446459226739465 42 \, n^{70}$$

$$- 2186127649303519014739427 9745 \, n^{69} - 19321018901017764725843597 75041 \, n^{68}$$

$$+ 29090834733991823039502368 02434 \, n^{67} + 243051862812515230980539915 367332 \, n^{66}$$

$$- 36603233595547243762278499 1452215 \, n^{65}$$

$$- 2886293069317778831942650313 6363251 \, n^{64}$$

$$+ 4347741220774441869795114720 0270984 \, n^{63}$$

$$+ 3230055660268902807513549466 488737612 \, n^{62}$$

$$- 48668221965072264206192997733 33241910 \, n^{61}$$

$$- 34003016312738717922679220921 1915749926 \, n^{60}$$

$$+ 51247865578933438205049796370 4540245844 \, n^{59}$$

$$+ 336064226177929603726076013609 95540511872 \, n^{58}$$

$$- 506658732545841077499366510233 45580890730 \, n^{57}$$

$$- 3111923382234379638442551702809 406920688946 \, n^{56}$$

$$+ 4693218009978861511538795879725 783171478784 \, n^{55}$$

$$+ 26938959239565503884026030955976 9656193162372 \, n^{54}$$

$$- 40643099759847198901615986227951 7375875482950 \, n^{53}$$

$$- 217497505619654338144373327950517 65479199518118 \, n^{52}$$

$$+ 328278413417473867161640791237174 06906737018652 \, n^{51}$$

$$+ 1633619219198311283298222452217353 322620952726056 \, n^{50}$$

$$- 2466842749468340618305415717887888 687384797598410 \, n^{49}$$

$$- 11383838966168266777486894336602667 5656073626247298 \, n^{48}$$

$$+ 17199100586725817197145612290798395 7827802838170152 \, n^{47}$$

$$+ 733825975917648631963230052242247023 3810458933023796 \, n^{46}$$

$$- 110933851416983585654341788450876973 29629589818620770 \, n^{45}$$

$$- 436198455402778011023790655947798359 09693941990059041 8 \, n^{44}$$

$$+ 659844375675016195818403073344241387 310223924760196012 \, n^{43}$$

$$+ 2382669067153908166146663463976072045 7957952612532672816 \, n^{42}$$

$$- 3606995819514613059010915349631320138 0592040881179107230 \, n^{41}$$

$$- 11915141087795865875653969958799347005 23735584946654419158 \, n^{40}$$

$$+ 18053061422669529466431500705680586514 75899397860571182352 \, n^{39}$$

$$+ 54325533463003260099826469483065648248 3820736210573649052 76 \, n^{38}$$

$$- 82390953265638366623061279259882501698 3110601305163329490 90 \, n^{37}$$

$$- 22480949887125130169702099664486313225 94630362788279161140258 \, n^{36}$$

$$+ 341333795970158870876684558930288823474110107424767690818493\,2\ n^{35}$$

$$+ 840160950628312207140805960662968527994218589591677814471700\,56\ n^{34}$$

$$- 127730811574097625425504316894096723316503338975875510624847\,550\ n^{33}$$

$$- 281993451417351366240188807603649656819995857744731413678067\,9566\ n^{32}$$

$$+ 429376717704731930631558427250179321395818953565890896048344\,3124\ n^{31}$$

$$+ 844797935673199737902860766214869950576363834851821126718445\,48022\ n^{30}$$

$$- 128866573939503620338586907068481389193433669995602623488008\,543595\ n^{29}$$

$$- 224323721210206160955969629842426325298811691581606804273902\,01812611\ n^{28}$$

$$+ 342928910512284422450883790117063557407889220872190337585253\,5490714\ n^{27}$$

$$+ 523805630020694340000482936565208211908644781759088484913917\,21653072\ n^{26}$$

$$- 802854890556655731123268594353665495733361633682242244250138\,50224965\ n^{25}$$

$$- 106586774006036133017353717190408097881260232111710224661925\,8335968521\ n^{24}$$

$$+ 163894435461837478181646918757380474300557156335976548214139\,4429065264\ n^{23}$$

$$+ 187033633077134745849983085754002574297200898363113282667490\,58890473202\ n^{22}$$

$$- 288745171388793992684056974568872885160829205361468751411942\,85550242435\ n^{21}$$

$$-$$
$$279560769568003778868824524692638546764491632213145890377819\,342725403411\ n^{20}$$

$$+$$
$$433778412921445367937439635767401464404778908587792273137326\,156863226334\ n^{19}$$

$$+$$
$$350755034774955662528253353250485581921957506019838032647180\,1894493593892\ n^{18}$$

$$-$$
$$547821472808505762189252011664098446103175204459146662627636\,5920172004005\ n^{17}$$

$$- 362875025665613621918041085189667109669030588165780748363612\,70981173291521$$
$$n^{16} +$$
$$571703612138845720986524228367705586808704642471628455676800\,89431845939284$$
$$n^{15} +$$
$$302738892469380306561721209321877422276996545976877650637677\,126472599244322$$
$$n^{14} -$$
$$482693519311012745891908025401201412755930051088897898740355\,734424821836125$$
$$n^{13} -$$
$$197898629768885245490640974420518741534826221244392344419031\,0616056937830805$$
$$n^{12} +$$
$$320982620618878505530556862900838182940035834421033411565564\,3791297817664270$$
$$n^{11} +$$
$$974794217650205904561070016311303013548288610400264512288608\,9773319673285160$$
$$n^{10} - 162268263678474810960688345591737361179245083281091347421569\,56555628418\backslash$$

$$759875\, n^9 - 341685403523492731624788681241305215509122630760449252038633662195\backslash$$
$$80459568015\, n^8 + 593662237124476502917527194657826503853306487781219551768735\backslash$$
$$27607184898731960\, n^7 + 774793482329988026187905348138566163130281673739989295\backslash$$
$$91594067075073018327030\, n^6 - 14590213420572202907406216195367624973461954949\backslash$$
$$51608164758278644162019768565255\, n^5 - 92440400598782278952153893180496956346975\backslash$$
$$0173752502192299793024143924327209485\, n^4 + 21161166800103443296526192074758355\backslash$$
$$593877650353763336966876034684239874792424905\, n^3 + 1876530728429683029668681095\backslash$$
$$33214481606976545446632692295524645476909613967220\, n^2 - 133953794926962461927\backslash$$
$$6611767736135021043473358581158867870887024273581605720755\, n + 4465126497565415\backslash$$
$$5397588705892453783403478244528603862892902956747578605352402555) \, n \, / \, 219450$$

$$s_{n,\,95} := (\,2310\, n^{92} + 106260\, n^{91} + 1540770\, n^{90} - 3187800\, n^{89} - 250956090\, n^{88}$$
$$+ 505099980\, n^{87} + 50228412850\, n^{86} - 100961925680\, n^{85} - 10058582819670\, n^{84}$$
$$+ 20218127565020\, n^{83} + 1943609457604430\, n^{82} - 3907437042773880\, n^{81}$$
$$- 359367692533116630\, n^{80} + 722642822109007140\, n^{79} + 6340439342779545795055\, n^{78}$$
$$- 12753142967769923040\, n^{77} - 1065793952059993989544855\, n^{76}$$
$$+ 21443410470877579713936\, n^{75} + 170463589941977297053677655\, n^{74}$$
$$- 3430715209310423520787488\, n^{73} - 259078213814539852642384992\, n^{72}$$
$$+ 52158714283383901288055557472\, n^{71} + 3736683705157929757189455604855\, n^{70}$$
$$- 7525526124599698527259466955688\, n^{69} - 510721002631120381994434314061255\, n^{68}$$
$$+ 1028967531386840462516128095555792\, n^{67} + 6605034582398903704329536871136685\, n^{66}$$
$$- 13312965917936491454910686555178128\, n^{65}$$
$$- 807003665542895032146143288450620925\, n^{64}$$
$$+ 1627320297003726555747197263453023125\, n^{63}$$
$$+ 92994990365613421306855576143706348268\, n^{62}$$
$$- 187617301028230569169458720137579988485\, n^{61}$$
$$- 10089221446069559306586869845322040825275\, n^{60}$$
$$+ 20366060193167349182343198410781661639025\, n^{59}$$
$$+ 102861260103124252616213125089048572675523\, n^{58}$$
$$- 2077591262255652401506605700191753115149485\, n^{57}$$
$$- 983492374623554359593458119677526601930676355\, n^{56}$$
$$+ 1987760661869665243201982296356970735012847455\, n^{55}$$
$$+ 8800056395113453271346995505345857905964106155\, n^{54}$$
$$- 17798888856413873067014189240327412885542949704\, n^{53}$$
$$- 73519207052795362525873950046553069093223064431\, n^{52}$$
$$+ 14881830299123211235844931901713887947498907855666\, n^{51}$$
$$+ 57207220888212632141298638613648212665557319556099\, n^{50}$$

$$- 115902624806337585406181770417467814125896281907764 \; n^{49}$$

$$- 41351316023999444666462750011663986899832703406153 1 \; n^{48}$$

$$+ 83861658296062265186987317727502651940924369631382 6 \; n^{47}$$

$$+ 276872959552014329587407540483560889078404849297105 19 \; n^{46}$$

$$- 562132084933634885693513812739872043350902135557348 64 \; n^{45}$$

$$- 171193849191142367866854108731725636013557435846627 9491 \; n^{44}$$

$$+ 348009019231621084590643355590849992460623893048829 3846 \; n^{43}$$

$$+ 974231654423059944724863268665048479239658184024191 03799 \; n^{42}$$

$$- 198326421076928199790879087288918195772537857533265 01444 \; n^{41}$$

$$- 508421338104788047177724775434477380488352135403935 9163871 \; n^{40}$$

$$+ 103667531831726891433453745959784658055395805838140 44829186 \; n^{39}$$

$$+ 242351784143908884932553909264000499909617182702410 787079099 \; n^{38}$$

$$- 495070321470990459008453193123979465624773945988635 618987384 \; n^{37}$$

$$- 105059006702233788638191130458812683735460176883580 59444038355 \; n^{36}$$

$$+ 215068716619177481866466792848865162127168093227047 54507064094 \; n^{35}$$

$$+ 412184643281683952901030644414770313372463702061787 239586342967 \; n^{34}$$

$$- 845876158225285653988707968114427142957644213446279 233679750028 \; n^{33}$$

$$- 145578866415150832592335040593565465702571894428041 27502863460847 \; n^{32}$$

$$+ 299616494412554521724557160868275202834720230990545 34239406671722 \; n^{31}$$

$$+ 460105517128251402925005416782951971477988090067924 154802945399803 \; n^{30}$$

$$- 950172683697758258022466549652731463239448203234902 843845297471328 \; n^{29}$$

$$- 129255432070986772059146740417749400804957870579994 03996923018581102 \; n^{28}$$

$$+ 268012590978951126698518146332026116242310223192337 10837691334633532 \; n^{27}$$

$$+ 320301444550055368712139398733834315551617938184674 67337866842004198 \; n^{26}$$

$$- 667404148207906186412279694379969474734554609956168 64551342501864192 8 \; n^{25}$$

$$- 694050360854713583496324870582569860261524795960773 5838101301950865022 \; n^{24}$$

$$+ 145484113653021778563387771060031366679965050529171 6403217160289203719 72 \; n^{23}$$

$$+ 130179774191992178960214879262813481339833302131731 596363788171121223878 \; n^{22}$$

$$- 274907959749286535776768535631658329359631654792634 8330492923711628197 28 \; n^{21}$$

$$- 208853076367314823123814914072165162921061773401462 41915093361639828122 \; n^{20}$$

$$+ 445196948709558299825306681707496158778086712947292 731687947909444247597 2 \; n^{19}$$

$$+ 282545925081530174733148290961831211149826807756175 276022409188099248981 78$$

$n^{18}-$

$6096115450340161794488272500944120381774622868070798252136131671429227232 8$
$n^{17}-$

$316742358555556816578764949582988007068362857206515974839901506957511559114$
$n^{16}+$

$6944458716145332511024126241754172179544719430937399320116433062931539055 6$
$n^{15}+$

$287825790073371513958778278728601922251782657570253513861750180416550191202$
$n^{14}-$

$645096167308196353027797819874745566299012509450424695992466469146241577296 0$
$n^{13}-2059622011961796512459276052154371777666148457976273719779648865098617 0\backslash$
$199260\ n^{12}+476434019123178937794634992418348912163130942540297213555176419 9\backslash$
$3434756171480\ n^{11}+11145599920159900592382596639704663255879130769359382579 1\backslash$
$095835323379709735100\ n^{10}-27055540031551590562711543203592815633389570964 12\backslash$
$17372937709312640194175641680\ n^{9}-428502139624478639210023420426544090451921 \backslash$
$20867047595600616930870480337549570 0\ n^{8}+11275596795644731840471622728890163\backslash$
$37237738126982169284950047930049800926633080\ n^{7}+1042715522105048587439158 32\backslash$
$9697445774322223995344096351001548815825556145736020\ n^{6}-321299072377457035 8\backslash$
$9254789322839078858821861176703619869531455617009132181051 20\ n^{5}-1158028341 5\backslash$
$3668997968064531402762254899580591104041567432096458596030850946396 0\ n^{4}+552\backslash$
$904740684795031828676956033915298387379793975119333559507473362153023703304\backslash$
$0\ n^{3}-62126298459257576830080595179176045826734159614574247920411812797303 45\backslash$
$49363320\ n^{2}-42865214376627987816851576567556320673391147474597083771868384 7\backslash$
$7675461138306400\ n+21432607188313993908425788283778160336695573737298541885\backslash$
$9341923883773056915320 0)\ (n+1)^{2}\ n^{2}\ /\ 221760$

$s_{n,96}:=(2n+1)(n+1)(255255\ n^{94}+11996985\ n^{93}+179954775\ n^{92}-275930655\ n^{91}$
$\qquad-29156672545\ n^{90}+43872974145\ n^{89}+5954386441575\ n^{88}-8953516149435\ n^{87}$
$\qquad-1219197115741605\ n^{86}+1833272431687125\ n^{85}+24105717326656751 5\ n^{84}$
$\qquad-362502396115694835\ n^{83}-4563340749046940246 1\ n^{82}+6863136243376195110 9\ n^{81}$
$\qquad+82480552326362808287 67\ n^{80}-124063985301713022187 05\ n^{79}$
$\qquad-142121092387779555009287 9\ n^{78}+213801958508177897624867 1\ n^{77}$
$\qquad+2331565527101831762126572 33\ n^{76}-350803838857815653807110185\ n^{75}$
$\qquad-363722288404735904837310722 15\ n^{74}+547337451801392935525001634 15\ n^{73}$
$\qquad+538836251909813618564657037086 5\ n^{72}-8109910651237273925246105638005\ n^{71}$
$\qquad-757027414726661962908210319676075\ n^{70}$
$\qquad+1139596077415611581324938532333115\ n^{69}$

$+ 1007176330585729458248346274054244405 \, n^{68}$

$- 1516462476265672245279144103740333165 \, n^{67}$

$- 12669936947021149143399532776279597275 \, n^{66}$

$+ 19080728544345007327363256369606412495 \, n^{65}$

$+ 150458223483844802186577086670320691 9295 \, n^{64}$

$- 226641371652984453646233792823961358519 0 \, n^{63}$

$- 1683780633850935690257155593743928474812 42 \, n^{62}$

$+ 2537003019359052758068045080257090780144 58 \, n^{61}$

$+ 1772527362076425221947414878362405097845 8534 \, n^{60}$

$- 2671476058211433096711462542944893100669 5030 \, n^{59}$

$- 1751853514732420117929026428311077692170 381578 \, n^{58}$

$+ 2641137652389687342377096955181341003758 919882 \, n^{57}$

$+ 1622200011264652274021016011494744752073 00908166 \, n^{56}$

$- 2446505705158926847743409502018023831283 0822190 \, n^{55}$

$- 1404288429197850162230173190011526152722 5614777490 \, n^{54}$

$+ 2118665172322569877583976832527379348249 4837577330 \, n^{53}$

$+ 1133782592733120983739126557393697192063 947989625710 \, n^{52}$

$- 1711267214961294324996609720253182684837 169403227230 \, n^{51}$

$- 8515817358940944256379711615956408121509 0513181516930 \, n^{50}$

$+ 1285928939915948110081939790994727131650 54354473889010 \, n^{49}$

$+ 5934228266919489886726910211166833324181 642676299398790 \, n^{48}$

$- 8965638847375032235594462306299986342854 991191686042690 \, n^{47}$

$- 3825327169711523240410264643504838709536 4087396890 6406910 \, n^{46}$

$+ 5782818948804160021793369276788757996018 88806549202631710 \, n^{45}$

$+ 2273838563362211851954368130900230510970 5868391798647181890 \, n^{44}$

$- 3439671939787338578040519042734289556435 9746990972572088690 \, n^{43}$

$- 1242050434044322310676012062894339368718 2442126175012709098 38 \, n^{42}$

$+ 1880274010765420158904220689555180500859 546192421738192409102 \, n^{41}$

$+ 6211188269411585913741724340653847501974 281998367011611690246 6 \, n^{40}$

$- 9410796104655649878557797545458530278004 4003071716043271558250 \, n^{39}$

$- 2831910370919087118913751297712389661587 855110263003722269414742 \, n^{38}$

$+ 4294919536901908927763415934295877143771 804666930363605039901238 \, n^{37}$

$+ 1171898940979905501856321789988880067272 621036677467895234011367 14 \, n^{36}$

$- 1779323009154367797423299764654799486627 790578350853660876216556 90 \, n^{35}$

$- 4379635794917465484654471098058586336396 54224496758756304455428339 0 \, n^{34}$

$$+ 66584198428339166168528716353206194789262028963689240276106422529930\ n^{33}$$

$$+ 146990497459014684044339118755106243552462072786384381145660018888140\ n^{32}$$

$$- 223827784540269160915077303630926246272332412366142119185654323958675\ n^{31}$$

$$- 440380771777093646999846095751810015258643749631799025092351806914580 5\ n^{30}$$

$$+ 671762546892653928545523008809261335201582245066005643597810426569804 5\ n^{29}$$

$$+ 116936665329015569165723596027525291379232324392152458734317038503750 995\ n^{28}$$

$$- 178763810727986623391313009085334243744856397813558716319464609888475 515\ n^{27}$$

$$- 273052191380985659150346174808518527286171710829100230267758820231405 0245\ n^{26}$$

$$+ 418516477607877819895084912667044503116500386134328281217611460841531 3125\ n^{25}$$

$$+ 555621217996992973471159744783727753410641861635291714746005932688490 68835\ n^{24}$$

$$- 854357650875883351201493862808943855271787811759653986179889472074812 59815\ n^{23}$$

$$- 974978893824436210544838927009141461157180991420660699587090685630843 688633\ n^{22}$$

$$+ 150518622328044848337733308365415938449936087771897374868963050205000 6162857\ n^{21}$$

$$+ 145730928382122189108477192928167974034167798891443590746981323469146 24\backslash 809031\ n^{20}$$

$$- 226122323689585526079602454810522757973748502725760254863920137 7\backslash 1396940294975\ n^{19}$$

$$- 182843454507030879167112924399363759116709510978704428356\backslash 679500776561476587665\ n^{18}$$

$$+ 285571297945025595054649509339571776573751691604 3\backslash 44655278215258050540685028985\ n^{17}$$

$$+ 189161427973793719555708625702550100333 37\backslash 846870738468753979665876372770584021 55\ n^{16}$$

$$- 2980207068579418590862954140208 0\backslash 37393287552876412942640736057510481185930117 725\ n^{15}$$

$$- 157813343857627163387 10\backslash 2808882707149411536871446456083557807518097251113388 07619\ n^{14}$$

$$+ 2516210511293\backslash 3783803496898402510091108374307155174883854039156469 8282597327029 1\ n^{13}$$

$$+ 103\backslash 16165278245218078123913620195639114777287021084272238253235711487 4233813692\backslash 893\ n^{12}$$

$$- 167323531730145163073607153504189632275846458893851525500818113907 2\backslash 25480707174485\ n^{11}$$

$$- 5081459266949583318037919861192599309875963753184508372 7\backslash 6649743043262903898829739\ n^{10}$$

$$+ 845880655907510079242491555930984712619317792\backslash 4246020186653836715185070962018318 51\ n^{9}$$

$$+ 1781155887753632286721158674521516 7\backslash 9837656992148395540566135415729395446924322 173\ n^{8}$$

$$- 3094674159584203469702 98\backslash 37897477675538745137784382341178247230717002467184873991 85\ n^{7} - 4038884771237\backslash$$

$567944031200331851677413820426689516909701795119064249052075412993679\ n^6 + 76\backslash$
$056642366484536508982923926513998976678969234944816116050401322237014723631\backslash$
$90111\ n^5 + 481878248514367150329306067299236253925873365646980320014614476169\backslash$
$6146131338437073\ n^4 - 11031005846039734080388737205814243757722048946451945600\backslash$
$6021737208656069933189250665\ n^3 - 9782079424599519280117570949811385221566284 0\backslash$
$144086766159734730656890913926363661 43\ n^2 + 698281483670979493322072502762419\backslash$
$971121045049483898772697107820228616720555495795 47\ n - 2327604945569931644073 5\backslash$
$750092080665704034834982796625756570260676205573518495 79849\ )\ n\ /\ 49519470$

$s_{n,97} := (\ 72930\ n^{94} + 3427710\ n^{93} + 50844365\ n^{92} - 105116440\ n^{91} - 8622064165\ n^{90}$

$+\ 17349244770\ n^{89} + 1801720225305\ n^{88} - 3620789695380\ n^{87} - 377117979821869\ n^{86}$

$+\ 757856749339118\ n^{85} + 76242421912670533\ n^{84} - 153242700574680184\ n^{83}$

$-\ 14764939694533897453\ n^{82} + 29683122089642475090\ n^{81}$

$+\ 2731460553029145298673\ n^{80} - 5492604228147933072436\ n^{79}$

$-\ 481984079335618183595063\ n^{78} + 969460762899384300262562\ n^{77}$

$+\ 81021604372590280998764269\ n^{76} - 163012669508079946297791100\ n^{75}$

$-\ 12958739560578817908113328845\ n^{74} + 26080491790665715762524448790\ n^{73}$

$+\ 1969531395828550455600653144905\ n^{72} - 3965143283447766626963830738600\ n^{71}$

$-\ 284065552726650943683986145756005\ n^{70} + 572096248736749653994936122250610\ n^{69}$

$+\ 38825406070327108608337197807218245\ n^{68}$

$-\ 78222908389390966870669331736687100\ n^{67}$

$-\ 5021198590064230577778593365589053973\ n^{66}$

$+\ 10120620088517852122427856062914795046\ n^{65}$

$+\ 613490459252451156215029463244838121281\ n^{64}$

$-\ 1237101538593420164552486782552591037608\ n^{63}$

$-\ 70695514648567658265675708372915695488226\ n^{62}$

$+\ 1426281308357287366959039035283839820140 60\ n^{61}$

$+\ 76699045838414191499176481739303931815 96181\ n^{60}$

$-\ 1548243729851856703653120025138917034520 6422\ n^{59}$

$-\ 781959296539479537843934540861092274577 529401\ n^{58}$

$+\ 157940103037747764272440028197357371950 0265224\ n^{57}$

$+\ 74765854964703513485523076215356387173442244313\ n^{56}$

$-\ 15111111095978450461377055271268634806638475 3850\ n^{55}$

$-\ 66898712903088418318673162931280936946364295 57085\ n^{54}$

$+\ 13530853691577468168348403138968873737339243868020\ n^{53}$

$+\ 5588987279931390460410381819809654719819966611 22525\ n^{52}$

$-\ 113132830967785556025042476710089981770133256611 3070\ n^{51}$

$$- 4348937409444437347487035696101835446800155768006 7825\ n^{50}$$

$$+ 881100764985666025099911386891376087537044479262 48720\ n^{49}$$

$$+ 3143559194003868220278051029748421788749189020284210145\ n^{48}$$

$$- 637522846450630304306609319818598118625208248849 4669010\ n^{47}$$

$$- 2104809764858941335067364981478127023475112461175 34904193\ n^{46}$$

$$+ 4273371814362945700565390894938113858812745747235 64477396\ n^{45}$$

$$+ 130142895154633187038557595090283264402035594427 09238072901\ n^{44}$$

$$- 264559162123629319777680581075504642662883934601 42040623198\ n^{43}$$

$$- 740618478157715655670748865482272434131062971857 7665166668 09\ n^{42}$$

$$+ 150769287252779424331926578907209533252844337175 675073956816\ n^{41}$$

$$+ 386505854106196542184003062751534029151388815849 57126313337577\ n^{40}$$

$$- 788088636937671026801198783393789011628061775070 89927700631970\ n^{39}$$

$$- 184237710545099474054043364239039451942155966315 66605756990940893\ n^{38}$$

$$+ 376356307459575658376098716312016794000592550382 0301441682513756\ n^{37}$$

$$+ 798666737087781905854927583921829979678978268346 23586022504310221\ n^{36}$$

$$- 163496910492152137754746503947486163875801579173 067473486691134198\ n^{35}$$

$$- 313345970479725051025768920447631970097309029577 31618518807 01346033\ n^{34}$$

$$+ 643041632008665315827012491290012556582198217071 9391177248093826264\ n^{33}$$

$$+ 110670186096718974650573204540480048884669201899 3295142330752822535825\ n^{32}$$

$$- 227770788513524602459416533993860223335160385969 3784196433988588297914\ n^{31}$$

$$- 349775791353540228034138951292771718676978267450 118320135957 3651402142\ n^{30}$$

$$+ 722328661558432916314219555984929459687425734971 7448223625461611 02198\ n^{29}$$

$$+ 982609843966966650806008077390211441576128156095 43298286604112713457261\ n^{28}$$

$$- 203745255408977659324343811037891582912100356954 05834139557077158801 6720\ n^{27}$$

$$- 243495648433358150563897801844936703827017194166 2339820099114060905985973\ n^{26}$$

$$+ 507365822407614067060229984793662565945244424027 8737981593798893399988666\ n^{25}$$

$$+ 527622480430938365142248956169119659707461168625 03458398437487354421022041\ n^{24}$$

$$- 110598154310263813699052091081760557600944677965 28565477846877360224 2032748\ n^{23}$$

$$- 989636764636686740171587657197864168067613605658 0720692164871199427811 57709\ n^{22}$$

$$+ 208987168358363729404222740547748889373617188928 142979321144301348780434 8166$$

$n^{21} + 15877173244726864604501982210168596421780372322109783388722472643003753\backslash$
$371877 \; n^{20} - 33844218173037366503046191825814681737296916533500996570656388829\backslash$
$9495311091920 \; n^{19} - 214793608987659075301210695071582218179449043818921246491\backslash$
$487730011539862529725 \; n^{18} + 463431436148355517105467581968791180961950041713\backslash$
$4348955363184832257503611537 0 \; n^{17} + 24079000358536922490067981213502871486077 5\backslash$
$3949206365056917593647421120076939505 \; n^{16} - 5279231507855740015241430008974 7\backslash$
$2209025127398829864462790553914316481519003 0380 \; n^{15} - 2188074034044001874 9825\backslash$
$828482329014795780962362154515125431760576957659015430 279 \; n^{14} + 4904071218873\backslash$
$5777514893086973632751681813198712607674878769060297080133220890938 \; n^{13} + 156\backslash$
$574066666168585775016005417839309526730167836683092342869892833617011755988\backslash$
$813 \; n^{12} - 362188845521072949064925097809311370735273534385973859564508845964 3\backslash$
$14156732868564 \; n^{11} - 8472971714219226596923321063483489357111423561723489018 8\backslash$
$707367159831671439574172 5 \; n^{10} + 20567831883649182684495893105060092421575582 4\backslash$
$6730671663338656189160947585524352014 \; n^9 + 32575065806493953487201535211725 12\backslash$
$8889567721817881492769256724770600457016615662 17 \; n^8 - 85717963496637089658898\backslash$
$9635285103502007110261030697021719000114328103898884748444 48 \; n^7 - 79268044681 8\backslash$
$5737031570631517265639816057521671184601734509944651324756127076005157 \; n^6 + 2\backslash$
$4425405286035183029031159387382314652186145952676173686209890445930551242 99\backslash$
$9494762 \; n^5 + 880342148685685900019545957865369338766150970516800204538667275 7\backslash$
$138756698695554533 \; n^4 - 420322482597489010294220785446897014275091653630121 77\backslash$
$7769832359602080646403906038 28 \; n^3 + 472288951088492900619601420788838472093 01\backslash$
$9819354845085889243550676013085724628297 1 \; n^2 + 32586469237979043017030050128 9\backslash$
$1293198564876897591527605919836494668780292589803788 6 \; n - 16293234618989521 50\backslash$
$8515025064456465992824384487957638025991824733439014629490189 43 \,) \, (n+1)^2 \, n^2 \, /$
$7147140$

$s_{n,\,98} := (2\,n+1)\,(n+1)\,(3315\,n^{96} + 159120\,n^{95} + 2439840\,n^{94} - 3739320\,n^{93}$

$\quad - 411574488\,n^{92} + 619231392\,n^{91} + 87718713576\,n^{90} - 131887686060\,n^{89}$

$\quad - 18763723162236\,n^{88} + 28211528586384\,n^{87} + 3879587921980992\,n^{86}$

$\quad - 5833487647264680\,n^{85} - 768802459409261784\,n^{84} + 1156120432937525016\,n^{83}$

$\quad + 145618201272210503148\,n^{82} - 219005362124784517230\,n^{81}$

$\quad - 2632346036558124257 1942\,n^{80} + 3959469322943452 6116528\,n^{79}$

$\quad + 4535910902508274315 760864\,n^{78} - 6823663700377128601 699560\,n^{77}$

$\quad - 7441444197551282396818 58440\,n^{76} + 1119628461482880923823 637440\,n^{75}$

$\quad + 1160861539244306948312 37281000\,n^{74} - 174689045113874827087 67740220\,n^{73}$

$\quad - 1719758798314920971229 1256275468\,n^{72} + 2588372649728250830979 1268283312\,n^{71}$

$\quad + 2416141677777607094129 731021285856\,n^{70}$

$$- 363715437991505189534949216607 0440\, n^{69}$$

$$- 3214521350926371315223814473545 9024\, n^{68}$$

$$+ 4839967798289132232312464631862 23756\, n^{67}$$

$$+ 404375896482000767667038805735717 84578\, n^{66}$$

$$- 608983828622145717616714440919507 88745\, n^{65}$$

$$- 480205066380330375497286620902187 7283645\, n^{64}$$

$$+ 723352518713606291834013503557879 1319840\, n^{63}$$

$$+ 537398337380265101233776403270590 513594560\, n^{62}$$

$$- 809714268663965683309834672423675 166051760\, n^{61}$$

$$- 565722896623138449375863714208526 55513887216\, n^{60}$$

$$+ 852632916278027502480344744674908 20853856704\, n^{59}$$

$$+ 559124595788347659298322901618937 0268941627472\, n^{58}$$

$$- 842950058263911626459886076151780 0813839369560\, n^{57}$$

$$- 517744159523401419395220517442453 492662540669432\, n^{56}$$

$$+ 780830989576421687225130206544439 139400730688928\, n^{55}$$

$$+ 448195122337643161166915732619920 77965812178456704\, n^{54}$$

$$- 676196838454346850186499249962603 36518418633029520\, n^{53}$$

$$- 361860012009657113002102103385617 5513139823915689008\, n^{52}$$

$$+ 546171002206757403754085651328239 343796894519004 8272\, n^{51}$$

$$+ 271792298764296108154062117742203 62841185580762663 7736\, n^{50}$$

$$- 410419303157477949249863604869946 6393367681840349 80740\, n^{49}$$

$$- 189397855082550080160190562557811 4446678308934019 3791764\, n^{48}$$

$$+ 286148879139612509986535161861066 90019843018102308 178016\, n^{47}$$

$$+ 122089803146124811112389568965294 7155433906365622 241660608\, n^{46}$$

$$- 184565449114885279218517029257247 4078160781057484 516579920\, n^{45}$$

$$- 725722246152073853811337987857029 2550314520493720 2135732368\, n^{44}$$

$$+ 109781164168385504467793283324840 6252937981979354 5461888512\, n^{43}$$

$$+ 396414963380677814092568587318712 2639417340034311 853869491856\, n^{42}$$

$$- 600111503279435996362242545144310 4271772909150435 053535182040\, n^{41}$$

$$- 198237358393888402905732202842521 8150223783554690 04468634 82232\, n^{40}$$

$$+ 300356595107229784340409516989504 2746694539877788 52599762814368\, n^{39}$$

$$+ 903837408864181148972601635867232 7104731313849870 756958101002144\, n^{38}$$

$$- 137077394305163321267592292965032 4279443169776869 5561737032910400\, n^{37}$$

$$- 374025291599251399306457562770632 3413711176057278 15868123298747496\, n^{36}$$

$$+ 567891807114135265023065958804200 1334538922574607 1583053464576444\, n^{35}$$

$$+\ 13978121303896737317141628308803824387235270193743803432741152737202\ n^{34}$$

$$-\ 2125112785940217360822397544260783664757985141935374094063846139402 5\ n^{33}$$

$$-\ 46916470800843928818176119938155613237937511798020550084426522530412 5\ n^{32}$$

$$+\ 7143726259423600190767537867936381168928526026799851217367170686532 00\ n^{31}$$

$$+\ 1405526883068112340545899664847195945874938380110452571281291583479 8560\ n^{30}$$

$$-\ 2144008955899286511772687186610475824657050200299678113008773228652 4440\ n^{29}$$

$$-\ 37321708228319902907232602109080725428582278663231868537793394694433 9384\ n^{28}$$

$$+\ 57054566820429497616735246756926326055201943094997641863194478655977 1296\ n^{27}$$

$$+\ 87147809364595161390457318203815544649774253609013214539433556509227 73448\ n^{26}$$

$$-\ 13357444238791421696652273964356963327742147756826970390231005869664 045820\ n^{25}$$

$$-\ 17733302831239604736712227892106755179066829849956685853745636565310 6018188\ n^{24}$$

$$+\ 27267826458798978189900955536377980934987352162776377300130005141449 1050192\ n^{23}$$

$$+\ 31117594897805547922697196919214918692838423159107605328339749847114 00997696\ n^{22}$$

$$-\ 48039783669648270793540843155641277086007002346800226857516125027743 47021640\ n^{21}\ -\ 4651173499446613445808655660342731235196614343798580506815266234053 6366\backslash$$
$$533240\ n^{20}\ +\ 7216959167518161522680687706292303238224956527431871894510479976\backslash$$
$$2191723310680\ n^{19}\ +\ 5835663298050486889554841803578821441934259645048183 17357\backslash$$
$$50404882442280259614 0\ n^{18}\ -\ 9114342905451638410466297090682847324812637293943\backslash$$
$$8683550880847311773006554955 0\ n^{17}\ -\ 6037308831260574866622701883970732063160 4\backslash$$
$$42062376444281957696339546024825788582\ n^{16}\ +\ 9511680392163444220457367680490 2\backslash$$
$$4046098129495826185984069094874587790227145764 8\ n^{15}\ +\ 5036797960176162399953 0\backslash$$
$$8040666906577769077028907864718733184919820554699657391 04\ n^{14}\ -\ 8030780959872\backslash$$
$$41581095248899402811068958522018153106377303232123460221560843374 80\ n^{13}\ -\ 329\backslash$$
$$2525141421645664496279560876633788411633258547113379431408686799929056880 10\backslash$$
$$344\ n^{12}\ +\ 5340326760126089287292043791016356217096710896897223257798729091930\backslash$$
$$0043657418425 6\ n^{11}\ +\ 1621807322806622083721122806411595254473256446534860384 7\backslash$$
$$41884172629442622014680968\ n^{10}\ -\ 2699727322216237589946286399168210692564720 2\backslash$$
$$146471517400027627135406645151309113580\ n^{9}\ -\ 5684767918159583590198969818824 24\backslash$$
$$39840798475553918633108976375681601385216793369 56\ n^{8}\ +\ 9877015538374941802 71\backslash$$

$$597927820471322402131440411370836347837709010539858173562224\, n^7 + 12890574446\backslash$$
$$9292945380356614357577564464248665867467492425753167194018110023 88594592\, n^6$$
$$-\ 2427436943956768889718929111754687033083836560032580928203689393360798643\backslash$$
$$2669673000\, n^5 - 1537970947095626581409633330766828149251485671043974097303133\backslash$$
$$695276429553244366 8608\, n^4 + 352067489262182431697391455202758574041914678 6582\backslash$$
$$251610056545239595043651500033 9412\, n^3 + 31220653772189123310064288345523 59069\backslash$$
$$4332430481653247885884745165691447 15205331006\, n^2 - 2228647252893749008137 9216\backslash$$
$$011966467306245598505159245233165437972 8935330308166215\, n + 742282417631249\backslash$$
$$66937930720039888224354151995017197484110551459909429784434360554 05)\, n\ /$$
$$656370$$

$$s_{n,\,99} := (1326\, n^{96} + 63648\, n^{95} + 965328\, n^{94} - 1994304\, n^{93} - 170294865\, n^{92} + 342584034\, n^{91}$$
$$+\ 37119923997\, n^{90} - 74582432028\, n^{89} - 8113039988001\, n^{88} + 16300662408030\, n^{87}$$
$$+\ 171441098329 6141\, n^{86} - 3445122629000312\, n^{85} - 347375482506735657\, n^{84}$$
$$+\ 698196087642471626\, n^{83} + 6730818362904 1401605\, n^{82} - 13531456334572527 4836\, n^{81}$$
$$-\ 1245342158461964518346 3\, n^{80} + 25042157732585015641762\, n^{79}$$
$$+\ 2197559941324149652265639\, n^{78} - 4420162040380884320173040\, n^{77}$$
$$-\ 369413215787782278725183233\, n^{76} + 743246593615945441770539506\, n^{75}$$
$$+\ 59084727416646892470636313421\, n^{74} - 118912701426909730383043166348\, n^{73}$$
$$-\ 8979984837122340899567070636145\, n^{72} + 18078882375671591529517184438638\, n^{71}$$
$$+\ 1295183557130079834455815023690269\, n^{70}$$
$$-\ 260844599663583126044114723181917 6\, n^{69}$$
$$-\ 177022622576981142768207141378897337\, n^{68}$$
$$+\ 3566536911505981167968554299896138 50\, n^{67}$$
$$+\ 22893920240502715215339958964818750037\, n^{66}$$
$$-\ 461444941721560285474767733596271139 24\, n^{65}$$
$$-\ 2797181074997676290021575491582633683644\, n^{64}$$
$$+\ 5640506644167508608590627756524894481212\, n^{63}$$
$$+\ 322332894955794503759991939808501185286370\, n^{62}$$
$$-\ 650306296555756516128574507373527265053952\, n^{61}$$
$$-\ 34970571487231569541461739942469492284033201\, n^{60}$$
$$+\ 70591449271018895590520543923125118331203 54\, n^{59}$$
$$+\ 3565306866812305904721736690628928450853578893\, n^{58}$$
$$-\ 72012051828956307050425254356501694135402781 40\, n^{57}$$
$$-\ 340891421446633281491721481870418928789658962637\, n^{56}$$
$$+\ 688984048076162193688485489176488026992858203414\, n^{55}$$
$$+\ 3050215548968551436848141140909503592899842818020 9\, n^{54}$$

$$- 616932950274471909306513083073665598849897145638832\, n^{53}$$

$$- 254827262954039184642907174816548332968369455334 1785\, n^{52}$$

$$+ 515828554108230883788794804638333219252378821247402\, n^{51}$$

$$+ 198287768660101103149490661394591876752096946876673381\, n^{50}$$

$$- 401733775874310437182770117593822086723446272574594164\, n^{49}$$

$$- 143329112273800009101815047371415343275992918728 29749673\, n^{48}$$

$$+ 290675562306343122575457795918768907419220300182 34093510\, n^{47}$$

$$+ 959678175228561363286183094196761047427987103208 47188053\, n^{46}$$

$$- 194842390668775703882991196798539898559787945065 9928469616\, n^{45}$$

$$- 593380448086864143805010311413385989451234059261 82481725441\, n^{44}$$

$$+ 120624513524060585799831974250662596875844691303 02489192 0498\, n^{43}$$

$$+ 337681533754469802968241003810221524398166422228 90084370045\, n^{42}$$

$$- 687425518861345664516465205045509308483917313576 0805060660588\, n^{41}$$

$$- 176225537802296469994537649649094674629147908112 149736355330109\, n^{40}$$

$$+ 359325330793206396634239951348644442343134989360 060277771320806\, n^{39}$$

$$+ 840023230679275221461535730193782732253745474253 9163307873148097\, n^{38}$$

$$- 171597899443787108258649545552242990874180444744 38386893517617000\, n^{37}$$

$$- 364148366118740322314618846727540271158100861249 58622177063 5526121\, n^{36}$$

$$+ 745456522181859355455102648010304841403619766973 61083043478866 9242\, n^{35}$$

$$+ 142868630783532474166917803566277815533030874084 53946532753374217237\, n^{34}$$

$$- 293191826788883541888386633612658679480643945838 81503895941537103716\, n^{33}$$

$$- 504595541215681183411821507401947129165120776032 206728228484012261540\, n^{32}$$

$$+ 103851026511025072101248167816516012627830594664 8294960352909561626796\, n^{31}$$

$$+ 159478637352192312590142072439794369564896561134 8197050751779291856 8298\, n^{30}$$

$$- 329342377355487132390408961661240340392576181736 1223597538849539876 3392\, n^{29}$$

$$- 448016366022056511088665953118891883785854229436 584491327564661441949969\, n^{28}$$

$$+ 928966969779661735416372802403907801610966077046 78121863051781828266333 0\, n^{27}$$

$$+ 111020702899618784919367170298978039187806462273 565243730029368345969709 09\, n^{26}$$

$$- 231331075497034187192898068621995156391722585317 598299646363914874766051 48\, n^{25}$$

$$- 240567004051059929399896851816087657711568711241 489226437747204769876180 33\, n^{24}$$

$$+ 504267115651823277519083510494374831062309681014 73828284018583244145184 1214$$

$$n^{23} +$$

$$45122025766033232842612265777606115454535931170622015270146290232965217262 05$$

$$n^{22} -$$

$$95286722688584698460415366660155979219694959151391413368694438790344952936 24$$

$$n^{21} - 72391229372258213060482858630417345668556722340201350237214230602464508\backslash$$

$$750377 \, n^{20} + 1543111310133748959670072539268502892590829405955418418112979050\backslash$$

$$83963512794378 \, n^{19} + 97934142156475941124765238286863137913369491514421492705\backslash$$

$$295463370261421824442 1 \, n^{18} - 21129939741428937184623120196641130475264727708 8\backslash$$

$$397169591720172489191949283220 \, n^{17} - 109787076776305645417528370583288414636\backslash$$

$$519987218075022883028543383450473653076 55 \, n^{16} + 2407040932940402280196798613 6\backslash$$

$$32175597483047021449897627252291584917928667989853 0 \, n^{15} + 997642129618814142 6\backslash$$

$$42735417755533660402120011750118597285371713757627205833 57895 \, n^{14} - 22359883 5\backslash$$

$$253166851330515069687428528055254472564522695729597258600704727846614320 \, n^{13}$$

$$- 713892139303972690638618124812545585670922498876676614774258757194349276 9\backslash$$

$$69246785 \, n^{12} + 1651383113861112232607751319312519699397099470317875925278114 7\backslash$$

$$7298940328178510 7890 \, n^{11} + 3863211853737328778537730906245269037211527610075 4\backslash$$

$$62491189809502869676496347465805 \, n^{10} - 937780682133576978968321331803057773 8\backslash$$

$$2015469046880090765773377828756274480039500 \, n^9 - 148524490113345523100307646\backslash$$

$$9054169480853888374278839643764380711926687540617294 2945 \, n^8 + 390827048440048\backslash$$

$$7440974474251288644739089792217604559378294534801726250708682592539 0 \, n^7 + 361\backslash$$

$$41894504807779543915573337074165778695599269368704238346129080092067347758 1\backslash$$

$$69965 \, n^6 - 111366649385362043349757588918703477894828912071478300225963760617 7\backslash$$

$$44664178234226532 0 \, n^5 - 40138788831782672003075991293553151709917216927578101\backslash$$

$$09306601419257563218434304338 5 \, n^4 + 1916440715171857775037278717741410823681 2\backslash$$

$$35545699392044457696345625979061510283520 90 \, n^3 - 2153373993546792181393321584\backslash$$

$$71823168299097822677721181123333573718691686411536219 95 \, n^2 - 1485764835262499\backslash$$

$$3387586144007977644870830399003439496822110291981885956886872110810 0 \, n + 7428\backslash$$

$$82417631249669379307200398882243541519950171974841105514599094297844343605 5\backslash$$

$$4050 ) \, (n + 1)^2 \, n^2 \, / \, 132600$$

$$s_{n, 100} := (2 \, n + 1) \, (n + 1) \, ( \, 36465 \, n^{98} + 1786785 \, n^{97} + 27992965 \, n^{96} - 42882840 \, n^{95}$$

$$- 4912467560 \, n^{94} + 7390142760 \, n^{93} + 1091693643040 \, n^{92} - 1641235535940 \, n^{91}$$

$$- 243731019191660 \, n^{90} + 366417146555460 \, n^{89} + 52646885033883140 \, n^{88}$$

$$- 7915536124102440 \, n^{87} - 1090967097436967560 \, n^{86} + 1640452741421 7502560 \, n^{85}$$

$$+ 216317145557513894894 0 \, n^{84} - 3252959447069817174690 \, n^{83}$$

$$- 40977982150585232495743 0 \, n^{82} + 6162962119823133960234 90 \, n^{81}$$

$$+ 74078548541948785665560730 \, n^{80} - 11142597091891433519 6352840 \, n^{79}$$

$$- 12764906296261161728292446840\, n^{78} + 1920307242985119976003 6846680\, n^{77}$$

$$+ 2094166705521363344351540460400\, n^{76} - 315085159449697061640 7329113940\, n^{75}$$

$$- 3266891190926262472668914148 6940\, n^{74} + 4916091044361878562085407 86287380\, n^{73}$$

$$+ 4839737908476596639663992318 3338580\, n^{72}$$

$$- 7284187317936704352306415516 8151560\, n^{71}$$

$$- 6799495871402762143259544027803046800\, n^{70}$$

$$+ 10235664743693826736650848119288645980\, n^{69}$$

$$+ 9046292703712916117465770 5626468495430\, n^{68}$$

$$- 13620617379287843309881910104 99347066135\, n^{67}$$

$$- 1137992979335732858688917 6531536435 2720125\, n^{66}$$

$$+ 171379977769324320968831743478296202613255\, n^{65}$$

$$+ 135139117055554193190752 1014869149304077555\, n^{64}$$

$$- 20356556744717975058345697 394042872057422960\, n^{63}$$

$$- 151234418458224992552167 88523525371393988 89360\, n^{62}$$

$$+ 22786945552457338758116911272582714512 7045520\, n^{61}$$

$$+ 1592055042430083008318361255565984151 3277078400\, n^{60}$$

$$- 239947603642135318185660033897 1176758424 79140360\, n^{59}$$

$$- 157348613144687599970620677009797554075 37972261720\, n^{58}$$

$$+ 23722265773524207654685931568418191949 228197962760\, n^{57}$$

$$+ 145703347837994083470384791048501114415 3234889555720\, n^{56}$$

$$- 21974113504566733558831148315117258122044664 33314960\, n^{55}$$

$$- 126130886477598047803142298690225613802274941514776720\, n^{54}$$

$$+ 19029503539162540838265500545109428360951464548 8822560\, n^{53}$$

$$+ 10183449533659608151394867193581475904463225005798387720\, n^{52}$$

$$- 15370321818185224931283628293077609984995948 1441992860\, n^{51}$$

$$- 764876766220237177859558917546216065303101615557369034420\, n^{50}$$

$$+ 115500031023944837925498019046587297845390222075177454 8060\, n^{49}$$

$$+ 53300266263328491874753411533014873578846789858820803568620\, n^{48}$$

$$- 8052789955011246200175760739475524685749713589860709262 6960\, n^{47}$$

$$- 3435846204746897126526870936166116324173489312315613897054320\, n^{46}$$

$$+ 519400332568954019207911852079465521096889825364227243918 94960\, n^{45}$$

$$+ 204232455200018192822908974332679519513176874627461457402466720\, n^{44}$$

$$- 308945699428474901947590541029925532460980320940354829964 7560\, n^{43}$$

$$- 1115589354998164768436987302067428643851591482524199035357 5343000\, n^{42}$$

$$+ 168883131746867090216521890580629259354361771394676873045128 38280\, n^{41}$$

$$+ 557878756394989885053602688131932866101718507947981441845842287080\ n^{40}$$

$$- 845262291179828182091230126726930762120295850491706006421019849760\ n^{39}$$

$$- 25435755082981609122815310164558052465087113005768832786262436915160\ n^{38}$$

$$+ 3857626377006232775268580310200544078690817433899102182604165297620\ n^{37}$$

$$+ 1052580433013360701980188407782238814094779989080202915088963631337150\ n^{36}$$

$$- 1598158781405072216857916901828458493181515392337253923724747529654535\ n^{35}$$

$$- 39337171322984292289094342954937793408984455184843277472007281758462285\ n^{34}$$

$$+ 59804836375178974542070472883320919360067440473433543169873296402520695\ n^{33}$$

$$+ 132032138628536091548029752231100160432102107836499821173664213478417875\ n^{32}$$

$$- 2010384497615630860491481519908310700328186881991466503345432968418887160\ n^{31}$$

$$- 39554279574678738614570499280395760701136988044850238713779690813620673160\ n^{30} +$$

$$603366116108259233521014896805477964018695755082710913223422525704640453320\ n^{29} +$$

$$1050305973689451338861178588767993476558956109218647601432920746720761445600\ n^{28} -$$

$$16056272663395899696967818627992264113039368951582106947810552246433462395060\ n^{27} - 24525100568716455451806935541642518785447002100363040463085417651285103\backslash 129660\ n^{26} + 375904644862444781626943126264599102346901872633561416853340260\backslash 0144385892020\ n^{25} + 49904987689609483547244911789065610829058703507118824016\backslash 870202158519341741300\ n^{24} - 7673700475872644922900208331492141175536331490738\backslash 46043109572004537851205557960\ n^{23} - 87571007222047609332298460774603358971533\backslash 59277177466974099643259639511075162472\ n^{22} + 13519336107100773645989779532765\backslash 110904506855490303123482704250891728192215522688\ n^{21} + 1308931327914933776090\backslash 03673267203929191104529686672423952883829307535429461796204\ n^{20} - 20309936724\backslash 079045323650039966718844923891022275160197670677869407167240300455650\ n^{19} - 164226995851487283542766998500543269496381288691953607817002902829643682305\backslash 3377462\ n^{18} + 25649546213927044797597517734174326706517444151688421609038247\backslash 7148238854730294018\ n^{17} + 169901696130550155806654348759300326917040747209379\backslash 63608850262541247887754166556567954\ n^{16} - 26767731730278875610878029902565920671\backslash 088699302165387521320585050445951058614982440\ n^{15} - 1417453604609033343311188\backslash 07378169744949291753257944278870526869938369071615616624792\ n^{14} + 22600190655\backslash 6494439302117226018537577759481979537999112066450597432776582952732428408\ n^{13} + 926581067351518340630773186608640926955750790920002040465288561579767668\backslash$$

$$54661513584\,n^{12} - 150287255430552473059721839292223017931336717614900261670\backslash$$
$$30185829535334417583584584580\,n^{11} - 4564083478967202039066647487232460464771877\backslash$$
$$782326063485808328055216762919003399351372\,n^{10} + 7597561495603565423898580427\backslash$$
$$3098057868145002615635965370640013743018210993842782693\,48\,n^{9} + 159980504294027\backslash$$
$$04337162462254908626404532209868182189517372899009981576854522644416164\,n^{8} - $$
$$2779585639190583921769298359601784250020556493305508254459134902123275831\,4\backslash$$
$$76105758920\,n^{7} - 3627660144351306799585865854543850283881488071409693508805\,18\backslash$$
$$0161465467428388336021155\,2\,n^{6} + 6831283036122252160263447961616667550832510\,35\backslash$$
$$37672943904373377023043649341563093196788\,n^{5} + 4328151496210539012347264625\,87\backslash$$
$$49373683427907246385772666502058912045568040474969991594\,n^{4} - 990786876237693\backslash$$
$$459865262091962073982793044126384151309519397768795017673149400158578\,5\,n^{3} - $$
$$878610350812886405164845568464709230796651841147720878697655209577557786301\backslash$$
$$9182959667\,n^{2} + 6271849907407796907073578812507433760160198393642337865643471\backslash$$
$$658345845516027577523239\,3\,n - 20906166358025989690245262708358112533867327978\backslash$$
$$80779288547823886115281838675859174413\,1\,)\,n\,/\,7365930$$